\documentclass[10pt]{amsart}

\usepackage[leqno]{amsmath}
\usepackage{amsthm}
\usepackage{amsfonts}
\usepackage{amssymb}
\usepackage{eucal}
\usepackage[all]{xy}
\CompileMatrices

\DeclareFontFamily{OT1}{rsfs}{}
\DeclareFontShape{OT1}{rsfs}{n}{it}{<-> rsfs10}{}
\DeclareMathAlphabet{\mathscr}{OT1}{rsfs}{n}{it}

\DeclareMathOperator{\Hom}{Hom}

\DeclareMathOperator{\Aut}{Aut}
\DeclareMathOperator{\an}{an}

\DeclareMathOperator{\Spec}{Spec}

\newcommand*{\Sp}{{{\rm{Sp}}}}

\DeclareMathOperator{\Spf}{Spf}

\DeclareMathOperator{\et}{\acute{e}t}

\def\bir{{\rm bir}}

\def\wt{\widetilde}

\def\oS{{\overline{S}}}
\def\ox{{\overline{x}}}
\def\oR{{\overline{R}}}
\def\oU{{\overline{U}}}
\def\op{{\overline{p}}}

\def\tilU{{\wt U}}

\def\bfP{{\mathbf P}}
\def\calO{{\mathcal O}}

\newcommand*{\R}{\ensuremath{\mathbf{R}}}                        
\newcommand*{\Z}{\ensuremath{\mathbf{Z}}}                        
\newcommand*{\Q}{\ensuremath{\mathbf{Q}}}                        
\newcommand*{\C}{\ensuremath{\mathbf{C}}}                        

\renewcommand*{\P}{\ensuremath{\mathbf{P}}}                        
\renewcommand*{\O}{\mathscr{O}}                                  
\newcommand*{\calA}{\mathscr{A}}

\newcommand*{\calH}{\mathscr{H}}                               

\newcommand*{\calX}{\mathscr{X}}
\newcommand*{\calU}{\mathscr{U}}

\theoremstyle{plain}
  \newtheorem{theorem}[subsubsection]{Theorem}
  
  \newtheorem{lemma}[subsubsection]{Lemma}
  \newtheorem{corollary}[subsubsection]{Corollary}

\theoremstyle{definition}
  \newtheorem{definition}[subsubsection]{Definition}

\theoremstyle{remark}
  \newtheorem{example}[subsubsection]{Example}
  \newtheorem{remark}[subsubsection]{Remark}

\numberwithin{equation}{subsection}

   \DeclareMathOperator{\calR}{\mathscr R}
   
   \DeclareMathOperator{\calM}{\mathscr M}
   \DeclareMathOperator{\calP}{\mathscr P}

\title{Non-archimedean analytification of algebraic spaces}
\author{Brian Conrad}
\address{Department of Mathematics\\
University of Michigan\\
Ann Arbor, MI 48109, USA}
\email{bdconrad@umich.edu}
\author{Michael Temkin}
\address{Department of Mathematics\\
University of Pennsylvania\\
Philadelphia, PA  19104 USA}
\email{temkin@math.upenn.edu}
\date{June 20, 2007}
\thanks{The work of B.C. was partially supported by NSF grant DMS-0600919, and both
authors are grateful to the participants of the Arizona Winter School for helpful feedback on
an earlier version of this paper.  Also, B.C. is very grateful to Columbia University for its
generous hospitality during a sabbatical visit.}
\begin{document}

\maketitle

\section{Introduction}\label{motive}

\subsection{Motivation}
This paper is largely concerned with constructing
quotients by \'etale equivalence relations.  We are inspired by
questions in classical rigid geometry, but to give satisfactory answers in that category
we have to first solve quotient problems within the framework of Berkovich's $k$-analytic spaces.
One source of motivation is the relationship between
algebraic spaces and analytic spaces over $\C$, as follows.
If $X$ is a reduced and irreducible proper
complex-analytic space then the
meromorphic functions on $X$ form a field
$\mathscr{M}(X)$ and this field is finitely generated
over $\C$ with transcendence degree at most $\dim X$ \cite[8.1.3,~9.1.2,~10.6.7]{grauert}.
A proper complex-analytic space $X$ is called
{\em Moishezon} if ${\rm{trdeg}}_{\C}(\mathscr{M}(X_i)) = \dim X_i$
for all irreducible components $X_i$ of $X$ (endowed with the reduced structure).
Examples of such spaces are analytifications of proper
$\C$-schemes, but Moishezon found more examples, and Artin
found ``all'' examples by analytifying
proper algebraic spaces over $\C$.
To be precise, the analytification $\mathscr{X}^{\rm{an}}$ of an algebraic space $\mathscr{X}$
locally of finite type over $\C$ \cite[Ch.~I,~5.17{\em{ff}}]{knutson} is the unique solution to an
\'etale quotient problem
that admits a solution if and only if  $\mathscr{X}$ is locally separated over $\C$
(in the sense that $\Delta_{\mathscr{X}/\C}$ is an immersion).
The functor $\mathscr{X} \rightsquigarrow \mathscr{X}^{{\rm{an}}}$ from the category of proper
algebraic spaces over $\C$ to the category of Moishezon spaces is fully faithful,
and it is a beautiful
theorem of Artin \cite[Thm.~7.3]{artin2} that this is an equivalence of categories.

It is natural to ask if a similar theory works over a non-archimedean
base field $k$ (i.e., a field $k$ that is complete with
respect to a fixed nontrivial non-archimedean absolute value).
This is a surprisingly nontrivial question.
One can carry over the definition of analytification of
locally finite type algebraic spaces $\mathscr{X}$ over $k$ in terms of uniquely solving
a rigid-analytic \'etale quotient problem; when the quotient
problem has a solution we say that $\mathscr{X}$ is {\em analytifiable}
in the sense of rigid geometry.
(See  \S\ref{algspacesec}
for a general discussion of definitions, elementary results,
and functorial properties of analytification over such $k$.)
In Theorem \ref{locsep} we show that
local separatedness is a necessary condition for analytifiability
over non-archimedean base fields,
but in contrast
with the complex-analytic case
it is not sufficient;
there are smooth 2-dimensional counterexamples
over any $k$ (even arising from algebraic
spaces over the prime field), as we shall explain in Example \ref{counter}.
In concrete terms, the surprising dichotomy
between the archimedean and non-archimedean worlds is due to
the lack of a Gelfand--Mazur theorem over non-archimedean fields.
(That is, any non-archimedean field $k$ admits nontrivial
non-archimedean extension fields with a compatible absolute value,
even if $k$ is algebraically closed.)
Since local separatedness fails to be a sufficient criterion for the
existence of non-archimedean analytification of an algebraic space,
it is natural to seek a reasonable salvage of the situation.  We view
separatedness as a reasonable additional hypothesis to impose
on the algebraic space.

\subsection{Results}
Our first main result is the following (recorded as Theorem \ref{makespace}):
\begin{theorem}\label{11} Any separated algebraic
space locally of finite type over a non-archimedean field is analytifiable in the sense of rigid geometry.
\end{theorem}
For technical reasons related to admissibility of coverings
and the examples of non-analytifiable algebraic spaces
in Example \ref{counter}, we do not think it is possible to prove this
theorem via the methods of rigid geometry.
The key to our success is to study quotient problems in the category of
$k$-analytic spaces in the sense of Berkovich, and only {\em a posteriori}
translating such results back into classical rigid geometry.
In this sense, the key theorem in this paper is an existence
result for quotients in the $k$-analytic category (recorded as Theorem \ref{newb}):

\begin{theorem}\label{berqt}
Let $R \rightrightarrows U$ be an \'etale equivalence relation
in the category of $k$-analytic spaces.  The quotient $U/R$ exists whenever
the diagonal $\delta:R \rightarrow U \times U$ is a closed immersion.
 In such cases $U/R$ is separated.  If $U$ is good $($resp. strictly $k$-analytic$)$ then so is $U/R$.
\end{theorem}
The quotient $X = U/R$ in this theorem represents a suitable quotient sheaf for the \'etale topology, and the
natural map $U \rightarrow X$ is an \'etale surjection with respect to which the natural map $R \rightarrow U
\times_X U$ is an isomorphism. (See \S\ref{prelimcon} for a discussion of basic notions related to quotients in
the $k$-analytic category.)  The hypothesis on $\delta$ in Theorem \ref{berqt} is never satisfied by the \'etale
equivalence relations arising from $k$-analytification of \'etale charts $\mathscr{R} \rightrightarrows
\mathscr{U}$ of {\em non-separated} locally separated algebraic spaces $\mathscr{X}$ locally of finite type over
$k$.

The locally separated
algebraic spaces in Example \ref{counter} that are not analytifiable in the
sense of rigid geometry are also not analytifiable in the sense of
$k$-analytic spaces (by Remark \ref{kcounter}).  Thus, some
assumption on $\delta$ is necessary in Theorem \ref{berqt} to avoid such examples.
We view quasi-compactness assumptions on $\delta$ as a reasonable way to proceed.
Note that since $\delta$ is separated (it is a monomorphism), it is quasi-compact
if and only if it topologically proper (i.e., a {\em compact} map in the terminology of Berkovich).

\begin{remark}\label{qcrem} There are two senses in which it is not possible to generalize
Theorem \ref{berqt}
by just assuming $\delta$ is compact (i.e., topologically proper).
First, if we assume that $U$ is locally separated
in the sense that each $u \in U$ has a separated open neighborhood
then compactness of $\delta$ forces it to be a closed immersion.
That is, if an \'etale equivalence
relation on a locally separated
$k$-analytic space $U$ has a compact diagonal then this diagonal is automatically
a closed immersion.  For example, a locally separated
$k$-analytic space $X$ has compact diagonal $\Delta_X$ if
and only if it is separated.  This is shown in \cite[2.2]{descent}.

Second, with no local separatedness assumptions on $U$ it does
not suffice to assume that $\delta$ is compact.  In fact, in Example \ref{ex54}
we give examples of (non-separated) compact Hausdorff $U$
and a free right action on $U$ by a finite group $G$ such that $U/G$ does not
exist.  In such cases the action map $\delta:R = U \times G \rightarrow U \times U$
defined by $(u,g) \mapsto (u, u.g)$
is a compact map.    We also give analogous such examples of non-existence
in the context of rigid geometry.
Recall that in the $k$-analytic category, separatedness is stronger than the Hausdorff
property since the theory of fiber products is finer than in topology.

The weaker Hausdorff property of $U/R$ in Theorem \ref{berqt}
can also be deduced from the general fact (left as an exercise)
that a locally Hausdorff and locally compact topological space $X$
is Hausdorff if and only if it is quasi-separated (in the sense that
the overlap of any pair of quasi-compact subsets of $X$ is quasi-compact).
\end{remark}

Beware that \'etaleness in rigid geometry is a weaker
condition than in the category of $k$-analytic spaces (aside from maps
arising from algebraic scheme morphisms), so our existence result for
$U/R$ in the $k$-analytic category does not yield a corresponding
existence result for  \'etale equivalence
relations in rigid geometry when the diagonal map
is a closed immersion: we do not have a satisfactory analogue of Theorem \ref{berqt}
in rigid geometry.

\begin{example} Let $U$ be a separated $k$-analytic space
and let $G$ be an abstract group equipped with a right action
on $U$ that is free in the sense that the action map $\delta:U \times G \rightarrow
U \times U$ is a monomorphism.  Equivalently, $G$ acts freely
on $U(K)$ for each algebraically closed analytic extension $K/k$. Assume that
the action is properly discontinuous with finite physical stabilizers in the sense
that for each $u \in U$ the group $G_u = \{g \in G\,|\,g(u) = u\}$
is finite and $u$ admits an open neighborhood $V$ such that
$g(V) \cap V = \emptyset$ for all $g \in G - G_u$.   In such cases
$\delta$ is a closed immersion, so by Theorem
\ref{berqt} the quotient $U/G$ exists as a separated $k$-analytic space;
it is good (resp. strictly $k$-analytic) if $U$ is.
If $U$ is paracompact and strictly $k$-analytic and
$U_0$ denotes the associated separated
rigid space in the sense of \cite[1.6.1]{berihes} then
the separated strictly $k$-analytic quotient $U/G$ is paracompact
if either $G$ is finite or $U$ is covered by countably many compacts,
and when $U/G$ is paracompact then the
associated separated rigid space $(U/G)_0$ is a quotient $U_0/G$ in
the sense of rigid geometry (as defined in \S\ref{topqt}).  We leave the details to the reader.
As a special case, if $k'/k$ is a finite Galois extension
then descent data relative to $k'/k$ on separated $k'$-analytic spaces is always effective.
(This allows one to replace $G$ with an \'etale $k$-group above, provided
that its geometric fiber is finitely generated.)
\end{example}

\subsection{Strategy of proof}
To construct $U/R$ under our assumptions in Theorem \ref{berqt}, we first observe
that the problem can be localized because representable functors on the category of
$k$-analytic spaces are sheaves for the \'etale topology (proved in
\cite{berihes} under some goodness/strictness hypotheses that we remove).
We are thereby able to use an \'etale localization argument
and the nice topological properties
of $k$-analytic spaces (via
compactness and connectedness arguments) to reduce to the case of
a free right action by a
finite group $G$ on a possibly non-affinoid space $U$.
Since the action map $\delta:U \times G \rightarrow U \times U$ is a closed
immersion (due to the initial hypotheses), such $U$ must at least be separated.

Our existence result for $U/G$
requires a
technique to construct coverings
of $U$ by affinoid domains that have good behavior with respect to the
$G$-action.   We do not know any technique of this type in  rigid geometry
(with admissibility of the covering),
but in the $k$-analytic category there is a powerful tool to do this, namely
the theory of reduction of germs of $k$-analytic spaces that was developed in
\cite{temkin2}.  Briefly, this technique often reduces difficult local construction problems
in non-archimedean geometry to more tractable problems
in an algebro-geometric setting (with Zariski--Riemann spaces).
Here we use the separatedness of $U$:
assuming only that $U$ is Hausdorff (and making no assumptions
on $\delta$) we formulate a
condition on reduction of germs that is always satisfied when $U$ is separated and is sufficient for
the existence of $U/G$.

A special case of Theorem \ref{berqt} is the $k$-analytification of
an \'etale chart of a separated algebraic space locally of finite type over $k$.
This yields an analogue of Theorem \ref{11} using analytification in
the sense of $k$-analytic spaces.  But
in Theorem \ref{berqt} if
$U$
is strictly $k$-analytic (resp. good) then the quotient space $U/R$ is as well.
Hence, we can make a link with quotient constructions in classical rigid geometry,
thereby deducing Theorem
\ref{11} from Theorem \ref{berqt}.

\subsection{Overview of paper}   In \S\ref{algspacesec} we gather some basic definitions,
preliminary results, and formalism related to flat equivalence relations
and analytification in the sense of rigid geometry.
Examples of non-analytifiable (but locally separated) algebraic
spaces are given in \S\ref{countersec}, where we also discuss some
elementary examples of \'etale quotients in the affinoid case
and we review GAGA for proper algebraic spaces over
a non-archimedean field.   This is all  preparation for
\S\ref{analberk}, where we adapt \S\ref{algspacesec}
to the $k$-analytic category and then carry out the preceding strategy to prove
the two theorems stated above.  More precisely, in \S\ref{analberk} we
reduce ourselves to proving the existence of $U/G$ when
$G$ is a finite group with a free right action on
a separated $k$-analytic space $U$, and in \S\ref{birsec}
we solve this existence problem by using the theory from \cite{temkin2}.

Since we can now always analytify
proper algebraic spaces, as an application of our results it makes sense to try to
establish a rigid-analytic analogue of Artin's equivalence between
proper complex-analytic spaces and Moishezon spaces.
This involves new difficulties
(especially in positive characteristic),
and so will be given in a subsequent paper \cite{moish}.

\subsection{Conventions}  The ground field $k$ in the rigid-analytic
setting is always understood to be a non-archimedean field,
and a field extension $K/k$ is called {\em analytic} if $K$ is non-archimedean
with respect to a fixed absolute value extending the one on $k$.
When we work with Berkovich's theory of analytic spaces we also allow
the possibility that the absolute value on the ground field is trivial,
and we define the notion of {\em analytic} extension similarly (allowing trivial
absolute values).  By abuse of notation, if $x$ is a point in a $k$-analytic space $X$ then
we write $(X,x)$ to denote the associated germ (denoted $X_x$ in \cite{temkin1},
\cite{temkin2}).

We require algebraic spaces to have quasi-compact diagonal over $\Spec \Z$, and an
{\em \'etale chart} $\mathscr{R} \rightrightarrows \mathscr{U}$ for
an algebraic space $\mathscr{X}$ is always understood to have
$\mathscr{U}$ (and hence $\mathscr{R} = \mathscr{U} \times_{\mathscr{X}} \mathscr{U}$) be a scheme.  Also,
throughout this paper, it is tacitly understood that
``algebraic space'' means ``algebraic space
locally of finite type over $k$'' (and maps between them are $k$-maps) unless we explicitly say
otherwise (e.g., sometimes we work over $\C$).

\tableofcontents

\section{\'Etale equivalence relations and algebraic spaces}
\label{algspacesec}

In this section we develop basic concepts related to analytification
for algebraic spaces $\mathscr{X}$ that are
locally of finite type over $k$.
An existence theorem for such analytifications
in the separated case is given in
\S\ref{mainsec}.

\subsection{Topologies and quotients}\label{topqt}

We give the category of rigid
spaces (over $k$) the {\em Tate-fpqc topology} that is generated by the Tate topology
and the class of faithfully flat maps that
admit local {\em fpqc} quasi-sections
in the sense of \cite[Def.~4.2.1]{relamp}; that is, a covering of $X$ is a collection of flat
maps $X_i \rightarrow X$ such that
locally on $X$ for the Tate topology there
exist sections to $\coprod X_i \rightarrow X$ after a faithfully flat and quasi-compact base change.
If we instead work with \'etale surjections and require
the existence of local \'etale quasi-sections in the sense
of \cite[Def.~4.2.1]{relamp} then we
get the {\em Tate-\'etale topology}.
By \cite[Cor.~4.2.5]{relamp}, all representable functors are sheaves
for the Tate-fpqc topology.

\begin{example}  If $\{X_i\}$ is a set of admissible opens
in $X$ that is a set-theoretic cover then it is a cover of $X$ for
the Tate-fpqc topology if and only if it is an admissible covering
for the usual Tate topology.
\end{example}

\begin{example}\label{coverex}
Let $h:\mathscr{X}' \rightarrow \mathscr{X}$
be a faithfully flat map of algebraic $k$-schemes.  By \cite[Thm.~4.2.2]{relamp}
(whose proof uses $k$-analytic spaces),
the associated map $h^{\rm{an}}:X' \rightarrow X$ between analytifications
is a covering map for the Tate-fpqc topology.  Likewise,
if $h$ is an \'etale surjection then $h^{\rm{an}}$ is a covering
map for the Tate-\'etale topology.
\end{example}

\begin{example}\label{topex}
We will later need GAGA for proper algebraic spaces, but algebraic
spaces have only an \'etale topology rather than a Zariski topology.
Thus, to make the comparison of coherent cohomologies,
on both the algebraic and rigid-analytic sides we want to use
an \'etale topology.  In particular, as preparation for this we need
to compare coherent
cohomology for the Tate topology and for the Tate-\'etale topology
on a rigid space.  We now explain how this goes.

For any rigid space $S$, let $S_{\rm{Tate}}$ denote the site defined by the Tate topology
(objects are admissible opens in $S$) and
let $S_{\et}$ denote the site defined by the Tate-\'etale topology
(objects are rigid spaces \'etale over $S$).  The evident left-exact pushforward
map from sheaves on sets on $S_{\et}$ to sheaves of sets on $S_{\rm{Tate}}$
has an exact left adjoint that is constructed by sheafification in  the usual
manner, so the continuous map of sites $S_{\et} \rightarrow S_{\rm{Tate}}$
may be uniquely enhanced to a map of topoi
$\widetilde{S_{\et}} \rightarrow \widetilde{S_{\rm{Tate}}}$.  Descent theory for coherent
sheaves on rigid spaces (see \cite[Thm.~4.2.8]{relamp}
for the formulation we need) shows that $\mathscr{O}_{S_{{\et}}}(U \rightarrow S) :=
\mathscr{O}_U(U)$ is a sheaf
on $S_{\et}$.  Using this sheaf of $k$-algebras
makes $\widetilde{S_{\et}} \rightarrow \widetilde{S_{\rm{Tate}}}$ a map of ringed topoi in the evident manner.
Let $\mathscr{F} \rightsquigarrow \mathscr{F}_{\et}$ denote the associated pullback operation
on sheaves of modules.   This is an exact functor because
for any \'etale map $h:U \rightarrow V$ between affinoid spaces
the image $h(U)$ is covered by finitely many admissible
affinoid opens $V_i \subseteq V$ and the pullback map of affinoids $h^{-1}(V_i) \rightarrow V_i$
is flat on coordinate rings.   It is then easy to check via descent theory (as for
schemes with the Zariski and \'etale topologies) that
$\mathscr{F} \rightsquigarrow \mathscr{F}_{\et}$
establishes an equivalence between
the full subcategories of coherent
sheaves on $S_{\et}$ and $S_{\rm{Tate}}$.

Exactness of pullback provides a canonical
$\delta$-functorial $\mathscr{O}_{Y_{\et}}$-linear comparison morphism
$$({\rm{R}}^j h_{\ast}(\mathscr{F}))_{\et} \rightarrow {\rm{R}}^j(h_{\et})_{\ast}(
\mathscr{F}_{\et})$$
for any map of rigid spaces $h:X \rightarrow Y$
and any $\mathscr{O}_X$-module $\mathscr{F}$, and we claim that this is an
isomorphism when $h$ is proper and $\mathscr{F}$ is coherent.
By coherence of ${\rm{R}}^j h_{\ast}(\mathscr{F})$ on
$Y$ for such $h$ and $\mathscr{F}$, this immediately reduces to the general claim that
for any rigid space $S$ the canonical $\delta$-functorial comparison map
${\rm{H}}^i(S,\mathscr{F})
\rightarrow {\rm{H}}^i(S_{\et},\mathscr{F}_{\et})$ for $\mathscr{O}_S$-modules
$\mathscr{F}$ is an isomorphism when $\mathscr{F}$ is coherent.  Using
a \v{C}ech-theoretic
spectral sequence exactly as in the comparison of \'etale and Zariski cohomology
for quasi-coherent sheaves on schemes, this reduces the problem  to checking
that if $\Sp(B) \rightarrow \Sp(A)$ is an \'etale surjection of affinoids
and $M$ is a finite $A$-module then the habitual complex
$$0 \rightarrow M \rightarrow M \widehat{\otimes}_A B \rightarrow
M \widehat{\otimes}_A (B \widehat{\otimes}_A B) \rightarrow \dots$$
is an exact sequence. By the theory of formal models over the valuation ring $R$ of $k$,
it suffices to consider the case when there is
a faithfully flat map $\Spf(\mathscr{B}) \rightarrow \Spf(\mathscr{A})$
of admissible formal $R$-schemes
and a finitely presented $\mathscr{A}$-module $\mathfrak{M}$
giving rise to $\Sp(B) \rightarrow \Sp(A)$ and $M$ on generic fibers.
Thus, it is enough to show that the complex
$$0 \rightarrow \mathfrak{M} \rightarrow \mathfrak{M} \widehat{\otimes}_{\mathscr{A}}
\mathscr{B} \rightarrow \mathfrak{M} \widehat{\otimes}_{\mathscr{A}}
(\mathscr{B} \widehat{\otimes}_{\mathscr{A}} \mathscr{B}) \rightarrow \dots$$
is exact.   For any fixed $\pi \in k$ with $0 < |\pi| < 1$, this is a complex
of $\pi$-adically separated and complete $R$-modules, so
by a simple diagram chase with compatible liftings it suffices to prove
exactness modulo $\pi^n$ for all $n \ge 1$.  This in turn is a special case of ordinary
faithfully flat descent theory since $\mathscr{A}/(\pi^n) \rightarrow \mathscr{B}/(\pi^n)$
is faithfully flat for all $n$.
\end{example}

Let $X' \rightarrow X$ be a flat surjection of rigid spaces
and assume that it admits local {\em fpqc} quasi-sections.
The maps $R = X' \times_X X' \rightrightarrows X'$
define a monomorphism $R \rightarrow X' \times X'$, and
we have an isomorphism $X'/R \simeq X$ as sheaves of sets
on the category of rigid spaces with the Tate-fpqc topology since the maps
$R \rightrightarrows X'$ are faithfully flat
and admit local {\em fpqc} quasi-sections in such cases.

Conversely,
given a pair of flat maps $R \rightrightarrows X'$
admitting local {\em fpqc} quasi-sections
such that $R \rightarrow X' \times X'$ is functorially an equivalence
relation (in which case we call
$R \rightarrow X' \times X'$ a {\em flat equivalence relation}),
consider the sheafification of
the presheaf $Z \mapsto X'(Z)/R(Z)$ with respect to the Tate-fpqc topology.  If this
sheaf is represented by some rigid space $X$ then we
call $X$ (equipped with the map $X' \rightarrow X$)
the {\em flat quotient} of $X'$ modulo $R$ and we denote it $X'/R$.
By the very definition of the Tate-fpqc topology  that is used to define
the quotient sheaf $X'/R$, if a flat quotient $X$ exists
then the projection map $p:X' \rightarrow X$ admits local {\em fpqc} quasi-sections.
Moreover, $p$ is automatically faithfully flat.  Indeed, arguing as in the case
of schemes, choose a faithfully
flat map $z:Z \rightarrow X$
such that there is a quasi-section $z':Z \rightarrow X'$ over $X$.
The map $p$ is faithfully flat if and only if the projection $q_2:X' \times_X Z \rightarrow Z$
is faithfully flat, and via the isomorphism $X' \times_X Z \simeq X' \times_X X' \times_{X',z'} Z =
R \times_{p_2,X',z'} Z$ the map $q_2$ is identified with a base change of
the projection $p_2:R \rightarrow X'$ that is faithfully flat.

When the flat quotient $X = X'/R$ exists, the map
$R \rightarrow X' \times_X X'$ is an isomorphism and
so for every property $\P$ in \cite[Thm.~4.2.7]{relamp} the map
$X' \rightarrow X$ satisfies $\P$ if and only if the maps
$R \rightrightarrows X'$ satisfy $\P$.
Likewise, $X$ is quasi-separated (resp. separated)
if and only if the map $R \rightarrow X' \times X'$ is quasi-compact (resp.
a closed immersion).
By descent theory for
morphisms, the diagram of sets
\begin{equation}\label{zxmap}
\Hom(X,Z) \rightarrow \Hom(X',Z) \rightrightarrows \Hom(R,Z)
\end{equation}
is left-exact for any rigid space $Z$ when $X = X'/R$ is
a flat quotient.

\begin{definition} An {\em \'etale equivalence relation}
on a rigid space $X'$ is a functorial equivalence relation
$R \rightarrow X' \times X'$ such that the maps
$R \rightrightarrows X'$ are \'etale
and admit local \'etale quasi-sections
in the sense of \cite[Def.~4.2.1]{relamp}.
If the flat quotient $X'/R$ exists, it is called
an {\em \'etale quotient} in such cases.
\end{definition}

\begin{example}\label{etqtex}
If $X' \rightarrow X$ is an \'etale surjection
that admits local \'etale quasi-sections then the
\'etale quotient of $X'$ modulo the \'etale equivalence relation
$R = X' \times_X X'$ exists: it is $X$.   Thus, the preceding
arguments with flat quotients work with ``faithfully flat'' replaced by ``\'etale surjective''
to show that in the definition of \'etale quotient it does not matter if
we form $X'/R$ only with respect to the Tate-\'etale topology.
\end{example}

\begin{lemma}\label{rxx} Let $R \rightarrow X' \times X'$ be a flat
equivalence relation on a rigid space $X'$, and assume
that the flat quotient $X'/R$ exists.
The equivalence relation $R \rightarrow X' \times X'$ is \'etale
if and only if
the map $X' \rightarrow X'/R$ is \'etale and admits local
\'etale quasi-sections.
\end{lemma}

\begin{proof}
Let $X = X'/R$.  Since $R = X' \times_X X'$
and the map $X' \rightarrow X$ is faithfully flat
with local {\em fpqc} quasi-sections, we may use
\cite[Thm.~4.2.7]{relamp} for the property $\P$
of being \'etale with local \'etale
quasi-sections.
\end{proof}

 Let $\mathscr{X}$ be an algebraic
space and let $\mathscr{R} \rightrightarrows
\mathscr{U}$ be an \'etale chart for $\mathscr{X}$.
By Example \ref{coverex}, the maps
$\mathscr{R}^{\an} \rightrightarrows \mathscr{U}^{\an}$
admit local \'etale quasi-sections.  Since
a map in any category with fiber products
is a monomorphism if and only if its
relative diagonal is an isomorphism, analytification
of algebraic $k$-schemes
carries monomorphisms to monomorphisms.  Thus, the morphism
$\mathscr{R}^{\an} \rightarrow \mathscr{U}^{\an}
\times \mathscr{U}^{\an}$ is a monomorphism and so
$\mathscr{R}^{\an}$ is an \'etale equivalence relation on
$\mathscr{U}^{\an}$.  It therefore
makes sense to ask if the \'etale quotient
$\mathscr{U}^{\an}/\mathscr{R}^{\an}$ exists.

We will show in Lemma \ref{analgspace} that such existence and the actual quotient
$\mathscr{U}^{\rm{an}}/\mathscr{R}^{\rm{an}}$ (when it exists!)
are independent of the chart $\mathscr{R} \rightrightarrows \mathscr{U}$
for $\mathscr{X}$ in a canonical manner,
in which case we define $\mathscr{U}^{\rm{an}}/\mathscr{R}^{\rm{an}}$
to be the {\em analytification} of $\mathscr{X}$.
The rigid-analytic \'etale equivalence relations
$\mathscr{R}^{\rm{an}} \rightrightarrows \mathscr{U}^{\rm{an}}$ that arise
in the problem of analytifying algebraic spaces
are rather special, and so one might hope that in such cases
the required quotient exists whenever
the algebraic space $\mathscr{X}$ is locally separated over $k$ (as
is necessary and sufficient for the existence of analytifications of algebraic spaces
over $\C$ in the complex-analytic sense).  However,
we will give counterexamples in Example \ref{counter}:  locally separated algebraic
spaces that are not analytifiable in the above sense defined via quotients.   In the positive direction,
the quotient $\mathscr{U}^{\rm{an}}/\mathscr{R}^{\rm{an}}$
will be shown to always exist when $\mathscr{X}$ is separated.

\subsection{Analytification of algebraic spaces}\label{analgsec}

Let $\mathscr{X}$ be an algebraic space, and $\mathscr{R} \rightrightarrows \mathscr{U}$
an \'etale chart for $\mathscr{X}$.
We now address the ``independence of choice'' and canonicity issues
for $\mathscr{U}^{\rm{an}}/\mathscr{R}^{\rm{an}}$ in terms of $\mathscr{X}$.
These will go essentially as in the complex-analytic case except
that we have to occasionally use properties related
to local \'etale quasi-sections for the Tate topology.
In the complex-analytic case it
does not seem that the relevant arguments are available in the literature,
so for this reason and to ensure that the Tate topology
presents no difficulties we have decided to give the arguments in detail
(especially so we can see that they carry over to $k$-analytic spaces, as we shall
need later).

Let $\mathscr{R}_1 \rightrightarrows \mathscr{U}_1$
and $\mathscr{R}_2 \rightrightarrows \mathscr{U}_2$ be
two \'etale charts for $\mathscr{X}$.  Let
$\mathscr{U}_{12} = \mathscr{U}_1\times_{\mathscr{X}} \mathscr{U}_2$
and let
$\mathscr{R}_{12} = \mathscr{R}_1 \times_{\mathscr{X}}
\mathscr{R}_2$, so
$\mathscr{R}_{12} \rightrightarrows \mathscr{U}_{12}$ is an \'etale
chart dominating each chart $\mathscr{R}_i \rightrightarrows
\mathscr{U}_i$.

\begin{lemma}\label{analgspace}
If $\mathscr{U}_1^{\an}/\mathscr{R}_1^{\an}$ exists
then so do $\mathscr{U}_2^{\an}/\mathscr{R}_2^{\an}$
and $\mathscr{U}_{12}^{\an}/\mathscr{R}_{12}^{\an}$,
and the natural maps $$\pi_i:\mathscr{U}_{12}^{\an}/\mathscr{R}_{12}^{\an}
 \rightarrow
\mathscr{U}_i^{\an}/\mathscr{R}_i^{\an}$$ are isomorphisms.
The induced isomorphism
$\phi = \pi_2 \circ \pi_1^{-1}:
\mathscr{U}_1^{\an}/\mathscr{R}_1^{\an} \simeq
\mathscr{U}_2^{\an}/\mathscr{R}_2^{\an}$
is transitive with respect to a third choice of
\'etale chart for $\mathscr{X}$.
\end{lemma}

\begin{proof}
The natural composite map
$\mathscr{U}_{12}^{\an} \rightarrow
\mathscr{U}_1^{\an} \rightarrow \mathscr{U}_1^{\an}/\mathscr{R}_1^{\an}$
is \'etale with local \'etale quasi-sections (as each
step in the composite has this property, due to
Example \ref{coverex} applied to $\mathscr{U}_{12} \rightarrow
\mathscr{U}_1$ and the defining properties
for $\mathscr{U}_1^{\an}/\mathscr{R}_1^{\an}$
as an \'etale quotient).  We claim that
this composite map serves as the \'etale quotient for
$\mathscr{U}_{12}^{\an}$ by $\mathscr{R}_{12}^{\an}$
(so the \'etale quotient
$\mathscr{U}_{12}^{\an}/\mathscr{R}_{12}^{\an}$ exists
and $\pi_1$ is an isomorphism).  The problem is to prove
$$\mathscr{R}_{12}^{\an} \stackrel{?}{=}
\mathscr{U}_{12}^{\an} \times_{\mathscr{U}_1^{\an}/\mathscr{R}_1^{\an}}
\mathscr{U}_{12}^{\an}$$ as subfunctors of
$\mathscr{U}_{12}^{\an} \times \mathscr{U}_{12}^{\an}$,
and this is easily proved via two ingredients: the analytic isomorphism
$$\mathscr{R}_1^{\an} \simeq
\mathscr{U}_1^{\an} \times_{\mathscr{U}_1^{\an}/\mathscr{R}_1^{\an}}
\mathscr{U}_1^{\an}$$ (that arises from the assumption of existence
for the \'etale quotient $\mathscr{U}_1^{\an}/\mathscr{R}_1^{\an}$)
and the analytification of the algebraic isomorphism
$$\mathscr{R}_{12} \simeq (\mathscr{U}_{12} \times \mathscr{U}_{12})
\times_{\mathscr{U}_1 \times \mathscr{U}_1} \mathscr{R}_1$$
as algebraic $k$-schemes and as subfunctors of
$\mathscr{U}_{12} \times \mathscr{U}_{12}$.

Now we address the existence of $\mathscr{U}_2^{\an}/\mathscr{R}_2^{\an}$
and the isomorphism property for $\pi_2$.
The map $\mathscr{U}_{12}^{\an} \rightarrow
\mathscr{U}_2^{\an}$ is an \'etale surjection with
local \'etale quasi-sections, and so by rigid-analytic descent theory
the \'etale quotient map
$\mathscr{U}_{12}^{\an} \rightarrow
\mathscr{U}_{12}^{\an}/\mathscr{R}_{12}^{\an}$
admits at most one factorization
through the map $\mathscr{U}_{12}^{\an} \rightarrow
\mathscr{U}_2^{\an}$, in which case
the resulting map $h:\mathscr{U}_2^{\an}
\rightarrow \mathscr{U}_{12}^{\an}/\mathscr{R}_{12}^{\an}$
is an \'etale surjection with local \'etale
quasi-sections.  To prove that $h$
exists, by (\ref{zxmap})
 it is necessary and sufficient to check that the
two maps $\mathscr{U}_{12}^{\an} \times_{\mathscr{U}_2^{\an}}
\mathscr{U}_{12}^{\an} \rightrightarrows \mathscr{U}_{12}^{\an}/
\mathscr{R}_{12}^{\an}$ coincide, and this equality
follows from the inclusion
$$\mathscr{U}_{12} \times_{\mathscr{U}_2} \mathscr{U}_{12}
\subseteq \mathscr{U}_{12} \times_{\mathscr{X}} \mathscr{U}_{12}
= \mathscr{R}_{12}$$ as subfunctors of $\mathscr{U}_{12}
\times \mathscr{U}_{12}$.
Thus, we have an \'etale surjection
$h:\mathscr{U}_2^{\an} \rightarrow
\mathscr{U}_{12}^{\an}/\mathscr{R}_{12}^{\an}$
with local \'etale quasi-sections.

Since $\mathscr{R}_{12} = (\mathscr{U}_{12} \times \mathscr{U}_{12})
\times_{\mathscr{U}_2 \times \mathscr{U}_2} \mathscr{R}_2$
as subfunctors of $\mathscr{U}_{12} \times \mathscr{U}_{12}$,
the natural map $\mathscr{R}_{12} \rightarrow \mathscr{R}_2$
is an \'etale surjection and hence
the two composite maps
$$\mathscr{R}_2^{\an} \rightrightarrows \mathscr{U}_2^{\an}
\stackrel{h}{\rightarrow} \mathscr{U}_{12}^{\an}/\mathscr{R}_{12}^{\an}$$
are equal if and only if equality holds
after composition with the map
$\mathscr{R}_{12}^{\an} \rightarrow \mathscr{R}_2^{\an}$.
Such equality holds after composition because
$\mathscr{R}_{12} \rightrightarrows \mathscr{U}_{12}$
is co-commutative over $\mathscr{R}_2 \rightrightarrows
\mathscr{U}_2$.
Hence, we have
\begin{equation}\label{rinc}
\mathscr{R}_2^{\an} \subseteq
\mathscr{U}_2^{\an} \times_{\mathscr{U}_{12}^{\an}/\mathscr{R}_{12}^{\an}}
\mathscr{U}_2^{\an}
\end{equation}
as subfunctors of
$\mathscr{U}_2^{\an} \times \mathscr{U}_2^{\an}$.
To conclude that $\mathscr{U}_2^{\an} \rightarrow
\mathscr{U}_{12}^{\an}/\mathscr{R}_{12}^{\an}$
is an \'etale quotient by $\mathscr{R}_2^{\an}$, we need
the inclusion of subfunctors (\ref{rinc}) to be an equality, which is
to say that the natural map
$$g:
\mathscr{U}_2^{\an} \times_{\mathscr{U}_{12}^{\an}/\mathscr{R}_{12}^{\an}}
\mathscr{U}_2^{\an} \rightarrow \mathscr{U}_2^{\an} \times
\mathscr{U}_2^{\an}$$ factors through the subfunctor $\mathscr{R}_2^{\an}$.
Since the map
$$f:\mathscr{R}_{12}^{\an} \simeq \mathscr{U}_{12}^{\an}
\times_{\mathscr{U}_{12}^{\an}/\mathscr{R}_{12}^{\an}}
\mathscr{U}_{12}^{\an} \rightarrow \mathscr{U}_2^{\an} \times_{
\mathscr{U}_{12}^{\an}/\mathscr{R}_{12}^{\an}} \mathscr{U}_2^{\an}$$
is an \'etale surjection with local \'etale quasi-sections
(as $\mathscr{U}_{12} \rightarrow \mathscr{U}_2$ is an \'etale
surjection of schemes, so Example \ref{coverex} applies), by descent theory for
rigid-analytic morphisms
it suffices to check that the composite map $g \circ f$
factors through $\mathscr{R}_2^{\an}$.
Clearly $g \circ f$ equals the canonical composite map
$$\mathscr{R}_{12}^{\an} \rightarrow \mathscr{U}_{12}^{\an} \times
\mathscr{U}_{12}^{\an} \rightarrow \mathscr{U}_2^{\an} \times
\mathscr{U}_2^{\an}$$
that is the analytification of the composite map
$\mathscr{R}_{12} \rightarrow \mathscr{U}_{12} \times \mathscr{U}_{12}
\rightarrow \mathscr{U}_2 \times \mathscr{U}_2$.  This latter
composite map
obviously factors through the inclusion $\mathscr{R}_2 \rightarrow
\mathscr{U}_2 \times \mathscr{U}_2$, thereby completing
the proof that $\mathscr{U}_2^{\an}/\mathscr{R}_2^{\an}$
exists and that $\pi_2$ is an isomorphism.

To check that the isomorphism $\pi_2 \circ \pi_1^{-1}:
\mathscr{U}_1^{\an}/\mathscr{R}_1^{\an} \simeq
\mathscr{U}_2^{\an}/\mathscr{R}_2^{\an}$
is transitive with respect to a third choice of \'etale chart
for $\mathscr{X}$, it suffices to note that in the preceding considerations
we only needed that the \'etale chart $\mathscr{R}_{12} \rightrightarrows
\mathscr{U}_{12}$ dominates the other two charts, and not
that it is specifically their ``fiber product''.
\end{proof}

Lemma \ref{analgspace} permits us to make the following definition.

\begin{definition}
An algebraic space $\mathscr{X}$
is {\em analytifiable} if the \'etale
quotient $\mathscr{U}^{\an}/\mathscr{R}^{\an}$
exists for some $($and hence any$)$ \'etale chart
$\mathscr{R} \rightrightarrows \mathscr{U}$ for $\mathscr{X}$.
\end{definition}

For an analytifiable $\mathscr{X}$,
the {\em analytification}
$\mathscr{X}^{\an}$ is defined to be $\mathscr{U}^{\an}/\mathscr{R}^{\an}$.
Up to unique isomorphism, this
\'etale quotient
is independent of the specific choice of
\'etale chart $\mathscr{R} \rightrightarrows
\mathscr{U}$ (by Lemma \ref{analgspace}).
We now express the functoriality of analytification
when analytifications exist.

\begin{theorem}\label{fsep}
Let $\mathscr{X}$ and $\mathscr{X}'$ be analytifiable algebraic spaces
and let $f:\mathscr{X}' \rightarrow \mathscr{X}$ be a $k$-morphism.
Let $\mathscr{R} \rightrightarrows \mathscr{U}$
and $\mathscr{R}' \rightrightarrows \mathscr{U}'$ be respective
\'etale charts such that $f$ lifts to a map
$F:\mathscr{U}' \rightarrow \mathscr{U}$ for which $F \times F$
carries $\mathscr{R}'$ into $\mathscr{R}$
$($such a pair of charts always exists$)$.  The map of rigid spaces
$$f^{\an}:(\mathscr{X}')^{\an} \simeq (\mathscr{U}')^{\an}/(\mathscr{R}')^{\an}
\rightarrow \mathscr{U}^{\an}/\mathscr{R}^{\an} \simeq
\mathscr{X}^{\an}$$
induced by $F^{\an}$ depends only on $f$ and not on the \'etale charts
or the map $F$ lifting $f$, and this procedure enhances
the construction $\mathscr{X} \rightsquigarrow \mathscr{X}^{\an}$
to be a functor from the category of analytifiable
algebraic spaces over $k$ to
the category of rigid spaces over $k$.  Moreover:
\begin{itemize}
\item  the category of analytifiable algebraic spaces is
stable under the formation of fiber products and passage
to open and closed subspaces,
\item the functor
$\mathscr{X} \rightsquigarrow \mathscr{X}^{\an}$
is compatible with the formation of fiber products
and carries open/closed immersions to Zariski-open/closed immersions,
\item
if $f:\mathscr{X}' \rightarrow
\mathscr{X}$ is a morphism between analytifiable
algebraic spaces then $f^{\an}$ is separated if and only if
$f$ is separated.
\end{itemize}
\end{theorem}

\begin{proof}
Any two \'etale charts are dominated by a third,
and any two lifts of $f$ with respect
to a fixed choice of charts are $\mathscr{R}$-equivalent.
Thus, the well-definedness of
$f^{\an}$ is an immediate consequence of Lemma \ref{analgspace}.
The compatibility with composition of morphisms follows
from the independence of the choice of charts, so analytification
is indeed a functor on analytifiable algebraic spaces over $k$.

Let $\mathscr{X}$ be an analytifiable algebraic space
over $k$ and let $\mathscr{X}' \rightarrow \mathscr{X}$
be an open (resp. closed) immersion.   We let
$\mathscr{R} \rightrightarrows \mathscr{U}$ be an \'etale
chart for $\mathscr{X}$, and let
$\mathscr{U}' \rightarrow \mathscr{U}$ denote
the pullback of $\mathscr{X}'$.  The analytification of
this inclusion is a Zariski-open (resp. closed) immersion.
Rigid-analytic descent theory with respect to
coherent ideal sheaves trivially implies that
the analytification of $\mathscr{U}'$ descends to
a Zariski-open (resp. closed) immersion into $\mathscr{X}^{\an}$,
and this descent is easily seen to be an analytification of
$\mathscr{X}'$.

Now we consider fiber products.
Let $\mathscr{X}$ and
$\mathscr{Y}$ be algebraic spaces over an algebraic
space $\mathscr{Z}$, and assume that all three are analytifiable.
We need to check that $\mathscr{P} =
\mathscr{X} \times_{\mathscr{Z}} \mathscr{Y}$
is analytifiable and that the natural map
$\mathscr{P}^{\an} \rightarrow \mathscr{X}^{\an} \times_{\mathscr{Z}^{\an}}
\mathscr{Y}^{\an}$
is an isomorphism.  We choose
\'etale charts
$\mathscr{X} = \mathscr{U}'/\mathscr{R}'$,
$\mathscr{Y} = \mathscr{U}''/\mathscr{R}''$,
and $\mathscr{Z} = \mathscr{U}/\mathscr{R}$,
and we consider the fiber product
$$\mathscr{R}' \times_{\mathscr{R}} \mathscr{R}'' \rightrightarrows
\mathscr{U}' \times_{\mathscr{U}} \mathscr{U}''$$
as an \'etale chart for
$\mathscr{X} \times_{\mathscr{Z}} \mathscr{Y}$.
There is an evident \'etale surjective map
$$\pi:(\mathscr{U}' \times_{\mathscr{U}} \mathscr{U}'')^{\an} \simeq
(\mathscr{U}')^{\an} \times_{\mathscr{U}^{\an}}
(\mathscr{U}'')^{\an} \rightarrow
\mathscr{X}^{\an} \times_{\mathscr{Z}^{\an}} \mathscr{Y}^{\an},$$
and this admits local \'etale quasi-sections because
$\mathscr{X}^{\an}$, $\mathscr{Y}^{\an}$, and $\mathscr{Z}^{\an}$
are \'etale quotients associated to the analytifications of the
chosen \'etale charts for $\mathscr{X}$, $\mathscr{Y}$, and $\mathscr{Z}$
respectively.
The quotient
properties for $\mathscr{X}^{\an}$, $\mathscr{Y}^{\an}$, and
$\mathscr{Z}^{\an}$ in the category of rigid spaces permit us
to identify the fiber
square of $\pi$
with ${\mathscr{R}'}^{\rm{an}} \times_{\mathscr{R}^{\rm{an}}}
{\mathscr{R}''}^{\rm{an}} = (\mathscr{R}' \times_{\mathscr{R}}
{\mathscr{R}''})^{\rm{an}}$ in the evident manner, and so
analytification is compatible with fiber products.

Finally, we show that if $f:\mathscr{X}' \rightarrow
\mathscr{X}$ is a morphism between analytifiable algebraic
spaces then $f^{\an}$ is separated if and only if
$f$ is separated.  The compatibility of
analytification and fiber products
identifies $\Delta_{f^{\an}}$ and $\Delta_f^{\an}$, so
upon replacing $f$ with $\Delta_f$ we are reduced to checking
that if $f:\mathscr{X}' \rightarrow
\mathscr{X}$ is a quasi-compact immersion
between algebraic spaces then $f$ is a closed immersion if
and only if $f^{\an}$ is a closed immersion.
Let $h:\mathscr{U} \rightarrow \mathscr{X}$
be an \'etale covering by a scheme
and let $F:\mathscr{U}' \rightarrow \mathscr{U}$
be the base change of $f$ along $h$.  Since $f$ is a
quasi-compact immersion, so is $F$.  In particular,
$\mathscr{U}'$ is a scheme and $F^{\an}$ is
identified with the base change of $f^{\an}$
by means of $h^{\an}$.  The map $h^{\an}$
is an \'etale surjection with local \'etale quasi-sections, so
by rigid-analytic descent theory with respect
to the property of being a closed immersion we
see that $f^{\an}$ is a closed immersion if and only if
$F^{\an}$ is a closed immersion.  Likewise,
on the algebraic side we see that $f$ is a closed immersion
if and only if $F$ is a closed immersion.  Thus, we may
replace $f$ with $F$ to reduce to the case when $f$ is a quasi-compact
immersion between algebraic $k$-schemes.  This case
is handled by \cite[5.2.1(2)]{irred}.
\end{proof}

\begin{corollary}\label{opencheck}
Let $\mathscr{X}$ be an algebraic space
and let $\{\mathscr{X}_i\}$ be an open covering.
Analytifiability of $\mathscr{X}$ is equivalent
to that of all of the $\mathscr{X}_i$'s.
\end{corollary}

\begin{proof}
By Theorem \ref{fsep}, if $\mathscr{X}^{\an}$ exists
then so does $\mathscr{X}_i^{\an}$ (as a Zariski-open in
$\mathscr{X}^{\an}$) for all $i$.  Conversely, assume
that $\mathscr{X}_i^{\an}$ exists for all $i$.
The algebraic space $\mathscr{X}_{ij} = \mathscr{X}_i \cap \mathscr{X}_j$
is identified with a Zariski-open subspace of
both $\mathscr{X}_i$ and $\mathscr{X}_j$, so by
Theorem \ref{fsep} the rigid space $\mathscr{X}_{ij}^{\an}$
exists and is identified with a Zariski-open locus in
$\mathscr{X}_i^{\an}$ and $\mathscr{X}_j^{\an}$.
Since $\mathscr{X}_{ij} \cap \mathscr{X}_{ij'} =
\mathscr{X}_{ij} \times_{\mathscr{X}_i} \mathscr{X}_{ij'}$
for any $i, j, j'$,
the fiber-product compatibility in Theorem \ref{fsep}
provides the triple-overlap compatibility that is
required to glue the $\mathscr{X}_i^{\an}$'s to
construct a rigid space $X$
having the $\mathscr{X}_i^{\an}$'s as an admissible covering
with $\mathscr{X}_i^{\an} \cap \mathscr{X}_j^{\an} =
\mathscr{X}_{ij}^{\an}$ inside of $X$ for all $i$ and $j$.
In particular, $X - \mathscr{X}_i^{\an}$ meets
every $\mathscr{X}_j^{\an}$ in an analytic set,
and hence $\mathscr{X}_i^{\an}$ is Zariski-open in $X$.

We now check that $X$ serves as an analytification of
$\mathscr{X}$.  Let $p_1, p_2:\mathscr{R} \rightrightarrows \mathscr{U}$
be an \'etale chart for $\mathscr{X}$ and let
$\mathscr{U}_i$ be the preimage of $\mathscr{X}_i$ in $\mathscr{U}$.
Let $\mathscr{R}_i = p_1^{-1}(\mathscr{U}_i) \cap p_2^{-1}(\mathscr{U}_i)$.
Clearly $\mathscr{U}_i$ is $\mathscr{R}$-saturated in $\mathscr{U}$
and $\mathscr{R}_i \rightrightarrows \mathscr{U}_i$
is an \'etale chart for $\mathscr{X}_i$ via the
natural map $\mathscr{U}_i \rightarrow \mathscr{X}_i$.
Let $\mathscr{U}_{ij} = \mathscr{U}_i \cap \mathscr{U}_j$
and $\mathscr{R}_{ij} = p_1^{-1}(\mathscr{U}_{ij}) \cap
p_2^{-1}(\mathscr{U}_{ij})$, so
$\mathscr{U}_{ij}$ is $\mathscr{R}$-saturated and
$\mathscr{R}_{ij} \rightrightarrows \mathscr{U}_{ij}$
is an \'etale chart for $\mathscr{X}_{ij}$.
The gluing construction of $X$ implies that the maps
$$f_i:\mathscr{U}_i^{\an} \rightarrow \mathscr{X}_i^{\an}
\hookrightarrow X$$
satisfy $f_i|_{\mathscr{U}_{ij}^{\an}} =
f_j|_{\mathscr{U}_{ij}^{\an}}$,
so the $f_i$'s uniquely glue to define an \'etale surjection
$f:\mathscr{U}^{\an} \rightarrow X$ that restricts
to the canonical map
$\mathscr{U}_i^{\an} \rightarrow \mathscr{X}_i^{\an}$
over each Zariski-open $\mathscr{X}_i^{\an} \subseteq
X$ and hence (by the definition of
$\mathscr{X}_i^{\an}$ as the quotient
$\mathscr{U}_i^{\an}/\mathscr{R}_i^{\an}$) the map $f$ admits
local \'etale quasi-sections.  Since
the $\mathscr{R}_i^{\an}$'s form a Zariski-open
covering of $\mathscr{R}^{\an}$, it is easy to check
that the composite maps $\mathscr{R}^{\an} \rightrightarrows
\mathscr{U}^{\an} \rightarrow X$ coincide.  Thus,
$\mathscr{R}^{\an}$ is naturally a rigid space
over $X$ and we obtain a canonical
$X$-map $h:\mathscr{R}^{\an} \rightarrow
\mathscr{U}^{\an} \times_X \mathscr{U}^{\an}$.
The restriction of $h$ over $\mathscr{X}_i^{\an}$
is the canonical map $\mathscr{R}_i^{\an} \rightarrow
\mathscr{U}_i^{\an} \times_{\mathscr{X}_i^{\an}} \mathscr{U}_i^{\an}$
that is an isomorphism (due to the definition of
$\mathscr{X}_i^{\an}$), so $h$ is an isomorphism.  Hence,
$X$ indeed serves as an analytification of $\mathscr{X}$.
\end{proof}

\subsection{Properties of analytification}\label{propsec}

We now summarize some basic observations concerning properties
of analytification of algebraic spaces, especially in connection
with properties of morphisms between algebraic spaces.

If $f:\mathscr{X} \rightarrow \mathscr{Y}$
is a faithfully flat map between analytifiable algebraic spaces over $k$,
then we claim that the induced faithfully flat
map $f^{\an}:\mathscr{X}^{\an}
\rightarrow \mathscr{Y}^{\an}$ has local {\em fpqc} quasi-sections,
and if $f$ is an \'etale surjection then we claim that
$f^{\an}$ has local \'etale quasi-sections.  This
generalizes Example \ref{coverex} to the case
of algebraic spaces.
To prove these claims, we pick a chart $\mathscr{R} \rightrightarrows
\mathscr{U}$ for $\mathscr{Y}$, so
$\mathscr{U}^{\an} \rightarrow \mathscr{Y}^{\an}$
has local \'etale quasi-sections (as
$\mathscr{Y}^{\an}$ is the \'etale quotient
$\mathscr{U}^{\an}/\mathscr{R}^{\an}$).  Hence,
we may replace $\mathscr{X} \rightarrow \mathscr{Y}$
with its base change by
the \'etale surjection $\mathscr{U} \rightarrow \mathscr{Y}$,
so we can assume that $\mathscr{Y}$ is an algebraic $k$-scheme.
Running through a similar argument with an \'etale chart
for $\mathscr{X}$ reduces us to the case when $\mathscr{X}$
is also an algebraic $k$-scheme, and so we are brought to the
settled scheme case.

\begin{theorem}\label{propermap} If
$f:\mathscr{X} \rightarrow \mathscr{Y}$ is a map between analytifiable algebraic spaces over $k$, then $f$ has
property $\P$ if and only if $f^{\an}$ has
property $\P$, where $\P$ is any of the following properties: separated, monomorphism, surjective, isomorphism,
open immersion, flat, smooth, and \'etale.  Likewise, if $f$ is finite type then we may take $\P$ to be: closed
immersion, finite, proper, quasi-finite $($i.e., finite fibers$)$.
\end{theorem}

\begin{proof}  By
\cite[Thm.~4.2.7]{relamp} and descent
theory for schemes, we may work
\'etale-locally on $\mathscr{Y}$ and so
we can assume that $\mathscr{Y}$ is a scheme of finite type over $k$.
Since flat maps locally of finite type between algebraic spaces are open,
the essential properties to consider are isomorphism and properness;
the rest then follow exactly as in the case of schemes.
By Chow's lemma for algebraic spaces and \cite[\S{}A.1]{relamp}
(for properness), the proper case is
reduced to the case of
quasi-compact immersions of schemes (more specifically,
a quasi-compact immersion into a projective space
over $\mathscr{Y}$), and this case follows from
\cite[5.2.1(2)]{irred}.  If $f^{\an}$ is
an isomorphism then $f$ is quasi-finite, flat,
and (by Theorem \ref{fsep}) separated, so
$\mathscr{X}$ is necessarily a scheme.  Thus,
we may use \cite[5.2.1(1)]{irred} to infer that $f$ is an
isomorphism.
\end{proof}

\begin{example}\label{stalkcom} Let $\mathscr{X}$ be an analytifiable algebraic space.
Recall that a {\em point} of an algebraic space is a monic morphism
from the spectrum of a field, and the underlying topological
space $|\mathscr{X}|$ of an algebraic space is its
set of (isomorphism classes of) points.   Every map to $\mathscr{X}$ from the spectrum
of a field factors through a unique point of $\mathscr{X}$
\cite[II,~6.2]{knutson}, so
distinct points
have empty fiber product over $\mathscr{X}$.  Moreover,
every point of an algebraic space admits
a residually-trivial pointed \'etale scheme neighborhood
\cite[II,~6.4]{knutson}.

 With respect to the usual Zariski topology
on $|\mathscr{X}|$, it follows by working with
an \'etale scheme cover that a point $\Spec(k') \rightarrow \mathscr{X}$ is closed
in $|\mathscr{X}|$
if and only if $[k':k]$ is finite.
 In particular, if $j:\Spec(k') \rightarrow \mathscr{X}$
is a closed point then analytification defines a map
$j^{\rm{an}}:\Sp(k') \rightarrow X := \mathscr{X}^{\rm{an}}$ that is a monomorphism
(since $\Delta_{j^{\rm{an}}} = \Delta_j^{\rm{an}}$ is an isomorphism).
The monic map $j^{\rm{an}}$ is easily seen to be an isomorphism
onto an ordinary point of $X$, so in this way we have defined a map
of sets $|\mathscr{X}|_0 \rightarrow X$, where $|\mathscr{X}|_0$ is the
set of closed points of $\mathscr{X}$.   By a fiber product argument we
see that this map is bijective, and by construction it preserves
residue fields over $k$.

Recall that at any point of an algebraic space there is a naturally associated
henselian local ring with the same residue field (use the limit
of residually-trivial pointed \'etale scheme neighborhoods;
for schemes this is the henselization of
the usual local ring).  Thus, for
any coherent sheaf
$\mathscr{F}$ on $\mathscr{X}$
and closed point $x_0 \in |\mathscr{X}|_0$,
we get a stalk $\mathscr{F}_{x_0}$ that is a finite module
over the henselian local noetherian $k$-algebra $\mathscr{O}_{\mathscr{X},x_0}$.
If $x \in X$ is the corresponding point of $X$ then every residually
trivial pointed \'etale neighborhood $(\mathscr{U},x_0) \rightarrow (\mathscr{X},x_0)$
with analytifiable $\mathscr{U}$ (e.g., a scheme) induces an \'etale
map $\mathscr{U}^{\rm{an}} \rightarrow \mathscr{X}^{\rm{an}}$
that is an isomorphism between small admissible opens around
the canonical copy of $x$ in each (due to residual triviality).   Hence,
we get a natural map of $k$-algebras
$\mathscr{O}_{\mathscr{X}}(\mathscr{U}) \rightarrow \mathscr{O}_{X,x}$
compatible with the $k(x)$-points on each, and passage to the direct limit
thereby defines a local $k$-algebra map $\mathscr{O}_{\mathscr{X},x_0}
\rightarrow \mathscr{O}_{X,x}$.  There is likewise a map
$\mathscr{F}_{x_0} \rightarrow \mathscr{F}^{\rm{an}}_x$
that is linear over this map of rings, and so we get a natural
$\mathscr{O}_{X,x}$-linear {\em comparison map}
$$\mathscr{O}_{X,x} \otimes_{\mathscr{O}_{\mathscr{X},x_0}}
\mathscr{F}_{x_0} \rightarrow \mathscr{F}^{\rm{an}}_x.$$
This is functorial in $(\mathscr{X},x_0)$, so it can be computed
on a residually-trivial pointed \'etale scheme neighborhood of $x_0$.
Hence, by reduction to the scheme case
we see that this comparison map between finite
$\mathscr{O}_{X,x}$-modules is an isomorphism because
the induced map on completions is an isomorphism
(due to the isomorphism between completions of algebraic and
analytic local rings at a common point of an algebraic $k$-scheme).
\end{example}

\begin{example}\label{ex27}
Since reducedness is inherited under analytification of reduced
algebraic $k$-schemes, by using an \'etale scheme cover we
see that if $\mathscr{X}$ is an analytifiable reduced algebraic space
then $\mathscr{X}^{\rm{an}}$ is reduced.  Thus, if $\mathscr{X}$
is any analytifiable algebraic space (so
the closed subspace $\mathscr{X}_{\rm{red}}$ is
also analytifiable) then there is a natural closed immersion
$(\mathscr{X}_{\rm{red}})^{\rm{an}} \hookrightarrow (\mathscr{X}^{\rm{an}})_{\rm{red}}$
whose formation is compatible with \'etale base change on
$\mathscr{X}$, so it is an isomorphism by reduction to the scheme case.

As an application of this compatibility, if $\mathscr{X}$ is any analytifiable algebraic space
and $\{\mathscr{X}_i\}$
is its set of irreducible components with reduced structure
then we get closed immersions $\mathscr{X}_i^{\rm{an}} \hookrightarrow
\mathscr{X}^{\rm{an}}$ which we claim are the irreducible components of
$\mathscr{X}^{\rm{an}}$ endowed with their reduced structure.  To check this,
since $(\mathscr{X}_i \times_{\mathscr{X}} \mathscr{X}_j)^{\rm{an}}
\simeq \mathscr{X}_i^{\rm{an}} \times_{\mathscr{X}^{\rm{an}}}
\mathscr{X}_j^{\rm{an}}$ we see by dimension reasons that it suffices
to show that $\mathscr{X}^{\rm{an}}$ is irreducible when $\mathscr{X}$ is irreducible
and reduced.
There is a dense Zariski-open subspace $\mathscr{X}_0 \subseteq
\mathscr{X}$ that is a scheme, and by \cite[2.3.1]{irred}
the rigid space $\mathscr{X}_0^{\rm{an}}$ is irreducible.
Thus, by the existence of global irreducible
decomposition for rigid spaces,
it suffices to show that the Zariski-open locus $\mathscr{X}_0^{\rm{an}}$
in $\mathscr{X}^{\rm{an}}$ is everywhere dense.   Consideration of
points valued in finite extensions of $k$ shows that
the analytic set in $\mathscr{X}^{\rm{an}}$ complementary to $\mathscr{X}_0^{\rm{an}}$
is $\mathscr{Z}^{\rm{an}}$, where $\mathscr{Z} = \mathscr{X} - \mathscr{X}_0$
(say, with its reduced structure).   Thus, we just need that the coherent ideal
of $\mathscr{Z}^{\rm{an}}$ in the reduced $\mathscr{X}^{\rm{an}}$ is nowhere zero,
and this is clear by comparing its completed stalks with those of $\mathscr{Z}$
in the irreducible and reduced algebraic space $\mathscr{X}$.
(Here we use the end of Example \ref{stalkcom} and the fact that
a proper closed subspace of an irreducible algebraic space has
nowhere-dense pullback to any \'etale scheme cover, even though
irreducibility is not local for the \'etale topology.)
\end{example}

If $X$ is a complex-analytic space
then since $X$ rests on an ordinary topological space,
we can find an open subset $V \subseteq X \times X$ such that the diagonal
$X \rightarrow X \times X$ factors through a closed immersion into $V$.
(Take $V$ to be a union $\cup (U_i \times U_i)$
where the $U_i$ are Hausdorff open sets that cover $X$.)
That is, $\Delta_X$ factors like an immersion in algebraic geometry.
This suggests that local separatedness of $\mathscr{X}$ (i.e., the diagonal $\Delta_{\mathscr{X}}$
being an immersion) should be a necessary condition for the analytifiability
of an algebraic space $\mathscr{X}$
over $\C$.   Due to lack of a reference, we now give a proof of this
fact and its non-archimedean counterpart.

\begin{theorem}\label{locsep} Let $\mathscr{X}$ be an algebraic space locally of finite
type over either $\C$ or a non-archimedean field $k$.   If
$\mathscr{X}^{\rm{an}}$ exists then $\mathscr{X}$ must be locally separated.
\end{theorem}

See Theorem \ref{berklocsep} for an analogous result for $k$-analytic spaces.

\begin{proof}
Choose a scheme $\mathscr{U}$ equipped
with an \'etale surjection $\mathscr{U} \twoheadrightarrow \mathscr{X}$.
Let $\mathscr{R} = \mathscr{U} \times_{\mathscr{X}} \mathscr{U}$, so
our aim is to prove that the monomorphism
$i:\mathscr{R} \rightarrow \mathscr{U} \times \mathscr{U}$ is an immersion.
The map $i$ is separated since it is a monomorphism.
Note that $i$ is a base change of the quasi-compact diagonal $\Delta_{\mathscr{X}}$.
By the compatibility of analytification and fiber products (for algebraic
spaces in both the complex-analytic and non-archimedean cases),
$\mathscr{X}^{\rm{an}} \times \mathscr{X}^{\rm{an}}$ is identified with
an analytification of $\mathscr{X} \times \mathscr{X}$, and in this
way $\Delta_{\mathscr{X}}^{\rm{an}}$ is identified with $\Delta_{\mathscr{X}^{\rm{an}}}$.
Thus, the map $i^{\rm{an}}$ is a base
change of $\Delta_{\mathscr{X}^{\rm{an}}}$.

We first treat the complex-analytic case, and then we adapt the argument
to work in the non-archimedean case.   The diagonal map of a
complex-analytic space is a topological embedding.  Indeed,
if $S$ is such a space then for each $s \in S$ there is a Hausdorff open neighborhood
$U_s \subseteq S$ around $s$, so $U = \cup (U_s \times U_s)$ is an open subset
of $S \times S$ through which $\Delta_S$ factors as a closed immersion
$S \hookrightarrow U$.  Hence, $\Delta_S$ is a topological embedding.
Applying this to $S = \mathscr{X}^{\rm{an}}$ gives that $\Delta_{\mathscr{X}^{\rm{an}}}$
is a topological embedding, so its topological base change $i^{\rm{an}}$ is
also a topological embedding.
It therefore suffices to show that if $f:\mathscr{Y} \rightarrow \mathscr{S}$
is a quasi-compact monomorphism between algebraic $\C$-schemes and
$f^{\rm{an}}$ is a topological embedding then $f$ is an immersion.
Since $f$ is monic, so $\Delta_f$ is an isomorphism, $f$ is certainly separated.
We can assume that $\mathscr{S}$ is separated
and quasi-compact (e.g., affine), so $\mathscr{Y}$ is separated.
By Zariski's Main Theorem \cite[IV$_3$,~8.12.6]{ega},
the quasi-finite separated map $f$ factors as a composition
$$\mathscr{Y} \stackrel{j}{\hookrightarrow} \mathscr{T}
\stackrel{\pi}{\rightarrow} \mathscr{S}$$
where $j$ is an open immersion and
$\pi$ is finite.   We can replace
$\mathscr{T}$ with the schematic closure of $\mathscr{Y}$ to arrange
that $j$ has dense image.
Thus, $j^{\rm{an}}$ is an open immersion with
dense image with respect to the analytic topology
\cite[XII,~\S2]{sga1}.  We can also replace $\mathscr{S}$
with the natural closed subscheme structure on the closed set
$\pi(\mathscr{T})$ to arrange that $\pi$ is surjective.

We now show that $j^{\rm{an}}(\mathscr{Y}^{\rm{an}})$ is the preimage of its image under $\pi^{\rm{an}}$, and
that the restriction of $\pi^{\rm{an}}$ to $j^{\rm{an}}(\mathscr{Y}^{\rm{an}})$ is injective. Choose $y \in
\mathscr{Y}^{\rm{an}}$ and a point $t \in \mathscr{T}^{\rm{an}}$ such that $\pi^{\rm{an}}(t) =
\pi^{\rm{an}}(j^{\rm{an}}(y)) = f^{\rm{an}}(y)$ in $\mathscr{Y}^{\rm{an}}$. We will prove $t = j^{\rm{an}}(y)$.
By denseness of $j^{\rm{an}}(\mathscr{Y}^{\rm{an}})$ in the Hausdorff analytic space $\mathscr{T}^{\rm{an}}$,
there is a sequence of points $y_n \in \mathscr{Y}^{\rm{an}}$ such that $j^{\rm{an}}(y_n) \rightarrow t$, so
applying $\pi^{\rm{an}}$ gives that $f^{\rm{an}}(y_n) \rightarrow f^{\rm{an}}(y)$ in the Hausdorff analytic
space $\mathscr{S}^{\rm{an}}$.  But $f^{\rm{an}}$ is a topological embedding by hypothesis, so $y_n \rightarrow
y$ in the Hausdorff analytic space $\mathscr{Y}^{\rm{an}}$. Applying $j^{\rm{an}}$ gives $j^{\rm{an}}(y_n)
\rightarrow j^{\rm{an}}(y)$, but $t$ is the limit of this sequence, so $t = j^{\rm{an}}(y)$ as required.
Returning to the algebraic setting, the Zariski-open set $j(\mathscr{Y}) \subseteq \mathscr{T}$ is the preimage
of its image under the finite map $\pi$ because this asserts an equality of constructible sets and it suffices
to check such an equality on closed points.   Since $\pi$ is a finite surjection, hence a topological quotient
map, it follows that $f(\mathscr{Y}) = \pi(j(\mathscr{Y}))$ is an open subset of $\mathscr{S}$ since its
preimage $j(\mathscr{Y})$ in $\mathscr{T}$ is a Zariski-open subset. We can replace $\mathscr{S}$ with the open
subscheme structure on this open set and replace $\mathscr{T}$ with its preimage under $\pi$, so the open
immersion $j$ is now surjective on closed points.  Thus, $j$ is an isomorphism, so the monomorphism $f = \pi
\circ j$ is finite and therefore is a closed immersion.

Next we consider the non-archimedean case.  We need to exercise some more care with respect to topological
arguments, due to the role of the Tate topology and the fact that fiber products of rigid spaces are generally
not fiber products on underlying sets (unless the base field is algebraically closed).    As a first step, if
$X$ is a rigid space and $\Delta:X \rightarrow X \times X$ is its diagonal, then for any rigid-analytic morphism
$h:Z \rightarrow X \times X$ we claim that the induced map of rigid spaces $h^{\ast}(\Delta):Z \times_{X \times
X} X \rightarrow Z$ is an embedding with respect to the canonical (totally disconnected) topology. To prove
this, let $\{X_j\}$ be an admissible affinoid open covering of $X$ and let $U = \cup_j (X_j \times X_j)$, so $U$
is an open subset of $X \times X$ with respect to the canonical topology (though it is unclear if $U$ must be an
admissible open subset).
The map $\Delta$ factors through a continuous map $\Delta':X \rightarrow U$
under which the preimage of $X_j \times X_j$ is $X_j \subseteq X$, and which induces the diagonal $\Delta_{X_j}$
that is a closed immersion (since $X_j$ is affinoid).  Hence, $\Delta'$
is a closed embedding with respect to the canonical topologies. The preimage $h^{-1}(U) \subseteq Z$ is open
with respect to the canonical topology on $Z$, and $h^{\ast}(\Delta)$ factors
continuously through $h^{-1}(U)$ (as a map with
respect to the canonical topologies). The resulting continuous map $Z \times_{X \times X} X \rightarrow
h^{-1}(U)$ with respect to the canonical topologies is a closed embedding because on each of the opens
$h^{-1}(X_j \times X_j) \subseteq h^{-1}(U)$ that cover $h^{-1}(U)$ it restricts to the rigid-analytic morphism
$h^{-1}(X_j \times X_j) \times_{X_j \times X_j} X_j \rightarrow h^{-1}(X_j \times X_j)$ that is a base change of
the closed immersion $\Delta_{X_j}$.
Hence, for any $h$ as above, $h^{\ast}(\Delta)$
is an embedding with result to the canonical topologies,
as claimed.  As a special case, the analytification of the monic map $i:\mathscr{R}
\rightarrow \mathscr{U} \times \mathscr{U}$ is a topological embedding with respect to the canonical topologies
since it is a base change of $\Delta_{\mathscr{X}^{\rm{an}}}$.

Thus, exactly as in the complex-analytic case,
it suffices to prove that if $f:\mathscr{Y} \rightarrow \mathscr{S}$ is a quasi-compact monomorphism between
algebraic $k$-schemes such that $f^{\rm{an}}$ is a topological embedding with respect to the canonical
topologies then $f$ is an immersion.
The denseness result cited above from \cite[XII,~\S2]{sga1} remains
valid in the non-archimedean case (with the same proof) when using the canonical topology,
so the preceding argument in the complex-analytic case carries over {\em verbatim} once
we show that a separated rigid space is Hausdorff with respect to its
underlying canonical topology.
(The preceding argument with sequences works in the rigid-analytic case because
any point in a rigid space has a countable base of open neighborhoods,
though this countability property is not needed; we could instead use nets.)
Let $S$ be a separated rigid space, and let $|S|$ denote the underlying set
with the canonical topology.   To prove that it is Hausdorff,  we have
to overcome the possibility that $|S \times S| \rightarrow |S| \times |S|$
may not be bijective (let alone a homeomorphism with respect to the canonical topologies).
Choose $s, s' \in |S|$
and let $U$ and $U'$ be admissible affinoid opens that contain $s$ and $s'$ respectively.
The overlap $U \cap U'$ is affinoid and the union $U \cup U'$ is an admissible open
subspace with $\{U, U'\}$ as an admissible cover (since $S$ is separated).
We can replace $S$ with $U \cup U'$, so we can assume that $S$ is quasi-compact
and separated.  Let $S^{\rm{an}}$ denote the associated separated
strictly $k$-analytic space (in the sense of Berkovich), so
the underlying topological space $|S^{\rm{an}}|$ (with its ``spectral'' topology)
is compact Hausdorff.  The subset of classical points
(i.e., those $s \in S^{\rm{an}}$ with $[\mathscr{H}(s):k] < \infty$)
is $|S|$, and its subspace topology is the canonical topology
because this is true when $S$ is replaced
with each of finitely many (compact Hausdorff) strictly $k$-analytic affinoid domains
that cover $S^{\rm{an}}$.
Hence, $|S|$ is also Hausdorff.
\end{proof}

In view of the known converse to Theorem \ref{locsep} in the complex-analytic case,
it is natural to ask if the property of being locally separated is sufficient
for the analytifiable of an algebraic space over $k$.
We will see that this is not true, via counterexamples
in \S\ref{countersec}.   Such counterexamples will be explained in two ways, using
rigid-analytic methods and using $k$-analytic spaces.  To carry
out a rigid-analytic approach, we need to consider how analytifiability behaves
with respect to change in the base field.
A relevant notion for this was used in \cite{relamp}:
a {\em pseudo-separated} map $f:X \rightarrow S$ between rigid spaces
is a map whose diagonal $\Delta_f:X \rightarrow X \times_S X$ factors
as the composite of a Zariski-open immersion followed by a closed immersion.
The reason that we choose this order of composition is that in the scheme
case it is available in a canonical manner (via scheme-theoretic closure)
and hence behaves well with respect to \'etale localization and descent for schemes
(so it generalizes to define the notion of quasi-compact immersion for algebraic spaces).
We note that a map of rigid spaces is pseudo-separated and quasi-separated
(i.e., has quasi-compact diagonal) if and only if it is separated.

\begin{lemma}\label{4310}
Let $\mathscr{X}$ be an analytifiable
algebraic space.  The analytification $\mathscr{X}^{\an}$ is pseudo-separated.
In particular, $k' \widehat{\otimes}_k \mathscr{X}^{\an}$ makes sense
as a rigid space over $k'$ for any analytic extension field $k'/k$.
\end{lemma}

\begin{proof}
Since $\mathscr{X}$ is analytifiable, so is
$\mathscr{X} \times \mathscr{X}$.  Clearly
$\Delta_{\mathscr{X}}^{\rm{an}} = \Delta_{\mathscr{X}^{\rm{an}}}$.
By Theorem \ref{locsep}, $\mathscr{X}$ is locally separated.  Thus,
by \'etale descent,
the quasi-compact immersion
$\Delta_{\mathscr{X}}$ uniquely
factors as a (schematically dense) Zariski-open immersion
followed by a closed immersion, and these intervening
locally closed subspaces of $\mathscr{X}$ must be analytifiable (by Theorem \ref{fsep}).
It follows that
$\mathscr{X}^{\an}$ must be pseudo-separated.
\end{proof}

Now we can address the question of compatibility of analytification
for algebraic spaces and extension of the base field,
generalizing the known case of analytification of algebraic $k$-schemes.
This compatibility will be useful for giving a rigid-analytic
justification of our counterexamples
to analytifiability.

\begin{theorem}\label{4311} Let $\mathscr{X}$ be an analytifiable
algebraic space over $k$, and let $k'/k$ be an analytic extension field.
The algebraic space $k' \otimes_k \mathscr{X}$
over $k'$ is analytifiable and there is
a natural isomorphism
$k' \widehat{\otimes}_k \mathscr{X}^{\an} \simeq
(k' \otimes_k \mathscr{X})^{\an}$ which is the usual isomorphism when
$\mathscr{X}$ is an algebraic $k$-scheme.  These
natural isomorphisms are
transitive with respect to further analytic extension of
the base field and are compatible with the formation of fiber products.
\end{theorem}

\begin{proof}  Let $\mathscr{R} \rightrightarrows \mathscr{U}$
be an \'etale chart for $\mathscr{X}$.  Since
$\mathscr{U}^{\an} \rightarrow \mathscr{X}^{\an}$
is \'etale and admits local \'etale quasi-sections,
the same holds for the map $k' \widehat{\otimes}_k
\mathscr{U}^{\an} \rightarrow k' \widehat{\otimes}_k
\mathscr{X}^{\an}$ because (as we see via formal models)
a faithfully flat map between $k$-affinoids induces
a faithfully flat map after applying $k' \widehat{\otimes}_k (\cdot)$.
 Moreover, the natural map
$$k' \widehat{\otimes}_k \mathscr{R}^{\an} \rightarrow
(k' \widehat{\otimes}_k \mathscr{U}^{\an})
\times_{k' \widehat{\otimes}_k \mathscr{X}^{\an}}
(k' \widehat{\otimes}_k \mathscr{U}^{\an})$$
is an isomorphism because it
is identified with the extension of scalars of
the map $\mathscr{R}^{\an} \rightarrow
\mathscr{U}^{\an} \times_{\mathscr{X}^{\an}} \mathscr{U}^{\an}$
that is necessarily an isomorphism (due to the defining
property of the \'etale quotient $\mathscr{X}^{\an}$
that we are assuming to exist).
Thus, we conclude that $k' \widehat{\otimes}_k \mathscr{X}^{\an}$
serves as an \'etale quotient for the diagram
\begin{equation}\label{kru}
k' \widehat{\otimes}_k \mathscr{R}^{\an} \rightrightarrows
k' \widehat{\otimes}_k \mathscr{U}^{\an}
\end{equation}
that is an \'etale equivalence relation, due
to its identification with
the analytification $(k' \otimes_k \mathscr{R})^{\an} \rightrightarrows
(k' \otimes_k \mathscr{U})^{\an}$ of an \'etale chart
for the algebraic space $k' \otimes_k \mathscr{X}$ over $k'$.
This shows that $k' \widehat{\otimes}_k \mathscr{X}^{\an}$
naturally serves as an analytification for the algebraic space
$k' \otimes_k \mathscr{X}$ over $k'$.  Moreover, it is clear
that this identification $k' \widehat{\otimes}_k \mathscr{X}^{\an}
\simeq (k' \otimes_k \mathscr{X})^{\an}$ is independent of
the choice of \'etale chart $\mathscr{U} \rightrightarrows
\mathscr{R}$ for $\mathscr{X}$ and that it is therefore
functorial in the analytifiable $\mathscr{X}$.
The compatibility
with fiber products (when the relevant analytifications
exist over $k$) is now obvious.
\end{proof}

We can define an exact analytification functor from
coherent
$\O_{\mathscr{X}}$-modules to coherent $\O_{\mathscr{X}^{\an}}$-modules
in two (equivalent) ways.   One method is to use
an \'etale chart $\mathscr{R} \rightrightarrows \mathscr{U}$ and
descent theory for coherent sheaves on rigid spaces \cite[4.2.8]{relamp};
it is easy to see that this approach gives a functor that is independent of
the choice of \'etale chart and that  is natural in $\mathscr{X}$
with respect to pullback along maps of algebraic spaces
$\mathscr{X}' \rightarrow \mathscr{X}$.  An alternative method
that avoids the crutch of an \'etale chart will be explained above
Example \ref{gagaex}.

\section{Analytification counterexamples
and constructions}\label{countersec}

To show that the theory in \S\ref{algspacesec} is not vacuous, we need
to prove the analytifiability of an interesting class of
algebraic spaces that are not
necessarily schemes.  We also explain
how analytification interacts with \'etale topoi,
as is required for GAGA on algebraic spaces over $k$
as well as for a natural definition of
analytification for coherent sheaves on algebraic spaces.

\subsection{Non-analytifiable surfaces}

By Theorem \ref{locsep}, it is necessary to restrict attention
to those algebraic spaces (locally of finite type over $k$)
that are locally separated.  It will be proved
in Theorem \ref{makespace} that
such algebraic spaces are analytifiable in the separated case,
but we first give locally separated examples where analytifiability fails.
Our construction will provide locally separated smooth algebraic
spaces $\mathscr{S}$ over $\Q$ such that $\mathscr{S}_{k}^{\rm{an}}$
does not exist for any non-archimedean field $k/\Q$, though of course
$\mathscr{S}_{\C}^{\rm{an}}$ does exist by local separatedness!

\begin{example}\label{counter}
Let $k$ be an abstract field, let
$\mathscr{T} \subseteq \mathbf{A}^2_k$ be a dense open subset of the $x$-axis,
and let $\mathscr{T}' \rightarrow \mathscr{T}$ be the geometrically connected
finite \'etale covering $\{u^d = f(x)\}$ with degree $d > 1$ given by
extracting the $d$th root of a monic separable polynomial $f \in k[x]$
whose zeros are away from $\mathscr{T}$; we assume $d$ is not divisible by
the characteristic of ${k}$.   By shrinking $\mathscr{T}$ near its
generic point, we can find an open
$\mathscr{X} \subseteq \mathbf{A}^2_k$ in which $\mathscr{T}$ is closed
and over which there is a quasi-compact \'etale cover $\mathscr{U} \rightarrow
\mathscr{X}$
restricting to $\mathscr{T}' \rightarrow \mathscr{T}$ over
$\mathscr{T}$.    There is a locally separated
algebraic space $\mathscr{X}'$ that is \'etale over $\mathscr{X}$
(hence $k$-smooth and 2-dimensional) and obtained from
the open $\mathscr{X} \subseteq \mathbf{A}^2_k$ by replacing the curve
$\mathscr{T}$ with the degree-$d$ covering $\mathscr{T}'$.

More generally, in \cite[Intro.,~Ex.~2,~pp.10--12]{knutson} Knutson gives the following construction. For any
quasi-separated scheme $\mathscr{X}$ equipped with a closed subscheme $\mathscr{T}$ and a quasi-compact \'etale
surjection $\pi:\mathscr{U} \rightarrow \mathscr{X}$, he builds  a locally separated algebraic space
$\mathscr{X}'$ equipped with a quasi-compact \'etale surjection $i:\mathscr{X}' \rightarrow \mathscr{X}$ such
that $i$ is an isomorphism over $\mathscr{X} - \mathscr{T}$ but has pullback to $\mathscr{T}$ given by the
\'etale covering $\mathscr{T}' = \pi^{-1}(\mathscr{T}) \rightarrow \mathscr{T}$. When $\mathscr{T}'$ and
$\mathscr{T}$ are irreducible and $\mathscr{T}'$ has generic degree $d > 1$ over $\mathscr{T}$ then the behavior
of fiber-rank for the quasi-finite map $\mathscr{X}' \rightarrow \mathscr{X}$ is opposite to
what happens for quasi-finite separated \'etale maps of schemes (via the structure theorem for such maps
\cite[IV$_4$,~18.5.11]{ega}) in the sense that the fiber-degree goes up (rather than down) at special points.
Hence, $\mathscr{X}'$ cannot have an open scheme neighborhood (or equivalently, a separated open neighborhood)
around any point of $\mathscr{T}'$ in such cases, since if it did then such a neighborhood would contain the
1-point generic fiber over $\mathscr{T}$, yet no fiber over $\mathscr{T}$ can have a separated open neighborhood
in $\mathscr{X}'$ (e.g., an affine open subscheme).

We now consider the construction over $\mathscr{X} \subseteq \mathbf{A}^2_k$
that replaces an open subset $\mathscr{T}$ in
the $x$-axis with a degree-$d$ finite \'etale covering $\{u^d = f(x)\}$ as above.
For $k = \C$, an analytification of $\mathscr{X}'$ does exist (since $\mathscr{X}'$ is locally
separated) and its local structure over an open neighborhood of
a point $t \in \mathscr{T}(\C)$ is very easy to describe: it is a product
of an open unit disc with a gluing of $d$ copies of an open disc
to itself via the identity map on the complement of the origin.   In particular,
this analytification is non-Hausdorff over such a neighborhood.
In the non-archimedean setting, if $\mathscr{T}' \rightarrow \mathscr{T}$
has a non-split fiber over some $t \in \mathscr{T}(k)$ then
we will see that no analogous such local gluing can be
done over an open neighborhood of $t$ in $\mathscr{X}^{\rm{an}}$.
If $k$ is algebraically closed then the local gluing can be done but
there are global admissibility problems with the gluing.  We shall show that
the admissibility problems are genuine.  The key is that there is an obstruction to analytifiability
caused by the failure of the Gelfand--Mazur theorem over non-archimedean fields:
$k$ admits analytic extension fields $k'/k$ such that  the \'etale cover
$\mathscr{T}' \rightarrow \mathscr{T}$ has a non-split fiber over some $t \in \mathscr{T}(k')$,
even if $k$ is algebraically closed.

Let us now prove that if $k$ is a non-archimedean field then
the smooth 2-dimensional locally separated algebraic space $\mathscr{X}'$ is
not analytifiable.  Assume that an analytification
$X'$ of $\mathscr{X}'$ exists.  By Lemma \ref{4310}
the rigid space $X'$ must be pseudo-separated, and by Theorem \ref{4311}
if $k'/k$ is any analytic extension field then
$k' \otimes_k \mathscr{X}'$ is analytifiable with analytification $k' \widehat{\otimes}_k X'$.
Hence, to get a contradiction it suffices to consider the situation after
a preliminary analytic extension of the base field $k \rightarrow k'$ (which is easily checked to
automatically commute with the formation of $\mathscr{X}'$ in terms of
$\mathscr{X}$ and $\mathscr{T}' \rightarrow \mathscr{T}$).
The extension $K' = k(\mathscr{T}')$ of $K = k(\mathscr{T}) = k(x)$
is defined by adjoining a root to the irreducible polynomial $u^d - f \in K[u]$.
We first increase $k$ a finite amount so that $f$ splits completely in $k[x]$.
We then make a linear change of variable on $x$
so that $f = x \prod_{i > 1} (1 - r_i x)$
with $|r_i| < 1$ for each $i$.  In case of mixed characteristic
we also require $|r_i|$ to be so small that
$1 - r_i x$ has a $d$th root as a power series (a condition that
is automatic in case of equicharacteristic $k$).    Thus,
$f$ lies in the valuation ring of the Gauss norm on $k(x)$ and its image in
the residue field $\widetilde{k}(x)$ for the Gauss norm
(with $\widetilde{k}$ the residue field of $k$) is $x$, which is not a $d$th power
in $\widetilde{k}(x)$ since $d > 1$.
Hence, $u^d - f$ has no root over the completion $\widehat{K}$ of
$K = k(x)$ with respect to the Gauss norm.    By taking
$k' = \widehat{K}$ and working with the canonical
point in $\mathscr{T}(\widehat{K})$ arising
from the generic point of $\mathscr{T}$ we put ourselves
in the situation (upon renaming $k'$ as $k$) where there exists $t_0 \in \mathscr{T}(k)$
with no $k$-rational point in its non-empty fiber in
$\mathscr{T}'$.  All we shall actually require is that there is some $t'_0 \in \mathscr{T}'_{t_0}$
with $k(t'_0) \ne k$.

Letting $X$, $T$, and $T'$ denote the analytifications of $\mathscr{X}$,
$\mathscr{T}$, and $\mathscr{T}'$, by Theorem \ref{fsep} the analytified map $X' \rightarrow X$
is an isomorphism over $X - T$ and restricts to a degree-$d$
finite \'etale covering $h:T' \rightarrow T$ over $T$.  Consider
the fiber over $t_0 \in \mathscr{T}(k) = T(k)$.
Since $T'_{t_0} = (\mathscr{T}'_{t_0})^{\rm{an}}$,
we can choose $t'_0 \in T'$
over $t_0$ with $k(t'_0) \ne k$.  Since $X' \rightarrow X$ is locally quasi-finite,
by the local structure theorem for
quasi-finite maps of rigid spaces \cite[Thm.~A.1.3]{relamp}
there are connected admissible opens $U' \subseteq X'$ around $t'_0$
and $U \subseteq X$ around $t_0$ such that
$U'$ lands in $U$ and the induced map $U' \rightarrow U$
is finite \'etale with fiber $\{t'_0\}$ over $t_0$.   Connectedness
of $U$ forces $U' \rightarrow U$ to have constant fiber-degree,
so this degree must equal $[k(t'_0):k(t_0)] > 1$.   But
the non-empty 2-dimensional admissible open $U$ in $X$ obviously
cannot be contained in the 1-dimensional analytic set $T$ in $X$, and
the fiber-degree of the finite \'etale covering
$U' \rightarrow U$ over any point of $U$ not in $T$
has to be 1 since $X' \rightarrow X$ restricts to an isomorphism over $X - T$.
This contradicts the constancy of fiber-degree over $U$, and so shows
that no such analytification $X'$ can exist.
\end{example}

\begin{remark}\label{kcounter}
The preceding example can be carried out in the category of
$k$-analytic spaces without requiring an extension of the base field
in the argument.
The reason is that on $k$-analytic spaces there are generally
many non-rational points even if $k$ is algebraically closed.
More specifically, the $k$-analytic space associated to $T$ has
a point $\xi$ with $\mathscr{H}(\xi) = \widehat{K}$, and over this point $\xi$ there
is a unique point $\xi'$ in the $k$-analytic space associated to $T'$ and
 $\mathscr{H}(\xi')$ is a separable degree-$d$ extension of $\widehat{K}$.
\end{remark}

\subsection{Finite \'etale quotients for affinoid spaces}

In \S\ref{analberk} we will take up a general study of the existence
of analytic \'etale quotients, and in \S\ref{51subsec}  we will give further examples
of the failure of existence for certain kinds of analytic \'etale quotient problems
(in both classical rigid geometry and in the category of $k$-analytic spaces).
Since Example \ref{counter} shows that the necessary
condition for analytifiability in Theorem \ref{locsep} is not sufficient
in the non-archimedean case, for a general study of
analytifiability of
algebraic spaces over non-archimedean fields we are now motivated to focus
attention on the problem of analytifying separated algebraic spaces.
The strategy in the proof that separated algebraic spaces
are analytifiable will be to show that
locally (in the rigid sense) we can describe
the quotient problem in such cases as that of forming
the quotient of an affinoid
by a {\em finite} \'etale equivalence relation.
This is difficult to carry out, for two reasons that have no counterpart
in the complex-analytic theory:
products for rigid spaces
(and $k$-analytic spaces) cannot be easily described set-theoretically,
and saturation with respect to an equivalence
relation is a problematic operation
with respect to the property of admissibility
for subsets of a rigid space.  In fact, we do not know how to carry
out the reduction to the finite \'etale case without leaving
the rigid-analytic category.

In the finite
\'etale case with affinoid spaces, the construction of
quotients goes as in algebraic geometry except that
there is the additional issue of checking that various $k$-algebras
are also $k$-affinoid:

\begin{lemma}\label{affdescend} Let $f:U \rightarrow X$ be a finite \'etale
surjective
map of rigid spaces.  The rigid space $X$ is affinoid if and only if
the rigid space $U$ is affinoid.  Moreover, if
$R' \rightrightarrows U'$ is a finite \'etale equivalence
relation on an affinoid rigid space $U'$
then the \'etale quotient $X' = U'/R'$ exists and $U' \rightarrow X'$
is a finite \'etale cover.
\end{lemma}

The \'etale hypothesis in the first part of the lemma is essential,
in contrast with a theorem of Chevalley \cite[II,~6.7.1]{ega}
in the case of finite surjections of schemes.  Indeed, in
\cite{liu} there is an example of a non-affinoid quasi-compact
separated surface (over any $k$) such that the normalization is affinoid.
The proof of Lemma \ref{affdescend} carries over verbatim to
the case of  affinoid $k$-analytic spaces that are not necessarily strictly $k$-analytic,
the key point being that if $A$ is a $k$-affinoid algebra
in the sense of Berkovich and it is endowed with a continuous
action by a finite group $G$ then the closed subalgebra $A^G$ is $k$-affinoid
and $A$ is finite and admissible as an $A^G$-module
\cite[2.1.14({\em{ii}})]{berbook}.  We will return to this
issue in \S\ref{analberk} when we address the problem of \'etale quotients of $k$-analytic spaces
(as the foundation of our approach to analytifying
separated algebraic spaces via rigid spaces).

\begin{proof}
Let $R = U \times_X U$, so the two projections
$R \rightrightarrows U$ are finite \'etale covers.
If $X$ is affinoid then certainly $U$ is affinoid, so
now assume that $U$ is affinoid.  Hence, the $U$-finite
$R$ is affinoid and we have to prove that $U/R$ is affinoid when
it exists.
More generally, we suppose that we
are given a finite \'etale equivalence relation $R' \rightrightarrows U'$
with $R'$ and $U'$ affinoid rigid spaces over $k$, and we
seek to prove that the \'etale
quotient $U'/R'$ exists as an affinoid rigid space
with $U' \rightarrow U'/R'$ a finite \'etale covering.

We have $U' = {\rm{Sp}}(A')$ for some $k$-affinoid $A'$, and
likewise $R' = {\rm{Sp}}(A'')$ for some
$k$-affinoid $A''$.  Since
the maps
$p_1, p_2:R' \rightrightarrows U'$ are finite,
the groupoid conditions may be expressed
in opposite terms using $k$-affinoid algebras with
only ordinary tensor products intervening in the description.
The resulting pair of maps of affine $k$-schemes
$\Spec(A'') \rightrightarrows \Spec(A')$ is therefore
a finite \'etale equivalence relation in the
category of $k$-schemes provided that the natural map
$$\Spec(A'') \rightarrow \Spec(A') \times_{\Spec k} \Spec(A')$$
is a monomorphism.  We claim that
it is a closed immersion.  By hypothesis, the map
$\delta:R' \rightarrow U' \times U'$ is a monomorphism between
rigid spaces, yet the first projection
$U' \times U' \rightarrow U'$ is separated and has composite
with $\delta$ that is finite, so $\delta$ is finite.
A finite monomorphism of rigid spaces is a closed immersion
(by Nakayama's lemma), so $\delta$ is a closed
immersion.  Hence, the natural map
$A' \widehat{\otimes}_k A' \rightarrow A''$ is
surjective.  The image of
$a'_1 \widehat{\otimes} a'_2$ is $p_1^{\ast}(a'_1)p_2^{\ast}(a'_2)$,
so every element of the Banach $A'$-module $A''$ (say via $p_1^{\ast}$)
is a convergent linear combination of elements in
the image of $p_2^{\ast}$.  Since $A''$ is
in fact $A'$-finite via $p_1^{\ast}$, so all $A'$-submodules
are closed,  it follows that
$p_2^{\ast}(A')$ algebraically spans $A''$ over $A'$ (via $p_1$).
Hence, the natural map $A' \otimes_k A' \rightarrow A''$
is indeed surjective, as desired.

By \cite[Exp.~V,~4.1]{sga3},
the \'etale quotient of $\Spec(A'') \rightrightarrows
\Spec(A')$ exists as an affine scheme
$\Spec(A)$ over $k$, with $\Spec(A') \rightarrow
\Spec(A)$ a finite \'etale covering (and
$A' \otimes_A A' \rightarrow A''$ an isomorphism).
If we can show that $A$ is a $k$-affinoid algebra,
then for $X' = {\rm{Sp}}(A)$ the
finite \'etale covering $U' \rightarrow X'$ yields
{\em equal} composites $R' \rightrightarrows X'$ and
the induced map $R' \rightarrow U' \times_{X'} U'$ is
an isomorphism since
$A' \otimes_A A' = A' \widehat{\otimes}_A A'$.
Thus, $X'$ would serve as the \'etale quotient $U'/R'$
in the category of rigid spaces.

To show that such an $A$ must be $k$-affinoid,
consider more generally an affine $k$-scheme
$\Spec A$
equipped with a finite \'etale covering
$\Spec(A') \rightarrow \Spec(A)$ with
$A'$ a $k$-affinoid algebra.  We claim that $A$ must be
$k$-affinoid.  Since $A'$ has only finitely many idempotents,
the same must hold for $A$, and so we may assume $\Spec(A)$ is connected.
The map $\Spec(A') \rightarrow \Spec(A)$ is a finite \'etale
covering, so each of the finitely many connected
components of $\Spec(A')$ is a finite \'etale cover
of $\Spec(A)$.  Thus, we may also assume that $\Spec(A')$ is connected.
By the theory of the \'etale fundamental group,
the connected finite \'etale covering
$\Spec(A') \rightarrow \Spec(A)$ is
dominated by a Galois finite \'etale
covering $\Spec(B) \rightarrow \Spec(A)$, say with Galois group $G$.
The faithful $G$-action on $B$ is $A$-linear, hence continuous, and
$A = B^G$.  Since $B$ is $A'$-finite, $B$ is a $k$-affinoid algebra.
Thus, by \cite[6.3.3/3]{bgr} the invariant subalgebra
$A = B^G$ is $k$-affinoid.
\end{proof}

To generalize beyond the case of finite \'etale
equivalence relations on affinoids as in Lemma
\ref{affdescend}, a fundamental issue is the possibility that
the rigid-analytic morphism $R \rightarrow U \times U$ may not be quasi-compact.
For example, if $\mathscr{X}$ is
a locally separated algebraic space
then its diagonal is a quasi-compact immersion that is
not a closed immersion if $\mathscr{X}$ is not separated, and so when working
over an \'etale chart of the algebraic space the
pullback of this diagonal morphism has analytification that is
not quasi-compact in the sense of rigid geometry when $\mathscr{X}$ is not separated.
Lack of such quasi-compactness on the rigid side
presents a difficulty because forming saturations
under the equivalence relation thereby involves the
image of a non-quasi-compact admissible open under
a flat morphism of rigid spaces, and the admissibility of such images is
difficult to control (even when the flat morphism
is quasi-compact).  This is what happens in Example \ref{counter}
if we try to use gluing to build the non-existent analytification there. We are therefore led to restrict
our attention to the analytic quotient problem when
$\delta:R \rightarrow U \times U$ is quasi-compact.
For reasons explained in Remark \ref{qcrem}, we will focus on the case
when $\delta$ is a closed immersion.

\subsection{GAGA for algebraic spaces}

We conclude this section with a  discussion of
cohomological issues related to the Tate-\'etale topology.  The GAGA theorems
aim to compare cohomology
of coherent sheaves on (analytifiable) proper algebraic spaces
and proper rigid spaces, so a basic fact that we must confront before contemplating
such theorems
is that algebraic spaces have only an \'etale topology rather than a Zariski topology
whereas rigid spaces have a Tate topology with respect to which
\'etale maps are {\em not} generally local isomorphisms.   This contrasts with the situation
over $\C$, where \'etale analytic maps are local isomorphisms.  Thus, it seems appropriate to
sketch how the GAGA formalism is to be set up for analytifiable algebraic spaces
over a non-archimedean field $k$.

 Let $\mathscr{X}$ be an analytifiable algebraic space, say with $X = \mathscr{X}^{\rm{an}}$, and let
$\mathscr{X}_{\et}$ denote the \'etale site
whose objects are schemes \'etale over $\mathscr{X}$.    For any sheaf of sets
$\mathscr{F}$ on $X_{\et}$, define the pushforward
$(f_{\mathscr{X}})_{\ast}(\mathscr{F})$ on $\mathscr{X}_{\et}$ by the formula
$((f_{\mathscr{X}})_{\ast}(\mathscr{F}))(\mathscr{U}) =
\mathscr{F}(\mathscr{U}^{\rm{an}})$.  This is a sheaf, due to Example \ref{coverex}.
It is easy to construct an exact left adjoint in the usual manner,
and this gives a map of ringed topoi $f_{\mathscr{X}}:\widetilde{X_{\et}}
 \rightarrow \widetilde{\mathscr{X}_{\et}}$
that is natural in $\mathscr{X}$.

 The pullback operation
$f_{\mathscr{X}}^{\ast}$ on sheaves
of modules is exact because  if $U \rightarrow \mathscr{U}^{\rm{an}}$
is an \'etale map from an affinoid space over $k$ to the analytification of an affine
algebraic $k$-scheme then the induced map on coordinate
rings is flat (due to the induced map on completed stalks at maximal
ideals being finite \'etale).  We call this the {\em analytification} functor
on sheaves of modules.  By exactness, this functor preserves coherence.
At the end of \S\ref{propsec} we noted that for any coherent
sheaf $\mathscr{G}$ on $\mathscr{X}_{\et}$ one can use descent theory to naturally construct
a coherent sheaf $\mathscr{G}^{\rm{an}}$ on $X_{\rm{Tate}}$ whose formation
is compatible with pullback in $\mathscr{X}$.
We claim that the associated coherent sheaf $(\mathscr{G}^{\rm{an}})_{\et}$
on $X_{\et}$ is naturally isomorphic to
$f_{\mathscr{X}}^{\ast}(\mathscr{G})$.
This follows from the commutative diagram of ringed sites
$$\xymatrix{
{{U}_{\et}} \ar[r] \ar[d] & {\mathscr{U}_{\et}} \ar[d]\\
{U_{\rm{Tate}}} \ar[r] & {\mathscr{U}_{\rm{Zar}}}}$$
for any scheme $\mathscr{U}$ \'etale over $\mathscr{X}$.
Since pullback along $X_{\et} \rightarrow X_{\rm{Tate}}$
defines an equivalence between categories of coherent
sheaves, we may therefore write $\mathscr{G}^{\rm{an}}$ to denote
$f_{\mathscr{X}}^{\ast}(\mathscr{G})$ for any $\mathscr{O}_{\mathscr{X}_{\et}}$-module
$\mathscr{G}$ without creating confusion.

\begin{example}\label{gagaex}
We can now state GAGA for algebraic spaces.
Let $h:\mathscr{X} \rightarrow \mathscr{Y}$ be a map
between analytifiable algebraic spaces, with associated
map $h^{\rm{an}}:X = \mathscr{X}^{\rm{an}} \rightarrow \mathscr{Y}^{\rm{an}} = Y$.
Since the commutative diagram of ringed topoi
$$\xymatrix{
{\widetilde{{X}_{\et}}} \ar[d] \ar[r] & {\widetilde{{\mathscr{X}}_{\et}}} \ar[d]\\
{\widetilde{{Y}_{\et}}} \ar[r] & {\widetilde{{\mathscr{Y}}_{\et}}}}$$
has exact pullback operations along the horizonal direction,
there is a natural $\delta$-functorial map of $\mathscr{O}_{Y_{\et}}$-modules
$$({\rm{R}}^j h_{\ast}(\mathscr{F}))^{\rm{an}} \rightarrow
{\rm{R}}^j h^{\rm{an}}_{\ast}(\mathscr{F}^{\rm{an}}).$$
GAGA for algebraic spaces is the assertion
that this comparison morphism is an isomorphism when $h$ is proper and $\mathscr{F}$
is coherent, from which the usual GAGA results
concerning full faithfulness on proper objects over $k$
and equivalences  between categories of coherent sheaves
on such objects follow exactly as over $\C$ (using Chow's Lemma for algebraic spaces).

To prove GAGA, first note that by Theorem \ref{propermap}, $h^{\rm{an}}$ is proper.  Also,
recall from Example \ref{topex} that via the pullback equivalence between
categories of coherent sheaves for
the Tate and Tate-\'etale topologies,  we get a $\delta$-functorial compatibility between
higher direct images for coherent sheaves in the proper case.  Consequently,
we can argue exactly as over $\C$ to reduce GAGA
for proper maps between analytifiable algebraic spaces
to the known case of schemes with the Zariski (rather than \'etale) topology
and rigid spaces with the Tate (rather than Tate-\'etale) topology.
\end{example}

\section{Analytification via $k$-analytic spaces}\label{analberk}

\subsection{Preliminary considerations}\label{prelimcon}

We are going to
now study analytification in the category of $k$-analytic spaces, and then use such spaces to
overcome admissibility problems in the rigid case. In order to make sense of this, we briefly digress to discuss
how the methods in \S\ref{algspacesec} carry over to the category of $k$-analytic spaces, endowed with their
natural \'etale topology. (As usual in the theory of $k$-analytic spaces, we allow the possibility that $k$ has
trivial absolute value.) An {\em \'etale equivalence relation} in the category of $k$-analytic spaces is a pair
of \'etale morphisms $R \rightrightarrows U$ such that the map $R \rightarrow U \times U$ (called the {\em
diagonal}) is a functorial equivalence relation; in particular, it is a monomorphism.   As one example, if
$\mathscr{R} \rightrightarrows \mathscr{U}$ is an \'etale chart for an algebraic space $\mathscr{X}$ over $k$
then the analytification functor \cite[2.6.1]{berihes} to the category of good strictly $k$-analytic spaces
yields an \'etale equivalence relation $R \rightrightarrows U$ on $k$-analytic spaces.   (By
\cite[4.10]{temkin2}, the category of strictly $k$-analytic spaces is a full subcategory of the category of
$k$-analytic spaces, so there is no ambiguity about where the morphisms $R \rightrightarrows U$ take place when
$R$ and $U$ are strictly $k$-analytic.)

\begin{definition} Let $R \rightrightarrows U$ be an \'etale equivalence relation
on $k$-analytic spaces.
A {\em quotient} of $R \rightrightarrows U$ is
a $k$-analytic space $X$
equipped with an \'etale surjection
$U \rightarrow X$ such that
the composite maps $R \rightrightarrows U \rightarrow X$
coincide and the resulting map $R \rightarrow U \times_X U$ is an isomorphism.
\end{definition}

In order to check that the quotient (when it exists) is unique up to
unique isomorphism (and in fact represents a specific sheaf functor),
we can use the usual descent theory argument
as in the case of schemes provided that
representable functors on the category of
$k$-analytic spaces
are \'etale sheaves.  This sheaf property
is true within the full subcategories of good $k$-analytic spaces
and strictly $k$-analytic spaces by  \cite[4.1.5]{berihes}, according to which
the general case holds once we prove the next result.

\begin{theorem}\label{nogood} Let $f:X' \rightarrow X$ be a finite \'etale map between $k$-analytic
spaces.  If $V' \subseteq X'$ is a quasi-compact $k$-analytic subdomain then $f(V') \subseteq X$ is
a finite union of $k$-affinoid subdomains in $X$.
In particular, $f(V')$ is a $k$-analytic domain in $X$. If $X$, $X'$, and $V'$ are strictly $k$-analytic
then so is $f(V')$.
\end{theorem}

The fiber product $R = X' \times_X X'$ is finite over $X'$ (under either projection), so
it is good (resp. strictly $k$-analytic) when $X'$ is.  Since $X = X'/R$, it will
follow from Theorem \ref{newb} below that if $X'$ is separated
then $X$ is necessarily good (resp. strictly $k$-analytic)
if $X'$ is so.

\begin{proof}  The image of $f$ is open and closed in $X$, so we may and do assume that
$f$ is surjective.   Since $X$ is locally Hausdorff and $V'$ is compact,
there is a finite collection of Hausdorff open subsets $U_1, \dots, U_n$ in $X$
that cover $f(V')$.  The open cover $\{f^{-1}(U_i)\}$ of the quasi-compact
$V'$ has a finite refinement consisting of $k$-affinoid subdomains
$V'_j \subseteq V'$, so if we can settle the  case of
a Hausdorff target then applying this to $f^{-1}(U_i) \rightarrow U_i$
and each $V'_j$ mapping into $U_i$ gives
the result for $f(V')$.
Hence, we now may and do assume that $X$ is Hausdorff, so $X'$ is also Hausdorff.

Let $W_1, \dots, W_m \subseteq X$ be a finite collection of $k$-affinoid subdomains
whose union contains $f(V')$ (with all $W_j$ strictly $k$-analytic when
$X'$, $X$, and $V$ are so).  The pullback subdomains $W'_j = f^{-1}(W_j)$
are $k$-affinoid in $X'$, and are strictly $k$-analytic when
$X'$, $X$, and $V'$ are so.  Moreover, $V' \cap f^{-1}(W_j)$ is quasi-compact
since the graph morphism $\Gamma_f:X' \rightarrow X' \times X$ is quasi-compact (as it is a base change of the diagonal morphism $\Delta_X:X \rightarrow X \times X$
that is topologically proper since $|X|$ is Hausdorff
and $|X \times X| \rightarrow |X| \times |X|$ is proper).   Hence, we may reduce to the case when
$f(V') \subseteq W$ for some $k$-affinoid subdomain $W \subseteq X$.
It is harmless
to make the base change by $W \rightarrow X$, so we can assume that $X$ and $X'$
are $k$-affinoid and even connected.  Say $X' = \mathscr{M}(A')$ and $X = \mathscr{M}(A)$.

By the theory of the \'etale fundamental group as in
the proof of Lemma \ref{affdescend}, now applied to $\Spec A' \rightarrow \Spec A$,
we may find a connected
finite \'etale cover $X'' \rightarrow X'$ that is Galois over $X$.  In particular, if $X'$
is strict then so is $X''$.  The preimage of $V'$ in $X''$ is quasi-compact
(and strict when $X'$ and $V'$ are strict), so we may assume that
$X'$ is Galois over $X$, say with Galois group $G$. The union
$W' = \cup_{g \in G} g(V')$ is a quasi-compact $k$-analytic subdomain
whose image in $X$ is the same as that of $V'$, so we can rename it as
$V'$ to get to the case when $V'$ is $G$-stable.

 For
each point $x' \in V'$ we let $G_{x'} \subseteq G$ denote the stabilizer
group of the physical point $x'$, so by the Hausdorff property of
our spaces we can find a $k$-affinoid
neighborhood $W' \subseteq V'$ around $x'$ in $V'$ that is disjoint from its
$g$-translate for each $g \in G - G_{x'}$.  Replacing $W'$ with the $k$-affinoid overlap
$\cap_{g \in G_{x'}} g(W')$ allows us to assume that $W'$ is $G_{x'}$-stable
and that the subdomains $g(W')$ for $g \in G/G_{x'}$ are pairwise disjoint.
Hence, $Y' = \coprod g(W')$ is a $G$-stable $k$-affinoid subdomain in $X'$
that is a neighborhood of $x'$ in $V'$.   By quasi-compactness of $V'$,
finitely many such subdomains $Y'_1, \dots, Y'_n$ cover $V'$.   Thus,
we can replace $V'$ with each of the $Y'_i$'s separately, so we can assume
that $V' = \mathscr{M}(B')$ is $k$-affinoid.
By \cite[2.1.14({\em{ii}})]{berbook}, the closed $k$-subalgebra
$B = {B'}^G$ is $k$-affinoid.
It is moreover a strict $k$-affinoid algebra if $V'$ is strict \cite[6.3.3]{bgr}.
The map $V' \subseteq X' \rightarrow X$ factors through
the surjection $V' = \mathscr{M}(B') \rightarrow \mathscr{M}(B)$, so
it suffices to check that the natural map $V = \mathscr{M}(B) \rightarrow \mathscr{M}(A) = X$
is a $k$-analytic subdomain.  This amounts to the property that if
$Z = \mathscr{M}(C)$ is $k$-affinoid and a morphism $h:Z \rightarrow X$ factors
through $V$ set-theoretically then it uniquely does so in the category of
$k$-analytic spaces.   It is therefore equivalent to prove that
the projection $Z \times_X V \rightarrow Z$ is an isomorphism.
But this is a map of $k$-affinoids, so it suffices to check the isomorphism assertion
after an analytic extension of the base field (an operation
which commutes with the formation of $B$ from $B'$).  We may therefore put ourselves in
the strictly $k$-analytic case (with $|k^{\times}| \ne \{1\}$), in which case the image
$f(V')$ is a $k$-analytic subdomain by Raynaud's theory,
and $V' \rightarrow f(V')$ is a finite mapping because
$V'$ is the full preimage of $f(V')$ in $X'$ (due to the $G$-stability of $V'$ in $X'$).
Hence, Lemma \ref{affdescend} gives that $f(V')$ is $k$-affinoid, and then
its coordinate ring is forced to be ${B'}^G = B$ since $X' \rightarrow X$
corresponds to $A = {A'}^G \rightarrow A'$.  That is,
the $k$-analytic subdomain $f(V') \subseteq X$ is
precisely $V$ equipped with its natural map to $X$, so $V$ is a $k$-analytic subdomain of $X$
as desired.
\end{proof}

\begin{example}\label{daff} In the setup of Theorem \ref{nogood}, if
$V' \subseteq X'$ is a quasi-compact $k$-analytic subdomain whose
two pullbacks to $X''$ coincide then it descends uniquely to
a $k$-analytic subdomain $V \subseteq X$.  Indeed,
if we let $V$ be the quasi-compact $k$-analytic subdomain $f(V') \subseteq X$
then to check that the preimage of $V$ in $X'$ is no larger than (and
hence is equal to) $V'$ it suffices
check this after base change on $X$ by geometric points of $V$.  This case is trivial.
\end{example}

By Theorem \ref{nogood}, if $R \rightrightarrows U$ is an \'etale equivalence
relation on $k$-analytic spaces and $X$ is a quotient  for this equivalence
relation in the sense
that we have defined for $k$-analytic spaces, then
$X$ represents the quotient sheaf
of sets $U/R$ on the \'etale site for
the category of $k$-analytic spaces.  Thus, such an $X$ is unique up to unique isomorphism.
We can also use descent arguments as in the classical rigid case
to run this in reverse:  if the quotient sheaf $U/R$ on
the \'etale site for the category of $k$-analytic spaces
is represented by a $k$-analytic space $X$ then
the natural map $U \rightarrow X$ is automatically an \'etale
surjection that equalizes the maps $R \rightrightarrows U$
and yields an isomorphism $R \simeq U \times_X U$.
In particular, the formation of the quotient is compatible with arbitrary
analytic extension of the base field (when the quotient exists over the initial base field).

In the $k$-analytic category, if the diagonal $R \rightarrow U \times U$ of an \'etale equivalence relation on a
locally separated $k$-analytic space $U$ is compact then it must be a closed immersion
\cite[2.2]{descent}.
This is why in Theorem \ref{berqt} we impose the requirement that $\delta$ be a closed
immersion rather than the apparently weaker condition that it be compact. We do not have an analogous such
result in the rigid-analytic case because \'etaleness and local separatedness in $k$-analytic geometry are
stronger conditions than in rigid geometry.

The arguments in \S\ref{algspacesec} carry over essentially {\em verbatim} to show that if $R \rightrightarrows
U$ arises from an \'etale chart $\mathscr{R} \rightrightarrows \mathscr{U}$ for an algebraic space $\mathscr{X}$
then whether or not an analytic quotient $X= U/R$ exists is independent of the choice of \'etale chart for
$\mathscr{X}$, and its formation (when it does exist) is Zariski-local on $\mathscr{X}$.  In particular, when
$X$ exists it is canonically independent of the chart and is  functorial in $\mathscr{X}$ in a manner that
respects the formation of fiber products and Zariski-open and Zariski-closed immersions.  We call such an $X$
(when it exists) the {\em analytification} of $\mathscr{X}$ in the sense of $k$-analytic spaces, and we say that
$\mathscr{X}$ is analytifiable (in the sense of $k$-analytic spaces); we write $\mathscr{X}^{\rm{an}}$ to denote
this $k$-analytic space if there is no possibility of confusion with respect to the analogous notion for rigid
spaces. In principle analytifiability in the sense of rigid spaces is weaker than in the sense of $k$-analytic
spaces since \'etaleness is a weaker condition in rigid geometry than in $k$-analytic geometry. It seems likely
(when $|k^{\times}| \ne \{1\}$) that analytifiability in the sense of $k$-analytic spaces implies it in the
sense of rigid spaces over $k$, but we have not considered this matter seriously because in the separated case
we will prove analytifiability in both senses (and the deduction of the rigid case from the $k$-analytic case
will use separatedness).

Since change of the base field is a straightforward operation for $k$-analytic spaces
(unlike for general rigid spaces), it is easy to see that if $K/k$ is an analytic
extension field and $\mathscr{X}$ is analytifiable in the sense of $k$-analytic spaces
then $K \otimes_k \mathscr{X}$ is analytifiable in the sense of $K$-analytic spaces
with $K \widehat{\otimes}_{k} \mathscr{X}^{\rm{an}}$ as its analytification.

\begin{theorem}\label{berklocsep} If $\mathscr{X}$ is analytifiable in the sense of
$k$-analytic spaces then $\mathscr{X}$ is locally separated.
\end{theorem}

\begin{proof}
We wish to carry over the method used to prove Theorem \ref{locsep} for rigid spaces,
so we need to recall several properties of $k$-analytic spaces that are relevant to
this method.
A separated $k$-analytic space has Hausdorff underlying topological space, and
by using rigid-analytic techniques
we see that  a dense open immersion of algebraic $k$-schemes
induces an open immersion of $k$-analytic spaces with dense image.
Also, any base change of the diagonal map $\Delta_S:S \rightarrow S \times S$
of a $k$-analytic space $S$ is
a topological embedding.  To see this, since $S$ is locally Hausdorff we may
(as in the complex-analytic case) reduce to the case when $S$ is Hausdorff.
Thus, $|S| \rightarrow |S| \times |S|$ is a closed embedding.
Letting $h:Z \rightarrow S \times S$ be any map of $k$-analytic spaces,
we want to show that the natural map $h^{\ast}(\Delta_S):Z \times_{S \times S} S \rightarrow Z$
is a topological embedding.   We will show that it is even a closed embedding.
On geometric points this map is clearly injective, so it suffices to prove
that it is topologically a proper map.  But the map from the underlying topological
space of a $k$-analytic fiber product to the fiber product of the underlying topological
spaces is always proper, so we are reduced to showing that
$$|Z| \times_{|S \times S|} |S| \rightarrow |Z|$$
is a closed embedding.   Thus, it is enough to prove that $\Delta_S$ is topologically
a closed embedding.
Since $|S \times S| \rightarrow |S| \times |S|$ is separated (even a proper surjection)
and $|S| \rightarrow |S| \times |S|$ is a closed embedding (as $|S|$ is Hausdorff),
$\Delta_S$ is indeed a closed embedding.

To use the proof of Theorem \ref{locsep} in the setting of $k$-analytic spaces,
it remains to show that a finite type map $f:\mathscr{V}
\rightarrow \mathscr{W}$ of algebraic $k$-schemes is
injective (resp. surjective) if its analytification $f^{\rm{an}}:V \rightarrow W$
(in the sense of $k$-analytic spaces) is injective (resp. surjective).
For surjectivity we use that the natural map
$W = \mathscr{W}^{\rm{an}} \rightarrow \mathscr{W}$ is surjective.
For injectivity, recall that
a map of algebraic $k$-schemes is injective if and only if it
is injective on underlying sets of closed points, and the closed
points of an algebraic $k$-scheme are functorially identified with
the set of points of the analytification with residue field
of finite degree over $k$.  Hence, we get the desired inheritance of injectivity
from $k$-analytic spaces to schemes.
\end{proof}

\subsection{Main results}\label{mainsec}

Here is the main existence result in the rigid-analytic setting; the proof
will occupy the rest of \S\ref{analberk}, and will largely
be taken up with the proof of an existence result for \'etale quotients in the
$k$-analytic category (modulo a crucial existence
theorem for quotients by free actions of finite groups, to be treated in Theorem \ref{ug}).

\begin{theorem}\label{makespace}
Assume $|k^{\times}| \ne \{1\}$.
If $\mathscr{X}$ is a separated algebraic
space over $k$ then $\mathscr{X}$ is analytifiable in
the sense of rigid spaces.  Moreover, the rigid space $\mathscr{X}^{\rm{an}}$ is separated.
\end{theorem}

Once $\mathscr{X}^{\rm{an}}$ is proved to exist, it must be separated since $\Delta_{\mathscr{X}^{\rm{an}}} =
\Delta_{\mathscr{X}}^{\rm{an}}$ is a closed immersion (as $\mathscr{X}$ is separated). For a separated algebraic
space, we will prove analytifiability in the sense of rigid spaces by deducing it from a stronger existence
theorem for \'etale quotients in the setting of $k$-analytic spaces (allowing $|k^{\times}| = \{1\}$). Consider
an \'etale chart $\mathscr{R} \rightrightarrows \mathscr{U}$ for $\mathscr{X}$.  Note that we can take
$\mathscr{U}$ to be separated. Let $U$ and $R$ be the good strictly $k$-analytic spaces associated to
$\mathscr{U}$ and $\mathscr{R}$ (so $U$ is separated when $\mathscr{U}$ is). The dictionary relating
$k$-analytic spaces and algebraic schemes  \cite[3.3.11]{berihes} ensures that $R \rightrightarrows U$ is an
\'etale equivalence relation on $U$ and that $R \rightarrow U \times U$ is a closed immersion. At the end of
\S\ref{etlocsec}, Theorem \ref{makespace} will be deduced from the following purely $k$-analytic result
(allowing $|k^{\times}| = \{1\}$).

\begin{theorem}\label{newb} Let $R \rightrightarrows U$ be an \'etale equivalence relation
on $k$-analytic spaces such that $R \rightarrow U \times U$ is
a closed immersion.   The quotient $U/R$ exists
and is a separated $k$-analytic space.
If $U$ is strictly $k$-analytic $($resp. good$)$ then so is $U/R$.
\end{theorem}

In Example \ref{nonsep} we will show that
it is insufficient in Theorem \ref{newb} to
weaken
the hypothesis on $\delta$ to compactness.
Before we proceed to global
considerations, let us first show that the existence problem for $U/R$ is local on $U$
(setting aside for now
the matter of proving separatedness of $U/R$).
To this end, suppose $U$ is covered by open subsets $\{U_i\}$ such that for $R_i = R
\times_{U \times U} (U_i \times U_i) = R \cap (U_i \times U_i)$ the quotient $X_i = U_i/R_i$ exists (with $X_i$
strictly $k$-analytic when $U_i$ is, and likewise for the property of being good); note that $R_i \rightarrow
U_i \times U_i$ is
a closed immersion.  (We do not assume
that each $X_i$ is known to be separated.)
We need to define ``overlaps'' along which we shall glue the $X_i$'s to build
a $k$-analytic quotient $U/R$. The open overlap $R_{ij} = p_1^{-1}(U_i) \cap p_2^{-1}(U_j)$ in $R$ classifies
equivalence among points of $U_i$ and $U_j$, so its open image $U_{ij}$ in $U_i$ under the \'etale morphism
$p_1:R \rightarrow U$ classifies points of $U_i$ that are equivalent to points of $U_j$.
 Let $X_{ij} \subseteq X_i$ be
the open image of $U_{ij}$, so
$p_1:R_{ij} \rightarrow X_{ij}$ is an \'etale surjection.
Geometrically, the points of $X_{ij}$ are the $R$-equivalence
classes that meet $U_i$ and $U_j$ (viewed within $X_i = U_i/R_i$).

The canonical involution $R \simeq R$ restricts to an isomorphism
$\phi_{ij}:R_{ij} \simeq R_{ji}$ such that
$\phi_{ji} = \phi_{ij}^{-1}$, and it is easy to check that
the resulting isomorphism
$R_{ij} \times R_{ij} \simeq R_{ji} \times R_{ji}$
restricts to an isomorphism
of subfunctors
$R_{ij} \times_{X_{ij}} R_{ij} \simeq R_{ji} \times_{X_{ji}} R_{ji}$.
Hence, since representable functors on the category of
$k$-analytic spaces
are \'etale sheaves (due to Theorem \ref{nogood} and \cite[4.1.5]{berihes}), the isomorphisms
$\phi_{ij}$ uniquely descend to isomorphisms
$X_{ij} \simeq X_{ji}$ between open subsets
$X_{ij} \subseteq X_i$ and $X_{ji} \subseteq X_j$.
These descended isomorphisms among opens in $X$ satisfy the triple overlap
condition, and so we can glue the $X_i$'s
along these isomorphisms to build a $k$-analytic
 space $X$.  Moreover, if $U$ is strictly $k$-analytic (resp. good) and the $U_i$'s can be
 chosen to be strictly $k$-analytic (resp. good) then
  so are all $X_i$ and hence so is the space $X$ that has an open covering by the $X_i$'s.
 The \'etale composites
$U_i \rightarrow X_i \subseteq X$
glue to define an \'etale morphism
$U \rightarrow X$ such that the two composite maps
$R \rightrightarrows U \rightarrow X$ coincide
and $R \rightarrow U \times_X U$ is an isomorphism
(as it is an \'etale monomorphism that is surjective
on geometric points).  It follows that as an \'etale sheaf
of sets on the category of
$k$-analytic spaces, $X$ represents
the sheafified quotient $U/R$.

To finish the localization argument, we now check the global property that diagonal map $X \rightarrow X \times
X$ is a closed immersion (i.e., $X$ is separated) when $X = U/R$ exists. \'Etale surjective base change by $U
\times U \rightarrow X \times X$ yields the map $R \rightarrow U \times U$ that is a closed immersion by
hypothesis. To deduce that $\Delta_X$ is a closed immersion it remains to show that the property of a
$k$-analytic morphism being a closed immersion is \'etale-local on the target.
That is, if $f:Y' \rightarrow Y$ is a map of $k$-analytic spaces and
$V \rightarrow Y$ is an \'etale cover such that the base change $F:V' \rightarrow V$ of $f$
is a closed immersion then we want to prove that $f$ is a closed immersion.
To prove this we require a straightforward descent theory for coherent sheaves
with respect to the $G$-topology $S_G$ \cite[\S1.3]{berihes} on $k$-analytic spaces $S$:

\begin{lemma}\label{cohdescend}
Let $f:S' \rightarrow S$ be a flat quasi-finite surjection of $k$-analytic spaces,
and let $p_1, p_2:S'' = S' \times_S S' \rightarrow S'$ be the canonical projections.
For any coherent sheaf $\mathscr{F}$ on $S_G$
define $\mathscr{F}'$ to be the coherent pullback $f^{\ast}(\mathscr{F})$
on $S'_G$ and define
$\varphi_{\mathscr{F}}:p_1^{\ast}(\mathscr{F}') \simeq p_2^{\ast}(\mathscr{F}')$
to be the evident isomorphism.  The functor
$\mathscr{F} \rightsquigarrow (\mathscr{F}', \varphi_{\mathscr{F}})$
from the category ${\rm{Coh}}(S_G)$ of coherent sheaves on $S_G$
to the category of pairs consisting of
an object $\mathscr{F}' \in {\rm{Coh}}(S'_G)$ equipped with a descent datum
$\varphi:p_1^{\ast}(\mathscr{F}') \simeq p_2^{\ast}(\mathscr{F}')$ relative to $f$
is an equivalence of categories.
\end{lemma}

The case of interest to use is when $f$ is \'etale.  The general notions of
quasi-finite and flat quasi-finite maps are discussed in \cite[\S3.1-\S3.2]{berihes}.

\begin{proof}
For faithfulness it suffices to work locally on $S'$ and $S$, so we can assume
$f$ is a flat finite surjection.  Once faithfulness is proved, for full faithfulness we
can work locally on $S'$ since flat quasi-finite maps are open, so we can
again reduce to the case when $f$ is a flat finite map.  Similarly, once
full faithfulness is proved the essential surjectivity holds in general if
it holds for flat finite $f$.  Hence, we can assume that $f$ is is a flat finite map.
Since the coherent sheaves are taken with respect
to the $G$-topology, for the proof of full faithfulness
we can work locally for the $G$-topology on $S$ and so we can assume
that $S$ is $k$-affinoid.  Similarly, once full faithfulness is proved we
can work locally for the $G$-topology on $S$ for the proof of essential surjectivity.
Hence, it suffices to prove the lemma when $S = \mathscr{M}(\mathscr{A})$
for a $k$-affinoid algebra $\mathscr{A}$ and $S' = \mathscr{M}(\mathscr{A}')$
for a finite and faithfully flat $\mathscr{A}$-algebra $\mathscr{A}'$ that is admissible
as an $\mathscr{A}$-module.  (The surjectivity of $\Spec(\mathscr{A}')
\rightarrow \Spec(\mathscr{A})$ follows from the surjectivity  of $f$ and the surjectivity
of $\mathscr{M}(\mathscr{B}) \rightarrow \Spec(\mathscr{B})$
for any $k$-affinoid algebra $\mathscr{B}$.
This is why the flat map of algebras $\mathscr{A} \rightarrow
\mathscr{A}'$ is faithfully flat.)   But in this special case
coherent sheaves correspond to finite modules over the coordinate ring
and completed tensor products are ordinary tensor products, so
the required result is a special case of faithfully flat descent for
quasi-coherent sheaves on schemes.
\end{proof}

To use this lemma, we first recall that there is a natural bijection between
(isomorphism classes of) closed immersions into a $k$-analytic space
$S$ and coherent ideal sheaves for the $G$-topology $S_G$ of $S$
(in the sense of \cite[\S3]{berihes}).    More specifically,
any closed immersion gives rise to such an ideal sheaf
via \cite[1.3.7]{berihes} and conversely any coherent ideal
sheaf arises from a unique closed immersion (up to unique isomorphism)
by \cite[1.3.7]{berihes} in the $k$-affinoid case
and a standard gluing procedure \cite[1.3.3({\em{b}})]{berihes}
in the general case.  (Strictly speaking there is a local finiteness
condition in this general gluing procedure,
but it is easily bypassed in the case of intended application by
local compactness considerations on the ambient space $S$.)
In this way, the closed immersion $F:V' \hookrightarrow V$
arising by base change from $f$
naturally corresponds to a coherent ideal sheaf $\mathscr{I}$ on $V_G$,
and on $(V' \times_V V')_G = ((Y' \times_Y Y') \times_V V')_G$ the pullback
coherent ideal sheaves $p_j^{\ast}(\mathscr{I})$
are equal since there is an equality of
closed immersions $p_1^{\ast}(F) = p_2^{\ast}(F)$ into $(Y' \times_Y Y' \times_Y Y') \times_Y V$.
By Lemma \ref{cohdescend}, $\mathscr{I}$ is therefore the pullback of a unique
coherent ideal sheaf $\mathscr{I}_0$ on $Y_G$.
For the associated closed immersion $\iota:Y_0 \hookrightarrow Y$ there is a
unique $V$-isomorphism $\phi:Y_0 \times_Y V \simeq V' = Y' \times_Y V$
(since $F$ and $\iota_V$ are closed immersions into $V$ with the same
associated coherent ideal sheaves).
The two pullback isomorphisms induced over $V' \times_V V'$ from $\phi$ coincide
because closed immersions have no nontrivial automorphisms,
and since representable functors in the $k$-analytic category are sheaves for the \'etale topology it
follows that $\phi$ uniquely descents to a $Y$-isomorphism $Y_0 \simeq Y'$.
Hence, $f:Y' \rightarrow Y$ is a closed immersion due to how $Y_0 \rightarrow Y$ was
constructed.  This completes the proof that the \'etale quotient $X = U/R$ is
separated when it exists (given that we are assuming that $\delta:R \rightarrow U \times U$
is a closed immersion).

The next lemma, which is an analogue of Lemma \ref{affdescend}, will be useful for analyzing properties of the
map $U \rightarrow U/R$ when the quotient has been constructed.

\begin{lemma}\label{affgood}  Let $f:X' \rightarrow X$ be a finite \'etale surjection between
$k$-analytic spaces.  If $X'$ is $k$-affinoid then so is $X$, and if in addition $X'$ is strictly $k$-analytic
then so is $X$.
\end{lemma}

\begin{proof}
Since $X'' = X' \times_X X'$ is finite over $X'$ under either projection, it is $k$-affinoid (and strict when
$X'$ is so). Also, the map $X'' \rightarrow X' \times X'$ between $k$-affinoid spaces is a closed immersion
because a finite monomorphism between $k$-analytic spaces is a closed immersion (as we may check after first
using analytic
extension of the base field to reduce to the strict case; the monomorphism property is preserved by such
extension since it is equivalent to the relative diagonal map being an isomorphism). The method of proof of
Lemma \ref{affdescend} therefore carries over (using \cite[2.1.14({\em{i}})]{berbook} to replace
\cite[6.3.3]{bgr}) to construct a $k$-affinoid quotient for the finite \'etale equivalence relation $X''
\rightrightarrows X'$, and this quotient is (by construction) even strict when $X'$ is strict. But $X$ is also
such a quotient, so it must be $k$-affinoid.
\end{proof}

\begin{lemma}
\label{finitelem} Assume that $R\rightrightarrows U$ is a finite \'etale
equivalence relation on $k$-analytic spaces
such that the quotient $X=U/R$ exists.  The map
$U \rightarrow U/R$ must be finite, and if $U$ is Hausdorff then
$X$ is Hausdorff.  Moreover, if $U$ is $k$-affinoid then so is $X$.
\end{lemma}

\begin{proof}
The base change of $\pi:U\to X$ by the \'etale covering $\pi:U \rightarrow U/R$ is a finite map (it is a
projection $R \rightarrow U$), so to prove that $\pi$ is finite we just have to show that if a map $h:Y'
\rightarrow Y$ between $k$-analytic spaces
becomes finite after an \'etale surjective base change on $Y$ then it is finite.   We can work
locally on $Y$, so since \'etale maps are open and are finite locally on the source we can assume that there is
a finite \'etale cover $Z \rightarrow Y$ such that $Y' \times_Y Z \rightarrow Z$ is finite. To prove finiteness
of $h$ from this we can easily reduce to the case when $Y$ is $k$-affinoid, so $Z$ is $k$-affinoid and hence $Y'
\times_Y Z$ is $k$-affinoid. The map $Y' \times_Y Z \rightarrow Y'$ is a finite \'etale cover with $k$-affinoid
source, so $Y'$ is $k$-affinoid by Lemma \ref{affgood}. Thus, $h:Y' \rightarrow Y$ is a map between
$k$-affinoids which becomes finite after the finite \'etale base change by $Z \rightarrow Y$.  The desired
finiteness of $h$ is therefore clear.

With $\pi$ now known to be finite \'etale (and surjective)
in general, if $U$ is $k$-affinoid then Lemma \ref{affgood}
ensures that $X$ must be $k$-affinoid. To see that $X$ must be Hausdorff
when $U$ is Hausdorff, the finite surjection $\pi \times
\pi:U \times U \rightarrow X \times X$ induces a closed map on topological spaces, so properness and
surjectivity of $|T_1 \times T_2| \rightarrow |T_1| \times |T_2|$ for $k$-analytic spaces $T_1$ and $T_2$
implies that $|U| \times |U| \rightarrow |X| \times |X|$ is closed. The diagonal $|U| \subseteq |U| \times
|U|$ is closed since $U$ is Hausdorff, so we conclude that $|X|$ has closed diagonal image in
$|X| \times
|X|$ as desired. That is, $X$ must be Hausdorff when it exists.
\end{proof}

\subsection{\'Etale localization and reduction to group quotients}\label{etlocsec}

Now we return to the global construction problem for $p_1,p_2:R \rightrightarrows U$
as in Theorem \ref{newb}. The aim of this section is
to reduce the problem to the particular case when $p:R \rightrightarrows U$
is induced by a free right action of a finite group
$G$ on $U$; i.e.,
 $R = U\times G$ with $p_1$ the canonical projection and $p_2(u,g) = u.g$.
  We emphasize that the hypothesis that $\delta:R \rightarrow U \times U$
 is a closed immersion will be preserved under this reduction step.
 The existence result for $U/R$ in this special case
 is proved in Theorem \ref{ug}.

  We have already seen that it suffices to work
  locally on $U$ to solve the existence problem for $U/R$.  By ``work locally''
we mean that we work with opens
$V$ that cover $U$ and the  \'etale equivalence relation $R_V = p_1^{-1}(V) \cap p_2^{-1}(V) =
\delta^{-1}(V \times V)$ on $V$ (for $\delta:R \rightarrow U \times U$).
Note that the map $R_V \rightarrow V \times V$ is still a closed immersion since
it is a base change of the map $\delta$ that we assume is a closed immersion.
 Localizing in this way  does not lose the property of the
new $U$ being strictly $k$-analytic or good when
the original $U$ is so.

It is possible to first use topological arguments (especially compactness and connectedness considerations) to
reduce the problem to the case when the maps  $R \rightrightarrows U$ are finite \'etale, and to then use \'etale
localization to split the equivalence relation, thus passing to the group action case. However, it turns out
that a shorter way to the same goal is to first apply \'etale localization
to split the equivalence relation
\'etale-locally around points of $U$ and to then use compactness and connectedness
considerations.
Given an \'etale morphism $f:U'\to U$ we define
$$R' = R \times_{U \times U} (U' \times U') = U' \times_{U,p_1} R \times_{p_2,U} U',$$
obtaining an \'etale
equivalence relation $p'_1, p'_2: R' \rightrightarrows U'$ induced from $R$.
Beware that even if the maps $R \rightrightarrows U$ are quasi-compact,
the maps $R' \rightrightarrows U'$ may fail to be quasi-compact.
(In the intended applications to algebraic spaces, such quasi-compactness
properties for the projections $R \rightrightarrows U$ are often not satisfied
even when the algebraic space is separated.)

\begin{lemma}  Let $R \rightrightarrows U$ be
an \'etale equivalence relation on $k$-analytic spaces
and let $f:U'\to U$ be an
\'etale surjection. The quotient $X=U/R$ exists if and only if the quotient $X'=U'/R'$
exists, and then $X'\simeq X$.
\end{lemma}

\begin{proof}
We will need the following set-theoretic analogue of the lemma: if $\calR\rightrightarrows\calU$ is
an equivalence relation on a set $\calU$ and there is given a surjective
map of sets $\calU'\twoheadrightarrow \calU$ then for the induced equivalence relation
$$\calR'=\calR \times_{\calU
\times\calU} (\calU' \times \calU')\rightrightarrows\calU'$$ on $\calU'$ we have that
the natural map $\calU'/\calR'\to\calU/\calR$ is bijective.
To apply this, we will work with \'etale sheaves of sets on the category of $k$-analytic spaces;
examples of such sheaves are representable functors $\underline{Z} =
\Hom( \cdot, Z)$ for $k$-analytic spaces $Z$ (as we indicated
above Theorem \ref{nogood}).   Consider the equivalence relations
$\underline{R} \rightrightarrows \underline{U}$
and $\underline{R}' \rightrightarrows \underline{U}'$, and let
$\calX$ and $\calX'$ be the corresponding quotients in the category of \'etale sheaves of sets.
It suffices to prove that $\calX'\simeq\calX$. By surjectivity of $f$,
the corresponding morphism of sheaves
$\underline{U}' \rightarrow \underline{U}$ is surjective.
In particular, this latter map of sheaves induces surjections on stalks at geometric points,
so by the above
set-theoretic result  we conclude that the natural map $\calX' \rightarrow \calX$
induces a bijection on geometric stalks.   Hence,
by \cite[4.2.3]{berihes} this natural map is an isomorphism of sheaves of sets.
\end{proof}

\begin{corollary}
If $R \rightrightarrows U$ is an \'etale equivalence
relation on $k$-analytic spaces then
the quotient $U/R$ exists if and only if for any point $u\in U$ there exists an \'etale morphism $U'\to U$ whose
image contains $u$ and such that the quotient $U'/R'$ exists, where $R' = R \times_{U \times U} (U' \times U')$.
\end{corollary}

Now we are in position to prove Theorem \ref{newb}, assuming that it is true in the case of a free action by a
finite group (with action map $U \times G \rightarrow U \times U$ a closed immersion!), a case that we will
settle in Theorem \ref{ug}.  We just have to prove the following lemma.

\begin{lemma}
Let $R \rightrightarrows U$ be an \'etale equivalence relation on $k$-analytic spaces, and assume that $R
\rightarrow U \times U$ is a closed immersion. For any point $u\in U$ there exists an \'etale morphism $W\to U$
whose image contains $u$ and such that the induced equivalence relation $R_W\rightrightarrows W$ is split; i.e.,
induced by a free right action of a finite group $G$ on $W$.
\end{lemma}

Note that if $U$ is strictly analytic $($resp. good$)$ then so is $W$.

\begin{proof}   By working locally on the topological
space of $U$ we can assume that $U$ is Hausdorff (so $U \times U$ and $R$ are Hausdorff, as $|U \times U|
\rightarrow |U| \times |U|$ is topologically proper and $R \rightarrow U \times U$ is monic and hence
separated).

To construct $W$ and $G$ we first want to split $R$ over $u$.  By hypothesis the diagonal $R \rightarrow U
\times U$ is a topologically proper map, so for any compact analytic domain $K \subseteq U$ that is a
neighborhood of $u$ in $U$ we have that $R_K = R \cap (K \times K)$ is compact.   But the closed subset
$p_1^{-1}(u) \subseteq R$ is \'etale over $\mathscr{H}(u)$ and hence is discrete, so it has finite overlap with
$R_K$.   Thus, $u$ has an open neighborhood ${\rm{int}}_U(K)$ in $U$ that contains only finitely many points
which are $R$-equivalent to $u$, so by localizing $U$ and using the Hausdorff property we can arrange that $u$
is not $R$-equivalent to any other points of $U$.  That is, now we have $p_1^{-1}(u)=p_2^{-1}(u)$ and this
common set is finite. For each $r$ in this finite set, the pullback maps $\calH(u) \rightrightarrows \calH(r)$
with respect to $p_1^{\ast}$ and $p_2^{\ast}$ are finite separable. Hence, we can choose a finite Galois
extension $\calH/\calH(u)$ that splits each $\calH(r)$ with respect to both of its $\mathscr{H}(u)$-structures
(via $p_1$ and $p_2$). By \cite[3.4.1]{berihes} we can find an \'etale morphism $U'\to U$ so that $u$ has a
single preimage $u'\in U'$ and $\calH(u')$ is $\calH(u)$-isomorphic to $\calH$. Shrinking $U'$ around $u'$ also
allows us to suppose that $U' \rightarrow U$ is separated (even finite over an open neighborhood of $u$ in $U$).
In particular, $U'$ is Hausdorff, so $R'$ is Hausdorff.

Since $u'$ is the only point in $U'$ over $u$ and no point in $U - \{u\}$ is $R$-equivalent to $u$,
no point in $U' - \{u'\}$ is $R'$-equivalent to $u'$.
That is, ${p'_1}^{-1}(u') = {p'_2}^{-1}(u')$ as (finite) subsets of $R'$.
Let $\{r'_1,\dots,r'_n\}$ be an enumeration of this set, so
the natural pullback maps $\calH = \calH(u') \rightrightarrows \calH(r'_j)$ via
$p'_1$ and $p'_2$ are isomorphisms for each $j$, due to the formula
$R' = U' \times_{U,p_1} R \times_{p_2,U} U'$ and the way that $\calH(u')/\calH(u)$ was chosen.
It follows from these residual isomorphisms and  \cite[3.4.1]{berihes}
that the \'etale maps  $p'_1$ and $p'_2$ are
local isomorphisms near each $r'_j$; i.e., for $i = 1, 2$
there exist open neighborhoods $U''_i$ of $u'$ and open neighborhoods
$R'_{ij}$ of $r'_j$ for each $j$ such that $p'_i$ induces isomorphisms $R'_{ij}\to U''_i$ for each $j$.
Since our original problem only depends on the \'etale neighborhood of $(U,u)$,
we can replace $R \rightrightarrows U$ and $u$ with $R' \rightrightarrows U'$ and $u'$
so reduce to the case when $p_1^{-1}(u) = p_2^{-1}(u) = \{r_1,\dots,r_n\}$
and there are open subspaces $U_i \subseteq U$ around $u$
such that $p_i^{-1}(U_i)$ contains an open neighborhood $R_{ij}$ around
$r_j$ mapping isomorphically onto $U_i$ under $p_i$.
However, there may be overlaps among the $R_{ij}$ for a fixed $i$
and $p_i^{-1}(U_i)$ may contain points not in any $R_{ij}$.

Choose a compact $k$-analytic domain
$\oU\subseteq U$ that is a neighborhood of $u$ in $U$,
 and let $\oR = \oU \times_{U,p_1} R \times_{p_2,U} \oU$.
 Note that $\oR$ is compact since $R \rightarrow U \times U$ is compact.
 Thus, $p:R \rightrightarrows U$ induces an
  equivalence relation $\op:\oR\rightrightarrows\oU$ that may not be \'etale (though it
is rig-etale in the strictly analytic case).  Since $\op$ is induced by $p$,
$\oR=p_1^{-1}(\oU)\cap p_2^{-1}(\oU)$ is a compact neighborhood of
$p_1^{-1}(u) = p_2^{-1}(u)$ in $R$.   Make an initial choice
of $U_1$ and let $\overline{R}_{1j} \subseteq \oR \cap p_1^{-1}(U_1)$
be an open neighborhood around $r_j$ in $R$ that is small enough so that
it maps isomorphically onto an open subspace of $U_1$ and
such that the $\oR_{1j}$'s are pairwise disjoint.
The image $p_1(\oR - (\cup \oR_{1j}))$ is  compact in $\oU$ and does not contain $u$, so
for any open subspace $U' \subseteq U_1$
around $u$ with $U' \subseteq \oU - p_1(\oR - \cup \oR_{1j})$
we have that $p_1^{-1}(U')$ contains pairwise disjoint open neighborhoods
$R'_{1j}$ around the $r_j$'s such that each $R'_{1j}$ maps isomorphically onto
$U'$ under $p_1$ and $p_2^{-1}(U')$ meets $p_1^{-1}(U')$ inside
of $\oR$.  Thus, $p_1^{-1}(U')$ is the disjoint union of the open subspaces
$R'_{1j}$ and $p_1^{-1}(U') \cap (R - \oR)$.
We can likewise choose $U_2$ and $\oR_{2j}$ adapted to $p_2$,
so by choosing
$$U' \subseteq \oU - (p_1(\oR - \cup \oR_{1j}) \cup p_2(\oR - \cup \oR_{2j}))$$
we also have that $p_2^{-1}(U')$ is the disjoint union of open subspaces
$R'_{2j}$ and $p_2^{-1}(U') \cap (R - \oR)$ with $r_j \in R'_{2j}$ and
$p_2:R'_{2j} \simeq U'$.

Let $R' = p_1^{-1}(U') \cap p_2^{-1}(U') \subseteq \oR$ be the \'etale
equivalence relation on $U'$ induced from $R$.   Thus,
the excess  open sets $p_i^{-1}(U') \cap (R - \oR)$ are disjoint from
$R'$, so $R'$ is the disjoint union of the overlaps $R'_{ij} \cap R'_{ij'}$
for $1 \le j, j' \le n$.   The initial choices of $\oR_{ij}$ above (prior to
the choice of $U'$) should be made so that not only are these
pairwise disjoint for fixed $i$ but also so that $\oR_{1j} \cap \oR_{2j'} = \emptyset$
for $j \ne j'$.  (This can be done since $r_j \ne r_{j'}$ and $R$ is Hausdorff.)
Thus, $R'_{ij} \cap R'_{ij'} = \emptyset$ for $j \ne j'$, so
$R'$ is the disjoint union of the $n$ open subspaces
$R'_j=R'_{1j}\cap R'_{2j}$ with $r_j \in R'_j$.
Each $p_i$ induces an open immersion $p'_{ij}:R'_j\to U'$ whose image is a
neighborhood of $u$.
We will prove that a careful choice of $U'$ leads to a split equivalence relation
$R'$, so in particular the $p'_{ij}$'s are isomorphisms in such cases.

Let $V$ be the connected component of $\cap_{i,j}p'_{ij}(R'_j)$ that contains $u$, so it is an open neighborhood
of $u$ in $U'$ and applying $p'_{2j}$ to ${p'_{1j}}^{-1}(V)$  for each
$1 \le j \le n$ yields open
immersions $$\phi_j:V \simeq {p'_{1j}}^{-1}(V)  \subseteq R'_j \to U'$$ that fix $u$.
Consider the equivalence relation of sets
$$\calP:\calR=\Hom(V,R')\rightrightarrows\calU=\Hom(V,U').$$
Since $V$ is connected and $R' = \coprod R'_j$ with $p'_{1j}:R'_j \rightarrow U'$ an open immersion that hits $u
\in V$, there exist exactly $n$ ways to lift the canonical open immersion $i_V:V\to U'$ with respect to $p'_1$
to a morphism $V\to R'$ (just liftings with respect to $p'_{1j}$). Therefore, $\{\phi_j\}_{1\le j\le n}$ is the
set of all elements of $\calU$ that are $\calP$-equivalent to $i_V$. In particular, for any pair $1\le a,b\le n$
the morphisms $\phi_a$ and $\phi_b$ are $\calP$-equivalent, and therefore (again using the connectedness of $V$)
the morphism $\phi_b\circ\phi_a^{-1}:\phi_a(V)\to\phi_b(V)$ is induced from $p'_{2c}\circ {p'_{1c}}^{-1}$ for
some $1\le c\le n$. Let $\psi_j$ denote the automorphism of the germ $(V,u)$ induced from $\phi_j$, so the above
argument implies that $\psi_b\circ\psi_a^{-1}=\psi_c$. Since the set $G=\{\psi_j\}_{1\le j\le n}$ includes the
identity (which corresponds to the identity point $r_{j_0}$ over $u$ in the identity part of the equivalence
relation), it follows that $G$ is actually a group of automorphisms of $(V,u)$.  Note that composition of these
automorphisms corresponds to a right action on this germ due to how each $\phi_j$ is defined in the form ``$p_2
\circ p_1^{-1}$''.

We identify the index set $\{1,\dots ,n\}$ for the $R'_j$'s with $G$
via $j \mapsto \psi_j$.
Since $G$ is a
finite group of automorphisms of $(V,u)$, there exists an open neighborhood $V'$ of $u$
in $V$ such that for
any $g \in G$ we have $g(V') \subseteq V$. It is clear that the open subspace
$W=\cap_{g\in G} g(V') \subseteq V'$ is taken to itself by each $g \in G$,
so $G$ thereby acts on $W$.  (This is a right action.)
It follows that the equivalence relation $R'$ induces on
$W$ the split equivalence relation corresponding to
this right action of $G$. Thus, $W$ and $G$ are as required.
\end{proof}

Granting Theorem \ref{newb} (whose proof rests on Theorem \ref{ug} that is proved below), we return to the
\'etale equivalence relation $\mathscr{R}^{\rm{an}} \rightrightarrows \mathscr{U}^{\rm{an}}$ on rigid spaces for
proving Theorem \ref{makespace} for an algebraic space $\mathscr{X}$ (with $|k^{\times}| \ne \{1\}$).  By
Corollary \ref{opencheck} we may assume that $\mathscr{X}$ is quasi-compact, so we can and do take $\mathscr{U}$
to be affine.   This forces $\mathscr{R}$ to be quasi-compact and separated since $\mathscr{R} \rightarrow
\mathscr{U} \times \mathscr{U}$ is a quasi-compact monomorphism. Hence, the associated $k$-analytic spaces $U$
and $R$ are paracompact.  By using analytification with values in the category of good strictly $k$-analytic
spaces, we conclude from Theorem \ref{newb} that $\mathscr{X}$ admits an analytification $X$ in the sense of
$k$-analytic spaces (as opposed to the sense of rigid spaces, which is what we want), and $X$ is separated,
good, and strictly $k$-analytic. Let us check that $X$ is also paracompact. The Hausdorff space $U$ is covered
by a rising sequence $\{K_j\}_{j \ge 1}$ of compact subsets such that $K_j \subseteq {\rm{int}}_U(K_{j+1})$ for
all $j$ because $U$ arises from an affine $\mathscr{U}$. (Choose a closed immersion of $\mathscr{U}$ into an
affine space over $k$ and intersect $U$ with closed balls of increasing radius centered at the origin in the
analytified affine space.) The images of the $K_j$'s under the open surjective map $U \rightarrow X$ are
therefore an analogous sequence of compact subsets in the Hausdorff space $X$, so $X$ is indeed paracompact.

It follows that under the equivalence of categories in \cite[1.6.1]{berihes} there
is a quasi-separated rigid space $X_0$ uniquely associated to $X$,
and the  \'etale surjective map $U \rightarrow X$ (which is in
the full subcategory of strictly $k$-analytic spaces) arises from a unique morphism
$\mathscr{U}^{\rm{an}} \rightarrow X_0$ that is necessarily \'etale (as we may check
using complete local rings) and surjective (due to $k$-strictness).
The two maps $\mathscr{R}^{\rm{an}} \rightrightarrows \mathscr{U}^{\rm{an}}$
are equalized by the map $\mathscr{U}^{\rm{an}} \rightarrow X_0$
because this equality holds back in the $k$-analytic setting (since
$X = U/R$).    The resulting map
$\mathscr{R}^{\rm{an}} \rightarrow \mathscr{U}^{\rm{an}} \times_{X_0} \mathscr{U}^{\rm{an}}$
is likewise an isomorphism due to
the isomorphism $R \simeq U \times_{X} U$ and
the functoriality and compatibility with fiber products for
the functor from paracompact strictly $k$-analytic spaces to quasi-separated rigid spaces.

Hence, by Example \ref{etqtex}
we see that $X_0$ represents $\mathscr{U}^{\rm{an}}/\mathscr{R}^{\rm{an}}$ as
long as the map $\mathscr{U}^{\rm{an}} \rightarrow X_0$ admits
local \'etale quasi-sections.  Since $X$ is paracompact and
Hausdorff, by \cite[1.6.1]{berihes} the rigid space $X_0$ has an admissible covering
arising from a locally finite collection of strictly $k$-analytic affinoid subdomains $D$ that cover
$X$.   The \'etale surjection $U \rightarrow X$
gives an \'etale cover $U \times_X D \rightarrow D$, so by quasi-compactness of $D$
and the local finiteness of \'etale maps we get a finite collection of
strictly $k$-analytic subdomains $Y_i \subseteq D$ that cover $D$ and over which there is
a finite \'etale cover $Y'_i \rightarrow Y_i$ that maps into $U \times_X Y_i$ over $Y_i$.
The rigid space associated to $\coprod Y'_i$ is then a quasi-compact \'etale cover of
$D_0 \subseteq X_0$ over which $\mathscr{U}^{\rm{an}} \rightarrow X_0$
acquires a section.

\section{Group quotients}\label{birsec}

\subsection{Existence result and counterexamples}\label{51subsec}

We shall now use the theory of reduction of germs as developed in
\cite{temkin1} for strictly analytic spaces
and in \cite{temkin2} in general to prove an existence
result for group-action quotients to which the proof of Theorem \ref{newb} was reduced in
the previous section.

\begin{theorem}\label{ug} Let $G$ be a finite group equipped with a free right action on a
$k$-analytic space $U$.  Assume that the action map $\alpha:U \times G \rightarrow U \times U$
defined by $(u,g) \mapsto (u, u.g)$ is a closed immersion.
The quotient $U/G$ exists as a separated $k$-analytic space. If $U$ is
strictly $k$-analytic $($resp. good, resp. $k$-affinoid$)$ then so is $U/G$.
\end{theorem}

The inheritance of strict
$k$-analyticity (resp. goodness) from that of $U$ will follow from our method of construction of $U/G$.
Note also that the action map $\alpha:U \times G \rightarrow U \times U$ is
the diagonal of an \'etale equivalence relation since the action is free, and since
it is a closed immersion the space $U$ must be separated
(since $\Delta_U$ is the restriction of $\alpha$ to the open
and closed subspace $U \times \{1\}$ in $U \times G$).
The separatedness of $U/G$ (once it exists) follows
by the same \'etale descent argument as we used for
a general quotient $X = U/R$ in the discussion following Theorem \ref{newb}.
Observe that if we assume $U$ is Hausdorff then
$\alpha$ is a compact topological map (since $|U \times U| \rightarrow |U| \times |U|$ is compact)
but it is not a closed immersion if $U$ is not separated,
and Example \ref{ex54} shows that $U/G$ can fail to exist in such cases.

\begin{example}\label{ugdomain}
Let us see how the formation of $U/G$ in Theorem \ref{ug} interacts
with passage to $G$-stable $k$-analytic subdomains in $U$.  In the setup of Theorem \ref{ug}, assume that $U/G$ exists and
let $\pi:U \rightarrow U/G$ be the quotient map, which must be finite \'etale
(by Lemma \ref{finitelem}).   Let $N \subseteq U$ be a $G$-stable quasi-compact
$k$-analytic subdomain,
so $\pi(N) \subseteq U/G$ is a $k$-analytic domain by Theorem
\ref{nogood} and it is good (resp. strictly $k$-analytic) if $U$, $N$, and $U/G$ are
good (resp. strictly $k$-analytic).  The inclusion of analytic domains $N \rightarrow
\pi^{-1}(\pi(N))$ in $U$ is bijective on geometric points due to $G$-stability of $N$,
so it is an isomorphism.  Thus, the map $\pi:N \rightarrow
\pi(N)$ serves as a quotient $N/G$.  Likewise,
if $N$, $U$, and $U/G$ are good (resp. strictly $k$-analytic) then so is $N/G$.
Note in particular that the natural map $N/G \rightarrow U/G$
identifies $N/G$ with a $k$-analytic domain in $U/G$.
\end{example}

Before we prove Theorem \ref{ug}, we give an example that shows that the closed immersion hypothesis on $\alpha$
in Theorem \ref{ug} cannot be replaced with a compactness hypothesis for the purposes of ensuring the existence
of $U/G$.

\begin{example}\label{nonsep}
We give an example of a 1-dimensional compact Hausdorff (and non-separated)
strictly $k$-analytic space $U$ equipped
with a free action of the group $G = \Z/2\Z$ such that $U/G$ does not exist,
assuming that there is a separable quadratic extension $k'/k$ that is ramified
in the sense that $k \rightarrow k'$ induces an isomorphism on residue fields (so in particular,
$|k^{\times}| \ne \{1\}$).   We also give
an analogous such non-existence result in the rigid-analytic category
(with quasi-compact quasi-separated rigid spaces that are not separated).
In Example \ref{ex54} this will be adapted to work in both
the rigid-analytic and $k$-analytic categories with an arbitrary $k$ (perhaps algebraically
closed or, in the $k$-analytic case, with trivial absolute value) using 2-dimensional
compact Hausdorff strictly $k$-analytic spaces.

Let $U$ be the strictly $k'$-analytic gluing of the closed unit ball $\calM(k' \langle t \rangle)$ to itself
along the identity map on the $k'$-affinoid subdomain $\{|t| = 1\}$.   Let $B_1$ and $B_2$ be the two canonical
copies of $\calM(k' \langle t \rangle)$  in $U$, so $B_1 \cap B_2$ is the affinoid space $\{|t| = 1\}$ over $k'$
whose (diagonal) map into $B_1 \times B_2$ is not finite (since this map of strict affinoids over $k'$ has
reduction over $\widetilde{k'}$ that is not finite:  it is the diagonal inclusion of $\mathbf{G}_m$ into
$\mathbf{A}^2_{\widetilde{k'}}$).  In particular, this map is not a closed immersion, so $U$ is not separated.
However, $U$ is compact Hausdorff since it is a topological gluing of compact Hausdorff spaces along a closed
subset. Under the equivalence of categories  between compact Hausdorff strictly $k'$-analytic spaces and
quasi-compact quasi-separated rigid spaces over $k'$, $U$ corresponds to the non-separated gluing $U_0$ of two
copies of $\mathbf{B}^1_{k'} = \Sp(k'\langle t \rangle)$ along the admissible open $\Sp(k' \langle t,
1/t\rangle)$ via the identity map. Now view this gluing in the $k$-analytic category. Let the nontrivial element
in $G = \Z/2\Z$ act on $U$ over $k$ by swapping $B_1$ and $B_2$ and then acting via the nontrivial element of
${\rm{Gal}}(k'/k)$ on each $B_j$.  This is seen to be a free action by computing on geometric points.  (The
action would not be free if we used only the swap.)

There is an analogous such action of $G$ on the rigid space $U_0$
over $k$.  Let us check that
non-existence of $U/G$ as a $k$-analytic space
implies non-existence of $U_0/G$ as a rigid space.  Assume that $U_0/G$ exists.
Since $U_0$ is quasi-compact and quasi-separated, the natural map
$U_0 \times G \rightarrow U_0 \times U_0$ is quasi-compact.  Thus,
$U_0/G$ is quasi-separated, so any set-theoretic union of finitely many
affinoid opens in $U_0/G$ is a quasi-compact
admissible open subspace of $U_0/G$.  It follows that $U_0/G$ must be quasi-compact.
In particular, $U_0/G$ corresponds to a compact Hausdorff $k$-analytic space.
Since the \'etale quotient map $\pi:U_0 \rightarrow U_0/G$ has local {\em fpqc}
quasi-sections and $U_0 \times G \simeq U_0 \times_{U_0/G} U_0$,
by \cite[Thm.~4.2.7]{relamp} the map $\pi$ must be finite \'etale (of degree 2).
The map of $k$-analytic spaces $U \rightarrow \overline{U}$ corresponding to $\pi$
is therefore finite \'etale (of degree 2), $G$-invariant, and gives rise to
a map $U \times G \rightarrow U \times_{\overline{U}} U$ that is an isomorphism.
Hence, $\overline{U}$ serves as a quotient $U/G$ in the category of $k$-analytic spaces.
Once it is proved that $U/G$ does not exist, it therefore follows that $U_0/G$ does not exist.

To see that $U/G$ does not exist as a $k$-analytic space, assume
to the contrary that such a quotient does exist.   Let
$\xi \in B_1 \cap B_2$ be the common Gauss point in the gluing $U$ of the
two unit disks $B_j$ over $k'$,
and let $\xi' = \pi(\xi) \in U/G$.  Since $U/G$ is a $k$-analytic space,
there exist finitely many $k$-affinoid domains $V'_1, \dots, V'_n$
in $U/G$ containing $\xi'$ with $\cup V'_j$ a neighborhood of
$\xi'$ in $U/G$.  The preimage $V_j = \pi^{-1}(V'_j)$ is
a $G$-stable $k$-affinoid in $U$ containing $\xi$ with $\cup V_j$ a neighborhood of
$\xi$ in $U$.

Choose some $j$ and let $V = V_j$. Since $G$ physically fixes the point $\xi \in U$, the action of $G$ on $U$
induces a $G$-action on the strictly $k$-affinoid germ $(V,\xi)$.  Let $\widetilde{V}_{\xi}$ be the reduction of
$(V,\xi)$ in the sense of \cite{temkin1}; this is a birational space over the residue field $\widetilde{k}$ (not
to be confused with the theory of reduction of germs and birational spaces over $\R^{\times}_{>0}$-graded fields
as developed in \cite{temkin2}). This reduction is a separated birational space over $\widetilde{k}$ since $V$
is $k$-affinoid, and it inherits a canonical $G$-action from the one on $(V,\xi)$. By separatedness of the
birational space, the $G$-action on this reduction is uniquely determined by its effect on the residue field
$\widetilde{\mathscr{H}(\xi)}$.   But this latter residue field is identified with $\widetilde{k'}(t)$
compatibly with the natural action of $G = {\rm{Gal}}(k'/k)$, so since $\widetilde{k} = \widetilde{k'}$ (by
hypothesis) we see that this $G$-action is trivial.  Hence, $G$ acts trivially on $\widetilde{V}_{\xi}$. Letting
$Y'$ (resp. $Y''$) denote the $k$-analytic domain $V \cap B_1$ (resp. $V \cap B_2$) in $U$, $(V,\xi)$ is covered
by $(Y', \xi)$ and $(Y'', \xi)$, so by \cite[2.3({\em{ii}})]{temkin1} the reductions $\widetilde{Y'}_{\xi}$ and
$\widetilde{Y''}_{\xi}$ are an open cover of $\widetilde{V}_{\xi}$ and are swapped by the $G$-action. But we
have seen that this $G$-action is trival, so we obtain the equality $\widetilde{Y'}_{\xi} =
\widetilde{V}_{\xi}$.  Thus, by \cite[2.4]{temkin1} it follows that $(Y', \xi) = (V,\xi)$ as germs at $\xi$.
Likewise, $(Y'', \xi) = (V, \xi)$, so $(Y' \cap Y'', \xi) = (V, \xi)$.  Hence, for a sufficiently small open $W
\subseteq U$ around $\xi$ we have $V \cap W \subseteq Y' \cap Y'' \subseteq B_1 \cap B_2$.  But $V = V_j$ for an
arbitrary choice among the finitely many $j$'s, so for any sufficiently small open $W$ around $\xi$ in $U$ we
have that $W \cap (\cup V_j) \subseteq B_1 \cap B_2$.  We can take $W$ so small that it is contained in the
neighborhood $\cup V_j$ of $\xi$ in $U$, so $W \subseteq B_1 \cap B_2$. Hence, $B_1 \cap B_2$ is a neighborhood
of $\xi$ in $U$, so $\calM(k' \langle t, 1/t\rangle)$ is a (closed) neighborhood of the Gauss point in
$\calM(k'\langle t \rangle)$.  This is a contradiction, since the closure in $\calM(k' \langle t \rangle)$ of
the open residue disc around any $k'$-point (such as the origin) contains the Gauss point.  Thus, $U/G$ does not
exist.
\end{example}

\begin{example}\label{ex54}  The preceding example can be adapted
in both the rigid-analytic and $k$-analytic categories to give similar
non-existence examples in the 2-dimensional case without restriction
on $k$ (e.g., $k$ may be algebraically closed or, in the $k$-analytic case, have trivial absolute
value).   Rather than work with a ramified
quadratic extension of the constant field, consider a finite \'etale degree-2 map
of connected smooth affinoid strictly $k$-analytic curves $C' \rightarrow C$ such that for some
$c \in C$ there is a unique point $c' \in C'$ over $c$
and $\mathscr{H}(c')/\mathscr{H}(c)$ is a quadratic ramified extension.
Such examples are easily constructed from algebraic curves (even if $k$ is
algebraically closed or has trivial absolute value) by
taking $c$ to correspond to a suitable completion of the function field
of an algebraic curve.   Let $D = \calM(k \langle t \rangle)$, and let
$U$ be the gluing of $S_1 = C' \times D$ to $S_2 = C' \times D$
along $C' \times \{|t| = 1\}$ via the identity map.   This is a compact Hausdorff
$k$-analytic space of dimension 2, and it is not separated.  Let the nontrivial element in
$G = \Z/2\Z$ act on $U$ by swapping $S_1$ and $S_2$ and then applying
the nontrivial automorphism of $C'$ over $C$.  This is a free action.
As in the previous example, if the analogous rigid space situation admits
a quotient then so does the $k$-analytic situation, so to show that no quotient
exists in either case it suffices to prove that $U/G$ does not exist
as a $k$-analytic space.  Assume $U/G$ exists.
The natural map $U \rightarrow C$
is invariant by $G$, so it uniquely factors through a map $U/G \rightarrow C$.
The fiber $U_c$ equipped with its $G$-action
is an instance of Example \ref{nonsep} with $k'/k$ replaced
with $\mathscr{H}(c')/\mathscr{H}(c)$.
Moreover, the finite \'etale map $U_c \rightarrow (U/G)_c$ clearly
serves as a quotient $U_c/G$ in the category of
$\mathscr{H}(c)$-analytic spaces, contradicting that such a quotient
does not exist (by Example \ref{nonsep}).
\end{example}

\subsection{Local existence criteria}

We now turn to the task of proving Theorem \ref{ug}.
It has already been seen that the closed immersion hypothesis on the action map
$\alpha:U \times G \rightarrow U \times U$ forces $U$ to be separated, and conversely it is clear
 that if $U$ is separated then for any finite group equipped
 with a free right action on $U$ the action  map is a closed immersion.
To keep the role of the separatedness conditions clear, we will initially consider free right actions without
any assumption on the action map, and assuming only that $U$ is Hausdorff rather than separated. For any point
$u \in U$, with $G_u \subseteq G$ its isotropy group, say that the $G$-action on $U$ is {\em locally effective}
at $u$ if there is a $G_u$-stable analytic domain $N_u \subseteq U$ that is a neighborhood of $u$ such that
$N_u/G_u$ exists; obviously $G_u$ acts freely on $N_u$.

\begin{lemma}\label{locu}
Assume that $U$ is Hausdorff.
The quotient $U/G$ exists if and only if
the $G$-action is locally effective at all points of $U$.  If the quotient
$N_u/G_u$ exists and is good for each $u \in U$
and some $G_u$-stable $k$-analytic domain $N_u \subseteq U$ that is a neighborhood
of $u$ then $U/G$ is also good.  The same holds
for the property of being strictly $k$-analytic.
\end{lemma}

\begin{proof}
First assume that $U/G$ exists.
By Lemma \ref{finitelem}, the \'etale surjective quotient map $\pi:U \rightarrow U/G$ must be finite.
It is clearly a right $G$-torsor, so for each $\overline{u} \in U/G$ the fiber
$\pi^{-1}(\overline{u}) = \{u_1, \dots, u_n\}$
has transitive $G$-action with $\pi^{-1}(\overline{u}) = \coprod \calM(\mathscr{H}(u_i))$
over $\calM(\mathscr{H}(\overline{u}))$ for a finite Galois
extension $\mathscr{H}(u_i)$
of $\mathscr{H}(\overline{u})$ with Galois group $G_{u_i} \subseteq G$.
The equivalence of categories of finite \'etale covers of the germ $(U/G, \overline{u})$
and finite \'etale covers of $\Spec(\mathscr{H}(\overline{u}))$ \cite[3.4.1]{berihes}
thereby provides a connected open $\overline{V} \subseteq U/G$ around $\overline{u}$
such that $\pi^{-1}(\overline{V}) = \coprod V_i$, with
$V_i$ the connected component through $u_i$
and $V_i \rightarrow \overline{V}$ a finite \'etale $G_{u_i}$-torsor.
Hence, $V_i/G_{u_i}$ exists and equals $\overline{V}$.    Varying $\overline{u}$, this
implies that $G$ acts locally effectively on $U$.

Conversely, assume that the $G$-action is locally effective at every $u \in U$.
Let $N_u \subseteq U$ be a $G_u$-stable analytic domain that is a neighborhood of $u$
in $U$ and for which $N_u/G_u$ exists.
 The quotient map $\pi_u:N_u \rightarrow N_u/G_u$
is a finite \'etale $G_u$-torsor, so $u$ is the only point over $\overline{u} = \pi_u(u)$
and any $G_u$-stable open subset $W_u \subseteq N_u$
admits the open subset $\pi_u(W_u) \subseteq N_u/G_u$ as a quotient $W_u/G_u$.
In particular, any open set around $u$ in $N_u$ contains $\pi_u^{-1}(\overline{V})$
for some open $\overline{V} \subseteq N_u/G_u$ around $\overline{u}$.
Since $N_u$ is a neighborhood of $u$ in $U$, it follows that $u$ has a base
of $G_u$-stable open neighborhoods in $U$ that admit
a quotient by their $G_u$-action.   Since $U$ is Hausdorff and $G$ is finite, we can
therefore shrink the choice of $N_u$ so that $g(N_u) \cap N_u = \emptyset$
for all $g \in G$ with $g \not\in G_u$.    Hence,
the overlap $R_{N_u}$ of the analytic domain $R = U \times G \subseteq U \times U$
with $N_u \times N_u$ is $N_u \times G_u$.
Let $N_u^0 = {\rm{int}}_U(N_u)$, an open subset
of $U$ containing $u$.  This is $G_u$-stable
by functoriality of the topological interior with respect to automorphisms,
and the open overlap $R_{N_u^0} = R \cap (N_u^0 \times N_u^0)$ inside
of $R = U \times G$ is exactly $N_u^0 \times G_u$.  Thus, the quotient
$N_u^0/G_u$ serves as the quotient $N_u^0/R_{N_u^0}$.   We can therefore
use the locality argument with open sets following the statement
of Theorem \ref{newb} to glue the $N_u^0/G_u$'s to get a global quotient $U/G$
(which is good when every $N_u/G_u$ is good, and likewise
for strict $k$-analyticity, since $N_u^0/G_u$ is an open subspace of $N_u/G_u$ for
every $u \in U$).
\end{proof}

By Lemma \ref{locu}, the global existence of $U/G$ for a free right action of a finite group $G$ on a Hausdorff
$k$-analytic space $U$ is reduced to the local effectivity of the $G$-action (and keeping track of the
properties of goodness and strict $k$-analyticity for the local quotients by isotropy groups).   We now use this
to give a criterion for the existence of $U/G$ in terms of good analytic domains when $U$ is Hausdorff.

\begin{lemma}\label{211} With notation and hypotheses as above,
$G$ is locally effective at $u \in U$ if and only if there are finitely many
$k$-affinoid domains $V \subseteq U$ containing $u$ such that
the germs $(V,u) \subseteq (U,u)$ cover $(U,u)$ and are $G_u$-stable as germs.
It is equivalent to take the $V$'s to merely be good analytic domains in $U$ with $u \in V$.
In such cases, if $U$ has a strictly $k$-analytic $($resp. good$)$ open
neighborhood of $u$ then $N_u/G_u$ exists and is strictly $k$-analytic $($resp. good$)$
for some $G_u$-stable $k$-analytic domain $N_u$ in $U$ containing $u$ that we may
take to be compact and as small as we please.
\end{lemma}

Note that $G_u$-stability of a germ $(V,u)$ is a weaker condition than $G_u$-stability of
the $k$-analytic domain $V$ in $U$.

\begin{proof}  First assume that $G$ is locally effective at $u$, so
$N_u/G_u$ exists for some $G_u$-stable $k$-analytic domain $N_u \subseteq U$
that is a neighborhood of $u$ in $U$.   We can replace $U$ and $G$ with $N_u$ and $G_u$
so that $G = G_u$ and $U/G$ exists.
Let $\pi:U \rightarrow U/G$ be the quotient
map and let $\overline{u} = \pi(u)$.  Let $\overline{V}_1, \dots, \overline{V}_n$
be $k$-affinoid domains in $U/G$ containing $\overline{u}$ such that
$\cup \overline{V}_i$ is a neighborhood of $\overline{u}$ in $U/G$.
Each $V_i = \pi^{-1}(\overline{V}_i)$
is a $k$-affinoid domain (since $\pi$ is finite), contains $u$, and is $G$-stable.
The finite collection of germs $(V_i,u)$ therefore satisfies
the required conditions.

Conversely, assume that there are good $k$-analytic
domains $V_1, \dots, V_n$  in $U$ containing
$u$ such that $(V_i,u) \subseteq (U,u)$ is $G_u$-stable as a germ and
such that these germs cover $(U,u)$.  We may and do shrink
each good domain $V_i$ around $u$ so that it is a separated
(e.g., $k$-affinoid) domain.  Thus, $\cup V_i$ is a neighborhood
of $u$ in $U$ and $g V_i \cap V_i$ is a neighborhood of $u$ in both $V_i$ and $gV_i$
for all $g \in G_u$.   In particular, $V'_i = \cap_{g \in G_u} g V_i$ is
a $G_u$-stable $k$-analytic domain with $(V'_i,u) = (V_i,u)$ as domains in the germ
$(U,u)$.   We may therefore assume that each $V'_i$ is separated, good, and $G_u$-stable.
Let $U_i \subseteq V_i$ be a $k$-affinoid neighborhood of $u$ in
the good $k$-analytic space $V_i$.  The overlap
$U'_i = \cap_{g \in G_u} gU_i$ is a $k$-affinoid neighborhood of $u$ in $V_i$ (hence in $U$)
since $V_i$ is separated.   But $(U_i,u) = (V_i,u)$ is a $G_u$-stable germ domain
in $(U,u)$, so $(U'_i,u) = (V_i,u)$ as well.  Thus, we can assume that
each $V_i = \calM(\calA_i)$ is $k$-affinoid and $G_u$-stable.
By \cite[2.1.14({\em{ii}})]{berbook},
$\calA_i^{G_u}$ is a $k$-affinoid algebra over which $\calA_i$ is finite and admissible,
so $\calM(\calA_i^{G_u})$ makes sense and $V_i \rightarrow
\calM(\calA_i^{G_u})$ is a finite map.   Moreover, if
$\calA_i$ is strictly $k$-analytic then so is $\calA_i^{G_u}$.  The action of $G_u$ on
$\Spec(\calA_i)$ is free because $V_i = \calM(\calA_i) \rightarrow
\Spec(\calA_i)$ is surjective and $G_u$ acts freely on $V_i$.
Hence, by \cite[Exp.~V,~4.1]{sga3}, the finite admissible map $\calA_i^{G_u} \rightarrow \calA_i$
is \'etale and the natural map $\calA_i \otimes_{\calA_i^{G_u}} \calA_i
\rightarrow \prod_{g \in G_u} \calA_i$ defined by $a \otimes a' \mapsto (a g(a'))_g$
is an isomorphism.  This latter tensor product coincides
with the corresponding completed tensor product, so
it exhibits $\calM(\calA_i^{G_u})$ as the quotient $V_i/G_u$ for all $i$.

Let $N_u = \cup V_i$, a compact $G_u$-stable $k$-analytic domain in $U$
that is a neighborhood of $u$.  Observe that by our initial shrinking of
the $V_i$'s we can take $N_u$ to be as small as we please.  If $U$ is strictly $k$-analytic then it
is clear that we could have taken each $V_i$ to be strictly $k$-analytic,
and every overlap $V_i \cap V_j = V_i \times_U V_j$ is then strictly
$k$-analytic.   If $U$ is good then we could have shrunk $U$ around $u$ to be
separated, so there is a $G_u$-stable $k$-affinoid neighborhood $V \subseteq U$ around $u$,
and hence $N_u = \cap_{g \in G_u} gV$ is a $G_u$-stable $k$-affinoid neighborhood
of $u$ in $U$.
It suffices to prove that $N_u/G_u$ exists
(and that it is strictly $k$-analytic, resp. $k$-affinoid, when $N_u$ is).
Thus, we can replace $U$ and $G$ with $N_u$ and $G_u$.
That is, we may assume that $U$ is covered by $G$-stable $k$-affinoid domains
$V_1, \dots, V_n$, and that $V_i/G$ exists as a $k$-affinoid domain
(even strictly $k$-analytic when $V_i$ is so).
Our aim is to construct $U/G$ in this case, with $U/G$ strictly $k$-analytic
when $U$ and all $V_i$ are strictly $k$-analytic, and with $U/G$ a
(strictly) $k$-affinoid space
when $U$ is (strictly) $k$-affinoid.   This latter (strictly) $k$-affinoid case
was already settled in the preceding paragraph. In general,
since $U$ is Hausdorff, each overlap $V_{ij} = V_i \cap V_j$ is
a $G$-stable compact $k$-analytic domain in $U$.  Moreover, if
$U$ and all $V_i$ are strictly $k$-analytic then so is each $V_{ij}$.
Applying Example \ref{ugdomain} to the finite \'etale quotient maps
$\pi_i:V_i \rightarrow V_i/G$ and $\pi_j:V_j \rightarrow V_j/G$,
$\pi_i(V_{ij})$ and $\pi_j(V_{ij})$ are each identified
as a quotient $V_{ij}/G$ and each is strictly
$k$-analytic when $U$, $V_i$, and $V_j$ are strictly
$k$-analytic.  Thus, $\pi_i(V_{ij})$ and $\pi_j(V_{ij}) = \pi_j(V_{ji})$ are uniquely isomorphic
(say via an isomorphism $\phi_{ji}:\pi_i(V_{ij}) \simeq \pi_j(V_{ji})$)
in a manner that respects the maps from $V_{ij}$ onto each, and
the analogous quotient conclusions hold for the triple
overlaps among the $V_i$'s.

It is easy to check (by chasing $G$-actions)
that the triple overlap cocycle condition holds.
Thus,
we can define a $k$-analytic space $\overline{U}$ that is covered
by the $k$-analytic domains $V_i/G$ with overlaps $V_{ij}/G$ (using the $\phi_{ji}$'s),
and this gluing $\overline{U}$ is strictly $k$-analytic when $U$ and all $V_i$
are strictly $k$-analytic.
By computing on geometric points
we see that the maps $V_i \rightarrow V_i/G_i \rightarrow \overline{U}$
uniquely glue to a morphism $\pi:U \rightarrow \overline{U}$ whose pullback over
each $k$-analytic domain $V_i/G \subseteq \overline{U}$ is the $k$-analytic domain
$V_i \subseteq U$.   Hence, $\pi$  is finite \'etale since each $\pi_i$ is finite \'etale,
and $\pi$ is $G$-invariant since every $\pi_i$ is $G$-invariant.
The resulting natural map $U \times G \rightarrow U \times_{\overline{U}} U$
over $U \times U$ restricts to the natural map $V_{ji} \times G
\rightarrow V_{ji} \times_{V_{ji}/G, \phi_{ji} \circ \pi_i} V_{ij}$ over $V_j \times V_i$ for all $i$ and $j$,
and these latter maps are isomorphisms.  Hence, $\pi:U \rightarrow \overline{U}$
serves as a quotient $U/G$.  By construction, if $U$ is strictly $k$-analytic then
so is $\overline{U}$.
\end{proof}

By Lemma \ref{211} and the gluing used in the proof of Lemma \ref{locu}, it remains to show that if $U$ is
separated then every germ $(U,u)$ is covered by finitely many $G_u$-stable good subdomain germs $(V,u) \subseteq
(U,u)$ with $V \subseteq U$ a $k$-analytic domain.  In particular, the finer claims concerning inheritance of
goodness and strict $k$-analyticity of $U/G$ are immediate corollaries of such an existence result in the
general (separated) case.   As above, we postpone the separatedness hypothesis on $U$ until we need it, assuming at the outset
just that $U$ is Hausdorff. We fix $u \in U$, and by Lemma \ref{locu} we may rename $G_u$ as $G$, so we
can assume $G = G_u$. To construct the required collection of subdomains around $u$ we will use the theory of
reduction of germs developed in \cite{temkin2}.    For ease of notation, we now write $U_u$ rather than $(U,u)$,
and we write $\widetilde{U}_u$ to denote the reduction of this germ (an object in the category
${\rm{bir}}_{\widetilde{k}}$ of birational spaces over the $\R_{>0}^{\times}$-graded field $\widetilde{k}$; see
\cite[\S1--\S3]{temkin2} for terminology related to graded fields and birational spaces over them,
such as the quasi-compact Zariski--Riemann space $\mathbf{P}_{L/\ell}$ of graded valuation rings
associated to an extension $\ell \rightarrow L$ of graded fields). The
reduction of germs is a useful technique for studying the local structure of a $k$-analytic space.  For example,
in \cite[4.5]{temkin2} it is shown that $V_u \mapsto \widetilde{V}_u$ is a bijection from the set of subdomains
of $U_u$ to the set of quasi-compact open subspaces of the birational space $\widetilde{U}_u$, and by
\cite[4.1({\em{ii}})]{temkin2} this bijection is compatible with inclusions and finite unions and intersections.
Moreover, by \cite[4.8({\em{iii}})]{temkin2} (resp. \cite[5.1]{temkin2}) the germ $V_u$ is separated (resp.
good) if and only if the reduction $\widetilde{V}_u$ is a separated birational space (resp. an affine birational
space) over $\widetilde{k}$.

At this point we assume that $U$ is separated.
To find a finite collection of $G$-stable good subdomains $(V_j)_u$ in $U_u$ that cover $U_u$, it is equivalent
to cover the separated birational space $\widetilde{U}_u$ over $\widetilde{k}$ by $G$-stable affine open
subspaces. Note that any intersection of finitely many affine subspaces of $\tilU_u$ is affine (because
separatedness of $\tilU_u$ allows us to identify $\tilU_u$ with an open subspace of the birational
space
$\bfP_{\widetilde{K}/\widetilde{k}}$). Choose a point $z\in\tilU_u$ and let $Gz$ denote its $G$-orbit. If $Gz$
is contained in an affine open subspace $V\subseteq\tilU_u$ then the intersection of all $G$-translates of $V$
is a $G$-stable affine neighborhood of $z$. Thus, by quasi-compactness of $\tilU_u$ our problem reduces to
proving the following statement.

\begin{theorem}\label{213}
Let $H$ be a commutative group and
let $k$ be an $H$-graded field and $U$ a separated object of $\bir_k$ provided with an action of a finite group
$G$. Any $G$-orbit $S\subseteq U$ admits an affine open neighborhood.
\end{theorem}

In the intended application of this theorem we have $H = \mathbf{R}^{\times}_{>0}$, so the group law on $H$ will
be denoted multiplicatively below. By the separatedness hypothesis, we can identify $U$ with an open subspace of
a Zariski--Riemann space of $H$-graded valuations $\bfP_K=\bfP_{K/k}$, where $K/k$ is an extension of $H$-graded
fields. In the classical situation when $K/k$ is finitely generated and gradings are trivial (i.e., $H = \{1\}$)
one can use Chow's lemma to prove the generalization of Theorem \ref{213} in which there is no
$G$-action and $S\subseteq U$ is an arbitrary finite subset of $\bfP_{K}$
: (i) find a model $\Spec(K)\to X$ in ${\rm{Var}}_k$ (see \cite[\S 1]{temkin1}), (ii) using Chow's lemma,
replace $X$ by a modification that is quasi-projective over $k$, and (iii) use that any finite set in a
quasi-projective $k$-scheme is contained in an open affine subscheme. Surprisingly, one has to be very careful
in the general graded case. For example, Theorem \ref{213} without the $G$-action is false for an arbitrary
finite $S$.

\subsection{Graded valuations on graded fields}

To prove Theorem \ref{213}, we first
need to generalize some classical results to the graded case. Throughout this section,
we assume all gradings are taken with respect to a fixed commutative group $H$,
and ``graded'' always means ``$H$-graded''.  We start with a certain portion of Galois
theory that is very similar to the theory of tamely ramified extensions of valued fields. Let $K/L$ be an
extension of graded fields.
Since $K$ is a free $L$-module by \cite[2.1]{temkin2}, we can define the {\em
extension degree} $n = [K:L]$ to be the $L$-rank of $K$;
we say that $K/L$ is a {\em finite} extension when the extension
degree is finite. Two more invariants of the extension $K/L$ are
analogues of the residual degree and ramification index in the classical theory,
defined as follows.   Writing $\oplus_{h \in H} K_h$ and $\oplus_{h \in H} L_h$ for
the decompositions of $K$ and $L$ into their graded components,
the components $K_1$ and $L_1$ indexed by the identity of $H$ are ordinary fields and
we define $f = f_{K/L} = [K_1:L_1]$.  For any nonzero graded ring $A$ we let
$A^{\times}$ denote the homogeneous
unit group of $A$ (i.e., the ordinary units of $A$ that are homogeneous with
respect to the grading); note that $A_1 \cap A^{\times}$ is the unit group of
the ring $A_1$.
We let $\rho:A^{\times} \rightarrow H$ denote the multiplicative map that
sends each homogeneous unit to its uniquely determined
grading index, so $\rho(A^{\times}) \subseteq H$ is a subgroup.
Define $e=e_{K/L}=\#(\rho(K^{\times})/\rho(L^{\times}))$.
The invariants $n$, $e$, and $f$ of $K/L$ may be infinite.

\begin{lemma}
The equality $n=ef$ holds, where we use the conventions
$\infty \cdot d = \infty \cdot \infty = \infty$ for any $d \ge 1$.

\end{lemma}
\begin{proof}
Let $B$ be a basis of $K_1$ over $L_1$ and $T \subseteq K^{\times}$ be any
set of representatives for the fibers of the surjective homomorphism
$K^{\times}\to\rho(K^{\times})\to\rho(K^{\times})/\rho(L^{\times})$.
Since $K$ is a graded field, in the decomposition $K = \oplus_{r \in \rho(K^{\times})} K_r$
the $K_r$'s are all 1-dimensional as $K_1$-vector spaces.  The same holds for $L$,
so the
products $bt$ for $(b,t) \in B \times T$ are pairwise
distinct and the set of these products is a basis of $K$ over $L$.
The lemma now follows.
\end{proof}

For any extension $K/L$ of graded fields, the set of intermediate graded fields
$F$ satisfying $e_{F/L} = 1$ is in natural bijection with the set of intermediate
fields $F_1$ in $K_1/L_1$ via the recipe
$F_1 \mapsto F_1 \otimes_{L_1} L$, so $F = K_1 \otimes_{L_1} L$ is
such an intermediate graded field and it contains all others.   It is clear
that this particular $F$ is the unique intermediate graded field in $K/L$ satisfying
$e_{F/L} = 1$ and $f_{K/F} = 1$.
For a fixed $K$ the graded subfields $L \subseteq K$ such that
$f_{K/L}=1$ are in natural bijection with subgroups of $\rho(K^{\times})$
via the recipe $L=\oplus_{r\in\rho(L^{\times})}K_r$, and in such
cases $K/L$ is finite if and only if $e_{K/L}$ is finite.
We say that an extension $K/L$ is
{\em totally ramified} if it is a finite extension and $f_{K/L} = 1$.

Let $K$ be a graded field and $G \subseteq \Aut(K)$ a finite subgroup.
The graded subring $L = K^G$ is obviously a graded field, and
$L_1 = K_1^G$.  The {\em inertia subgroup} $I \subseteq G$ is
the subgroup of elements that act trivially on $K_1$, so
$G/I \subseteq \Aut(K_1)$.
An element $\sigma\in I$ acts $K_1$-linearly on the 1-dimensional
$K_1$-vector space $K_r$ for each $r \in \rho(K^{\times})$, so this
action must be multiplication by an
element $\xi_{\sigma,r}\in K_1^{\times}$
that is obviously a root of unity.  In particular, $I$ is always abelian,
just like for tamely ramified extensions in classical valuation theory.
By a lemma of Artin, $K_1/L_1$ is a finite
Galois extension with $G/I \simeq {\rm{Gal}}(K_1/L_1)$.
In particular,
$[K_1:L_1]|\#G$.  Here is an analogue of the lemma of Artin for $K/L$.

\begin{lemma}\label{artin}
Let $K$ be a graded field and $G\subseteq\Aut(K)$ a finite subgroup. The
degree of $K$ over the graded
subfield $L=K^G$ is finite and equal to $\#G$.
\end{lemma}

\begin{proof}
Let $I$ be the inertia subgroup of $G$.
The graded subfield $F=K^I$ has $F_1 = K_1$ (i.e., $f_{K/F} = 1$), so $F$ is
uniquely determined by its value group $\rho(F^{\times}) \subseteq \rho(K^{\times})$.
Clearly $\rho(F^{\times})$ is the set of $r$'s such
that $I$ acts trivially on $K_r$, and this is a subgroup $V\subseteq\rho(K^{\times})$.
Since $f_{K/F} = 1$, we have $[K:F] = e_{K/F} = \#(\rho(K^{\times})/V)$.
For any $r \in \rho(K^{\times})$, the ratio $\sigma(x)/x \in K_1^{\times}$
for nonzero $x \in K_r$ is independent of $x$ and only depends on $r \bmod V$.
Thus, we can denote this ratio $\xi_{\sigma, \overline{r}}$ for $\overline{r} = r \bmod V$.
The pairing
$\xi:I\times(\rho(K^{\times})/V)\rightarrow  K_1^{\times}$ defined by
$(\sigma, \overline{r}) \mapsto \xi_{\sigma,\overline{r}}$ is biadditive and nondegenerate
in each variable.  In particular, $\rho(K^{\times})/V$ is finite and
$\xi$ takes values in a finite subgroup of $K_1^{\times}$ that
is necessarily cyclic with order that annihilates $I$ and $\rho(K^{\times})/V$.
Thus, by duality for finite abelian groups we get $\rho(K^{\times})/V = \#I$, so $[K:F]$ is finite
and equal to $\#I$.  That is,
the result holds with $I$ in the role of $G$.

The group $G'=G/I$ naturally acts on $F$.  Since $F^{G'}=L$, it suffices to prove
$[F:L] = \#G'$.
Viewing $F$ as a vector space over $F_1 = K_1$, its $G'$-action is
semilinear.  Artin's lemma in the classical case gives that
$F_1$ is a finite Galois extension of $L_1$ with Galois group $G'$.
Hence, Galois descent for
vector spaces provides an identification $F = F_1 \otimes_{L_1} L$,
so $[F:L] = [F_1:L_1] = \#G'$ (and $e_{F/L} = 1$).
\end{proof}

We now  study extensions of graded valuation rings. Let $L$ be a graded field,
and let $R$ be a graded valuation ring of $L$, meaning that $R \subseteq L$
is a graded subring that is a graded valuation ring with
$L$ equal to its graded fraction field.    Let $K/L$ be an extension of graded
fields.  For any graded valuation ring $A$ of $K$, $A_1$ is a valuation
ring of $K_1$
and $A \cap L$  is a graded valuation ring of $L$ whose inclusion into $A$ is a graded-local map.
If $A \cap L = R$ then we say that $A$ {\em extends} $R$ (with respect to $K/L$).
Note that in such cases, the valuation ring $A_1$ of the field $K_1$ extends the valuation
ring $R_1$ of the field $L_1$, and $R_1 \rightarrow A_1$ is a local map
(i.e., $A_1^{\times} \cap L_1 = R_1^{\times}$).

\begin{lemma}\label{extendv}
Let $K/L$ be a finite extension
of graded fields and fix a graded valuation ring $R$ of $L$. The correspondence
$A\mapsto A_1$ is a bijection between
the set of extensions of $R$ to a graded valuation ring of $K$ and the set of
extensions of $R_1$ to a valuation ring of $K_1$.  Moreover, if
$R \subseteq R'$ is an inclusion of graded valuation rings of $L$,
and $A$ and $A'$ are graded valuation rings of $K$ that respectively
extend $R$ and $R'$, then $A \subseteq A'$ if and only if
$A_1 \subseteq A'_1$.
\end{lemma}

\begin{proof}
The intermediate graded field $F = K_1 \otimes_{L_1} L$ satisfies
$e_{F/L} = 1$ and $f_{K/F} = 1$,
so it suffices to consider separately the cases $e = 1$ and $f = 1$.
First assume $f = 1$, so we have to show that $R$ admits a unique
extension to $K$ (and that inclusions among such extensions over
an inclusion $R \subseteq R'$ in $L$ can
be detected in $K_1$).  An extension of $R$ to $K$ exists without any
hypotheses on the graded field extension $K/L$ (by
Zorn's Lemma), so the main issue is to prove uniqueness.   Due to the grading,
we just have to check that the homogeneous elements in such an $A \subseteq K$
are uniquely determined.
Since $L_1 = K_1$, we have
$L = \oplus_{r \in \rho(L^{\times})} K_r$
with $\rho(L^{\times})$ an index-$e$ subgroup of
$\rho(K^{\times})$.  Hence,
any homogeneous element $a\in K$ is an $e$th root of a homogeneous element
of $L$. If $a \in K^{\times}$ satisfies $a \not\in A$
then $1/a$ lies in the unique graded-maximal ideal of $A$, so
$(1/a)^e \in A \cap L = R$ lies in the unique graded-maximal ideal of $R$
(since $R \rightarrow A$ is a graded-local map).  Hence, $a^e \not\in R$ in such cases,
so we have the characterization that a homogeneous $a \in K$ lies in $A$ if
and only if $a^e \in R$.   This gives the desired uniqueness, and also
shows that if $R \subseteq R'$ is a containment of graded valuation rings of $L$
then their unique respective extensions $A, A' \subseteq K$ satisfy $A \subseteq A'$;
observe also that in this case $A_1 = R_1$ and $A'_1 = R'_1$.

It remains to analyze the case $e=1$, so
$K= K_1 \otimes_{L_1} L$.  In particular, $K_r = K_1 \cdot L_r$ for all
$r \in \rho(K^{\times}) = \rho(L^{\times})$.
We have to prove that for any extension $A_1\subseteq K_1$ of $R_1 \subseteq L_1$,
there exists a unique graded valuation ring $B\subseteq K$ that extends $R$ and satsifies $B_1=A_1$,
and that if $A'_1$ is an extension of $R'_1$ where
$R' := R'_1 \cap L$ is a graded valuation ring of $L$ containing $R$ then the corresponding
$B'$ in $K$ satisfies $B \subseteq B'$ if and only if $A_1 \subseteq A'_1$.
As a preliminary step, we recall a fact from classical valuation theory:
if $k'/k$ is a degree-$d$ extension of fields and $V'$ is a valuation ring of $k'$ then
$(x')^{d!} \in k \cdot {V'}^{\times}$ for all $x' \in k'$.
This follows from the fact that all ramification indices are finite and bounded above by
the field degree.
In our initial situation, we conclude that
for $N = f!$ and any homogeneous element $x\in K_r = K_1 \cdot L_r$, we have $$x^N \in (K_1 \cdot
L_r)^N\subseteq K_1^N \cdot L_{r^{N}}
\subseteq (L_1 \cdot A_1^{\times}) \cdot L_{r^{N}}=A_1^{\times} \cdot L_{r^{N}}.$$
Given such a factorization $x^N=a_1 \ell$ with $a_1 \in A_1^{\times}$
and (necessarily)
homogeneous $\ell\in L$, the element $\ell$ is uniquely determined up to a factor lying in
$L\cap A_1^{\times}=R_1^{\times}$. In particular, it is
an intrinsic property of $x$ whether or not $\ell \in R$, so if $B$ exists then it
is contained in the set $A$ of
$x\in K$ each of whose homogeneous components $x_h \in K_h$ lies in the subset
$$A_h := \{y \in K_h\,|\,
y^N\in A_1^{\times} \cdot R_{h^N}\}.$$
We will show that this set $A$ is a graded valuation ring of $K$, so
its homogeneous units are those homogeneous $x \in K^{\times}$
such that $x^N \in A_1^{\times} \cdot R^{\times}$, and
moreover $A \cap L = R$ and $A \cap K_1 = A_1$.  Hence,
this $A$ works as such a $B$, and any valuation ring $B$
of $K$ that extends $R$ and satisfies $B \cap K_1 = A_1$ must have graded-local
inclusion into $A$, thereby
forcing $B = A$.   This would give the required uniqueness, and also implies
the desired criterion for containment of extensions (over a containment
$R \subseteq R'$ in $L$) by checking in $K_1$.

Since it is clear that $A_h \cdot A_{h'} \subseteq A_{hh'}$ for all
$h, h' \in H$, to show that $A$ is a graded subring of $K$ we just need to prove
that if $x,y$ are two nonzero elements of $A_h=A\cap K_h$
for some $h \in H$ then $x+y\in A_h$.
The ratios $x/y$ and $y/x$ in $K$ both lie in $K_1$ and are inverse to each
other, so at least one of them lies in the valuation ring $A_1$ of  $K_1$.  Thus, by switching
the roles of $x$ and $y$ if necessary we may assume $x/y \in A_1$.
We then have $(x+y)/y = 1 + x/y \in A_1$, so
$((x+y)/y)^N \in A_1^{\times} \cdot R_1$ by the classical valuation-theoretic
fact recalled above.  Hence,
$$(x+y)^N = \left(\frac{x+y}{y}\right)^N \cdot y^N \in (A_1^{\times} \cdot R_1) \cdot
(A_1^{\times} \cdot R_{h^{N}}) = A_1^{\times} \cdot R_{h^{N}},$$
so $x+y \in A$ as desired.   Finally, to show that $A$
is a graded valuation ring of $K$ it suffices to prove that for any homogeneous
nonzero $t\in K$ either $t$ or $1/t$ lies in $A$.
We have $t^N = a_1 \ell$ with $a_1 \in A_1^{\times}$ and a homogeneous nonzero $\ell \in L$.
But $L$ is the graded fraction field of the graded valuation ring $R$, so
$\ell \in R$ or $1/\ell \in R$, and hence $t \in A$ or $1/t \in A$ respectively.
\end{proof}

\begin{corollary}\label{trans}
Let $K$ be a graded field, $G\subseteq\Aut(K)$ a finite subgroup, and $R$ a graded valuation ring of $L=K^G$. Then

$($i$)$ $G$ acts transitively on the non-empty set of extensions of $R$ to $K$;

$($ii$)$ if $R'$ is a graded valuation ring of $L$ containing $R$
as a graded subring and $A'$ is an extension of $R'$ to $K$, then there
exists an extension $A$ of $R$ to $K$ with $A \subseteq A'$ as graded subrings of $K$.
\end{corollary}

\begin{proof}
Let $I \subseteq G$ be the inertia subgroup,
so $G/I = {\rm{Gal}}(K_1/L_1)$.  By classical valuation theory,
$G/I$ acts transitively on the set of extensions of $R_1$ to
$K_1$, so (i) follows from Lemma \ref{extendv}.  For (ii),
classical valuation theory also
gives that $R_1$ admits an extension $A_1$ to $K_1$ with $A_1 \subseteq A'_1$,
so consider the corresponding $A$ extending $R$ to $K$.
To see $A \subseteq A'$ we can use the containment criterion in
Lemma \ref{extendv}.
\end{proof}

\begin{corollary}
\label{transitcor} If $K/k$ is an extension of graded fields, $G\subseteq\Aut_k(K)$
is a finite subgroup, and $L=K^G$,
then

$($i$)$ $G$ acts transitively on the fibers of the natural surjective map $\psi:\bfP_{K/k}\to\bfP_{L/k}$;

$($ii$)$
if $x\in\bfP_{L/k}$ is a point, $S=\psi^{-1}(x)$ is the fiber, and $\ox$ is the set of generalizations of $x$
in $\bfP_{L/k}$
then $\psi^{-1}(\ox)$ is the set of generalizations of the points of $S$.
\end{corollary}
\begin{proof}
The first part is Corollary \ref{trans}(i). To prove the second part we note that
a point $y\in\bfP_{L/k}$ is a generalization of $x$
(i.e. is contained in every
open neighborhood of $x$) if and only if its associated graded valuation ring of $L$ contains the
one associated to
$x$. Hence (ii) follows from Corollary \ref{trans}(ii).
\end{proof}

Another classical notion whose graded analogue will be used in the proof of Theorem \ref{213} is the
constructible
topology on a Zariski--Riemann space (also sometimes called
the patching topology). Let $K/k$ be an extension of graded fields,
and let $\bfP_{K/k}$ be the associated Zariski--Riemann space. Consider subsets
$$\bfP_{K/k}\{F\}\overline{\{G\}}=\{\calO\in\bfP_{K/k}\,|\,F\subseteq \calO, G \cap \calO = \emptyset\}
\subseteq \bfP_{K/k}$$
for $F, G \subseteq K^{\times}$.  Such subsets with finite $F$ and empty $G$
form a basis of the usual topology on $\bfP_{K/k}$.
Since $\cap_i \bfP_{K/k}\{F_i\}\overline{\{G_i\}} =
\bfP_{K/k}\{\cup F_i\}\overline{\{\cup G_i\}}$,
the subsets $\bfP_{K/k}\{F\}\overline{\{G\}}$ with finite $F$ and $G$
satisfy the requirements to be a basis of open sets for a finer
topology on $\bfP_{K/k}$ called the {\em constructible topology}.

\begin{lemma}\label{patchtop}
Let $F, G \subseteq K^{\times}$ be arbitrary subsets. The subset $X =
\bfP_{K/k}\{F\}\overline{\{G\}} \subseteq \bfP_{K/k}$ is compact and Hausdorff with respect to the
constructible topology on
$\bfP_{K/k}$.
\end{lemma}

\begin{proof}
Let $S_{K}$ be the set of all subsets of $K^{\times}$.
Any graded valuation ring of $K$ is uniquely determined by its intersection
with $K^{\times}$, so there is a natural injection
$i:\bfP_{K/k} \rightarrow S_{K}$ given by $i(\calO) = \calO \cap K^{\times}$.
Consider the bijection $S_{K} \rightarrow \{0,1\}^{K^{\times}}$
that assigns to each $\Sigma \in S_{K}$ the
characteristic function of $\Sigma$.
The discrete topology on $\{0,1\}$ thereby endows $S_{K}$ with a compact
Hausdorff product topology, and a basis of open sets for this topology is given by the sets
$$S_{K}\{F'\}\overline{\{G'\}} = \{\Sigma \in S_{K}\,|\,F' \subseteq \Sigma, G' \cap \Sigma =
\emptyset\}$$
for finite subsets $F', G' \subseteq K^{\times}$.
The sets $S_K\{F'\}\overline{\{\emptyset\}}$ and $S_K\{\emptyset\}\overline{\{G'\}}$
will be denoted $S_K\{F'\}$ and $S_K\overline{\{G'\}}$ respectively.
Clearly the subspace
topology induced on $\bfP_{K/k}$ via $i$ is the constructible topology, so it suffices
to prove that $i(X)$ is closed in $S_K$.
 We will show that $S_K - i(X)$ is open
in $S_K$ by covering it by members of the above basis of open sets in $S_K$.

Choose a point $\Sigma \in S_K$.  The condition that $\Sigma \not\in i(X)$ means
that the subset $\Sigma \subseteq K^{\times}$ cannot be expressed as $\calO \cap K^{\times}$
for a graded valuation ring  $\calO$ of $K$ containing $k$.  The only possibility for such
an $\calO$ is the graded additive subgroup $\calO_{\Sigma} \subseteq K$ generated by
$\Sigma$.  Thus,
$\Sigma \not\in i(\bfP_{K/k})$ if and only if either $\calO_{\Sigma} \cap K^{\times}$
is strictly larger than $\Sigma$ or $\calO_{\Sigma} \cap K^{\times} =
\Sigma$ but $\calO$ fails to be a graded valuation ring
of $K$ containing $k$.  If $\Sigma \in i(\bfP_{K/k})$
(i.e.,  $\calO_{\Sigma}$ is a graded valuation ring
of $K$ containing $k$ with $\calO_{\Sigma} \cap K^{\times} = \Sigma$)
then the finer condition $\Sigma \not\in i(X)$ says exactly that $F \not\subseteq \calO_{\Sigma}$
or $G \cap \calO_{\Sigma} \ne \emptyset$.  In other words,
$\Sigma \not\in i(X)$ if and only if $\Sigma$ satisfies at least one of the following seven properties:
(i) $a+b \not\in\Sigma$ for some $a, b \in \Sigma$ lying in the same
graded component of $K$ with $a + b \in K^{\times}$ (i.e., with $a+b \ne 0$), (i$'$) $-a \not\in \Sigma$
for some $a \in \Sigma$, (ii) $a \not\in \Sigma$ for some $a \in k^{\times}$,
(iii) $ab \not\in \Sigma$ for some $a, b \in \Sigma$, (iv) $a \not\in \Sigma$
and $1/a \not\in \Sigma$ for some $a \in K^{\times}$, (v) $f \not\in \Sigma$
for some $f \in F$, (v$'$) $g \in \Sigma$ for some $g \in G$.
Hence, it suffices to show that if $\Sigma \subseteq K^{\times}$
is a subset satisfying one of these conditions then it has an open
neighborhood in $S_K$ satisfying the same condition.  For each respective condition,
use the following neighborhood (with notation as above):
(i) $S_K\{a,b\}\overline{\{a+b\}}$, (i$'$) $S_K\{a\}\overline{\{-a\}}$,
(ii) $S_K\overline{\{a\}}$, (iii) $S_K\{a,b\}\overline{\{ab\}}$, (iv) $S_K\overline{\{a, 1/a\}}$,
(v) $S_K\overline{\{f\}}$, (v$'$) $S_K\{g\}$.
\end{proof}

Now we have all of the necessary tools to prove Theorem \ref{213}.

\begin{proof}
By definition, the data of the separated birational space $U$ consists of the specification
of a connected quasi-compact and
quasi-separated topological space equipped with an open embedding into
$\bfP_{K} = \bfP_{K/k}$.  Thus, although (by abuse of notation) we shall write $U$ to
denote the open subspace, the action of $G$ on $U$ as a birational space really
means an action $\alpha$ of $G$ on both the open subspace and on the graded field
$K$ over $k$ such that the induced action on $\bfP_K$ carries the open subspace
back to itself via $\alpha$.  That is, the action by $G$ on $\bfP_{K}$ arising from the $G$-action on
$K$ over $k$ induces the original action of $G$ on $U\subseteq\bfP_K$. Let $L=K^G$, so
$L/k$ is a graded field extension as well.  By
Corollary \ref{transitcor}(i), $S$ is a fiber of the induced map $\psi:\bfP_K\to\bfP_L$ over a point
$x\in\bfP_L$.   Our problem is therefore to find a $G$-stable affine open neighborhood of
$S$ in $\bfP_{K}$ that contains $\psi^{-1}(x)$ and is contained in $U$.

Let $\{U_i\}_{i\in I}$ be the family of all affine open neighborhoods of $x$, so
$\{V_i=\psi^{-1}(U_i)\}_{i\in I}$ is a family of $G$-stable affine
open neighborhoods of $S$. Clearly $\ox:=\cap_{i\in
I}U_i$ is the set of all generalizations of $x$, and by  Corollary \ref{transitcor}(ii) we see that
$\oS:=\cap_{i\in I}V_i$ is the set of generalizations of points of $S$.  In particular, $\oS\subseteq U$.
Note that $U$ is open and the
$V_i$'s are closed in the constructible topology on $\bfP_K$, so $U$ is an open
neighborhood of the intersection of closed subsets $V_i$ in the compact Hausdorff space $\bfP_K$ (provided with
the constructible topology). It follows that $U$ contains the intersection of finitely many $V_i$'s,
and that intersection is the required $G$-stable affine open neighborhood of $S$.
\end{proof}

\end{document}